\documentclass[a4paper,11pt]{article}
\usepackage{amsthm}
\usepackage{tikz}
\usetikzlibrary{decorations.markings,decorations.pathreplacing,positioning,cd,fit,shapes.geometric}
\tikzset{radial vector/.style={blue, -latex}}
\tikzset{radial vector 2/.style={red, -latex}}
\tikzset{midarrow/.style={postaction=decorate,decoration={markings,mark={at position #1 with {\arrow{latex}}}}}}
\tikzset{invmidarrow/.style={postaction=decorate,decoration={markings,mark={between positions 0 and 1 step #1 with {\arrow{latex reversed}}}}}}
\usepackage{cancel}
\usepackage[top=2cm, bottom=2.5cm, outer=2cm, inner=2cm, marginparsep=0.7cm, marginparwidth=1.5cm]{geometry}
\usepackage[read,fonts=true,save=true,mdots=true,ptmx={all=upnewtx},lang=english,xtras=true,eseese=Terragnino,mfonts={lsym=ptm,cal=dutch,scr=true},bbold=true,textcomp=no]{mworksx}
\usepackage{newtxtextM,color}

\newcommand{\nipar}{\par\noindent}
\setint[2](2)
\makeatletter

\let\eatup\@gobble
\newcommand*\MapsTo{%
  \@ifstar{\xrightarrow[\raisebox{0.25 em}{\smash{\ensuremath{\sim}}}]{}}{\xrightarrow{\raisebox{-0.25 em}{\smash{\ensuremath{\sim}}}}}%
}
\newcommand{\fadecol}[1]{\textcolor{gray}{#1}}
\newif\ifmidv\midvtrue
\newcommand{\midv}{
    \newcommand{\fade}[1]{\fadecol{##1}}
    \newcommand{\fadees}{\fade}
    \newcommand{\fadecap}[3]{\fadees{##1}{##2}{\MakeUppercase##3}}
    \newcommand{\lac}[1]{\@ifstar{<>##1}{<##1>}}
    \newcommand{\slac}[1]{<##1> }
    \newcommand{\xlac}[1]{\textcolor{red}{<}##1\textcolor{red}{>}\@ifstar{\ }{}}
    \newcommand{\xsqb}{\sqb}
    \midvtrue
}
\newcommand{\finv}{
    \newcommand{\fade}[1]{}
    \newcommand{\fadees}{\fade}
    \newcommand{\fadecap}[3]{\fadees{##1}{##2}{\MakeUppercase##3}}
    \newcommand{\lac}{\@ifstar{\expandafter\@gobble\@gobble}{}}
    \newcommand{\slac}{\lac}
    \newcommand{\xlac}{\@ifstar{\lac}{\lac}}
    \newcommand{\xsqb}{}
    \finvtrue
}
\makeatother
\upcaps

\newcommand{\opn}{\operatorname}

\setdiff[0](.25)
\setpicpath{../DRCP/}

\newcommand{\vq}{\vartheta}

\newcommand{\lsim}{\lesssim}
\midv
\renewcommand{\thesection}{\arabic{section}}

\makeatletter
\newcommand{\frl}{\@ifstar{\@frl}{\frl@}}
\newcommand{\@frl}[1]{(-\Dg)^{\xfr{#1}{2}}}
\newcommand{\frl@}[1]{(-\Dg)^{#1}}
\makeatother
\definethm{probl}{Problem}[s]
\definethm{psprobl}{Pseudo-problem}[s]
\definethm{typo}{Probable typo/misprint}[s]
\definethm{genoss}{Navier-Stokes vs. Euler Remark}[s]
\newcommand{\gsim}{\gtrsim}
\makeatletter
\newcommand{\xtag}{\@ifstar{\@tag}{\tag@}}
\newcommand{\@tag}{\stepcounter{equation}\tag{\theequation}}
\newcounter{spccnt}
\newcommand{\sxtag}[1]{\xtag{#1}\label{eq:#1}}
\newcommand{\tag@}[1]{\@tag}
\makeatother
\numberwithin{equation}{section}
\newcommand{\DSz}{\cite{DSz}}
\newcommand{\BDLSzV}{\cite{BDLSzV}}
\newcommand{\DRSz}{\cite{DRSz}}

\makeatletter
\newcommand{\xeqref}[2]{
\expandafter\ifx\csname r@#2\endcsname\relax
  \protect\G@refundefinedtrue
  \nfss@text{\reset@font\bfseries (??) R#1}%
  \@latex@warning{%
      Reference `#2' on page \thepage \space undefined%
  }%
\else
  \expandafter\expandafter\expandafter\Hy@setref@link\csname r@#2\endcsname\@empty\@empty\@nil{{\textup{\tagform@{\ref*{#2}}} R#1}\expandafter\@gobble\@firstoffive}
  \fi}
\newcommand{\EXP}{\qg}
\newcommand{\xfrl}{\@ifstar{\frl*\EXP}{\frl\EXP}}
\makeatother
\newcommand{\fD}{\mathfrak{D}}
\newcommand{\potimes}{\mathbin{\pint\otimes}}

\providecommand{\Cinf}{\ensuremath{\s{C}^\8}}
\renewcommand{\ng}{\nu}

\makeatletter
\renewcommand{\mw@proofstring}{\emph{Proof.}}
\makeatother
\newcommand{\adH}{\hat\bg}
\newcommand{\solH}{\bg}
\newcommand{\strsub}{(\str v,\str p,\str R)}
\newcommand{\strngsub}{(\strng v,\strng p,\strng R)}
\newcommand{\adsub}{(\ad v,\ad p,\ad R)}
\newcommand{\wsol}{(\sol v,\sol p)}
\newcommand{\solu}{\wsol}
\newcommand{\ad}{\hat}
\newcommand{\str}{\tilde}
\newcommand{\strng}{\check}
\newcommand{\sol}{}

\title{$L^2$-Density of Wild Initial Data for the Hypodissipative Navier-Stokes Equations}
\author{Michele Gorini\fn{GSSI - Gran Sasso Science Institute, Viale F. Crispi, 7, 67100 L'Aquila, Italy; \hsp{2cm} E-mail: michele.gorini@gssi.it}}
\date{}

\begin{document}
\maketitle
\begin{abstract}
In this paper we deal with the Cauchy problem for the hypodissipative Navier-Stokes equations in the three-dimensional periodic setting. For all Laplacian exponents $\EXP<{}^1\!/_3$, we prove non-uniqueness of dissipative $L^2_tH^\EXP_x$ weak solutions for an $L^2$-dense set of $\s{C}^\bg$ Hölder continuous wild initial data with $\EXP<\bg<\xfr13$. This improves previous results of non-uniqueness for infinitely many wild initial data (\cite{CDLDR,DR}) and generalizes previous results on density of wild initial data obtained for the Euler equations (\cite{DSz,DRSz}).
\end{abstract}
\tbold{Keywords}: Convex integration, hypodissipative Navier-Stokes equations, wild initial data
\sect{Introduction}
\numberwithin{defi}{section}
\numberwithin{lemma}{section}
\numberwithin{oss}{section}
\numberwithin{propo}{section}
\numberwithin{cor}{section}
\numberwithin{teor}{section}
The existence of dissipative solutions for the Euler equations
\[\begin{sistema}
\pd_tv+\opn{div}(v\otimes v)+\grad p=0 \\
\opn{div}v=0
\end{sistema} \xtag{Eul}\label{eq:Eul}\]
with regularity lower than $\s{C}^0_t\s{C}^{\xfr13}_x$ has been investigated deeply. After the pioneering works of Scheffer \cite{Sche} (on the plane $\R^2$) and Shnirelman \cite{Shn} (on the periodic torus $\T^2$), De Lellis and Székelyhidi, in \cite{DLSz09}, introduced the convex integration technique (first used by Nash in \cite{N} and Kuiper in \cite{K} in the context of isometric embeddings, and formalized in a more general setting by Gromov in \cite{G}) in this setting, proving the existence of nontrivial compactly supported $L^\8_tL^2_x$ weak solutions of \eqref{eq:Eul} in $\R^n$ for any $n$. The subsequent paper \cite{DLSz10} provided a proof of the non-uniqueness of weak solutions satisfying the weak energy inequality
\[\xints{\T^3}{}\2\abs{v(x,t)}^2\diff x\leq\xints{\T^3}{}\2\abs{v(x,0)}^2\diff x, \xtag{EulAdm}\label{eq:EulAdm}\]
i.e. dissipating the total kinetic energy. We call such solutions \emph{dissipative} or \emph{admissible}.
\nipar Both of these papers use a Baire category argument, proving that such solutions constitute the set of continuity points of a Baire-1 map. This implies that such solutions not only exist, but are ``typical'' in the sense of category.
\nipar These results were the first steps in the resolution of the second part of Onsager's conjecture from \cite{O}.
\begin{teorspec}[Onsager's]{Conjecture}[thm:spec:Onsager]
Let $(v.p)$ be a weak solution of \eqref{eq:Eul} and define the total kinetic energy as
\[E(t)\coloneq\xfr12\xints{\T^3}{}\2\abs{v(x,t)}^2\diff x\]
If $v\in\s{C}^\bg$ for $\bg>\xfr13$, then the energy is a conserved quantity, i.e. $E(t)\equiv E(0)$.
\nipar By contrast, for any $\bg<\xfr13$, there exist $\s{C}^\bg$ weak solutions of \eqref{eq:Eul} which do not conserve the energy.
\xend{teorspec}
The first part was proved in \cite{CET}, and refined for more general spaces in \cite{CCFS,FW}.
\nipar In \cite{DLSz}, De Lellis and Székelyhidi introduced a more constructive approach which allowed to obtain infinitely many dissipative continuous solutions (see also \cite{DLSzC0}), and then infinitely many $\s{C}^{\xfr{1}{10}-\eg}$-Hölder solutions in \cite{DLSzC1/10}.
\nipar In \cite{BDLISz}, Isett was able to improve the Hölder exponent for the existence of non-conservative solutions to $\xfr15-\eg$ and, after the introduction of Mikado flows by Daneri and Székelyhidi in \cite{DSz}, he completed the proof of Onsager's conjecture above in \cite{I} by showing existence of infinitely many Hölder continuous dissipative solutions in the class $\s{C}^0_t\s{C}^\bg_x$, for all $\bg<\xfr13$. This result was later improved to dissipative solutions in the same class in \cite{BDLSzV}.
\nipar In the class of admissible solutions, weak-strong uniqueness holds, as proved in \cite{W}: $\s{C}^1$ solutions are unique, and moreover, if such a solution exists, any $L^\8_tL^2_x$ solution with the same initial data which is admissible coincides with the $\s{C}^1$ solution.
\nipar However, cleverly adapting and improving the convex integration technique of the above-mentioned papers, the existence of infinitely many initial data giving rise to admissible solutions in the class $\s{C}^\bg$ for $\bg<{}^1\!/_{16}$ was proved in \cite{D}. In \cite{DSz} and \cite{DRSz}, the following topologically stronger statement was proved: the set of $\s{C}^\bg$ initial data giving rise to admissible solutions is dense in $L^2(\T^3)$. It was proved in \cite{DSz} for $\bg<{}^1\!/_5$, and in \cite{DRSz} for $\bg<{}^1\!/_3$.
\nipar Removing the admissibility condition \eqref{eq:EulAdm} leads to non-uniqueness for any $\Cinf$ initial datum, as proved in \cite{DH}.
\nipar For the Navier-Stokes equations
\[\begin{sistema}
\pd_tv+\opn{div}(v\otimes v)+\grad p=\Dg v \\
\opn{div}v=0
\end{sistema}, \xtag{NS}\label{eq:NS}\]
in \cite{L}, Leray proved the existence of global weak solutions satisfying the following energy inequality:
\[\xints{\T^3}{}\2\abs{v(x,t)}^2\diff x+2\per\e\xints 0t\xints{\T^3}{}\2\abs{\grad v(x,s)}^2\diff x\diff s\leq\xints{\T^3}{}\2\abs{v(x,0)}^2\diff x\qquad\qquad\text{a.e. }t>0. \xtag{NSAdm}\label{eq:NSAdm}\]
Such solutions are called \emph{Leray solutions} or \emph{Leray-Hopf solutions}, and we will call them \emph{admissible solutions of \eqref{eq:NS}}. The strategy employed in \cite{L} can easily be adapted to prove the existence of solutions to \eqref{eq:Eul} satisfying \eqref{eq:EulAdm}, which are therefore called by the same names.
\nipar Thanks to the Ladyzhenskaya-Prodi-Serrin regularity theory, weak-strong uniqueness for \eqref{eq:NS} holds for $L^r_tL^s_x$ solutions, where $\xfr ds+\xfr2r=1$, as proved in \cite{P} for $d=3$ and in \cite{Se} for the general case.
\nipar The uniqueness or non-uniqueness of solutions to \eqref{eq:NS} satisfying \eqref{eq:NSAdm} is still a long-standing open problem. The latest step in this regard is \cite{ABC} where, introducing a body force in the equation, the authors exhibit two distinct admissible solutions on $\R^3$.
\nipar Several non-uniqueness results have been obtained for non-admissible solutions, i.e. in the absence of the energy inequality \eqref{eq:NSAdm}. In \cite{BV}, the authors prove the existence, for any smooth energy profile $e$, of $\s{C}^0_tH^\bg_x$ solutions with kinetic energy $e$, i.e.
\[\xints{\T^3}{}\abs{v(x,t)}^2\diff x=e(t),\]
for some parameter $\bg$. This implies non-uniqueness for the zero initial datum. Choosing a non-increasing $e$, this also implies the existence of solutions of \eqref{eq:NS} satisfying \eqref{eq:EulAdm}.
\nipar It is known that such $\bg$ cannot be too large, since $\bg={}^1\!/_2$ implies weak-strong uniqueness by \cite{CF}. In arbitrary dimensions, $\bg={}^1\!/_2$ is in fact a threshold for weak-strong uniqueness. Indeed, in Terence Tao's blog post \cite{Tao}, the non-uniqueness of $H^1_tH^s_x$ for any $s<{}^1\!/_2$ has been proven to hold on $\T^d$ where $d=d(s)$ is sufficiently large.
\nipar In the subsequent work \cite{BCV}, a ``gluing'' theorem is proved: given any two strong solutions $u_1,u_2\in\s{C}^0_t\dot H^3_x([0,T],\T^3)$, there exists a weak solution $v\in\s{C}^0_t(H^\bg_x\cap W^{1,1+\bg}_x)$, with a set of singular times $\Sg$ having dimension $\dim_H\Sg\leq1-\bg$, which coincides with $u_1$ on $[0,{}^T\!/{\!}_3]$ and with $u_2$ on $[{}^2\!/{\!}_3T,T]$. The parameter $\bg$ is not quantified in \cite{BCV}, but the strategy therein allows it to reach at most a value slightly above $10^{-3}$.
\nipar Reducing the regularity in time can lead to better spatial regularity. Indeed, in \cite{CL}, the authors prove an approximation result: given a smooth divergence-free field $v$, we can approximate it in $L^p_tL^\8_x\cap L^1_tW^{1,\8}_x$ with a solution of \eqref{eq:NS}, for any $p<2$. The singular set of these solutions is also of low dimension. The strategy of \cite{CL} can be extended to $L^p_tL^\8_x\cap L^s_tW^{1,q}_x$, where once again $p<2$, and we have that $s<2$ and $q<q_{max}(s,p)$, where $q_{max}(s,p)\to1$ if $s\to2$ or $p\to2$.
\nipar In this paper we are interested in studying how introducing a fractional dissipative term in \eqref{eq:Eul} may affect the uniqueness or ill-posedness of the Cauchy problem. More specifically, we consider the hypodissipative Navier-Stokes equations
\[\begin{sistema}
\pd_tv+\opn{div}(v\otimes v)+\grad p+(-\Dg)^\EXP v=0 \\
\opn{div}v=0
\end{sistema}, \xtag{FNS}\label{eq:FNS}\]
with exponent $\EXP<\xfr13$ and admissibility condition
\[\xfr12\xints{\T^3}{}\2\abs{v(x,t)}^2\diff x+\xints0t\xints{\T^3}{}\2\abs{\xfrl*v(x,s)}^2\diff x\diff s\leq\xfr12\xints{\T^3}{}\2\abs{v(x,0)}^2\diff x. \xtag{FNSAdm}\label{eq:FNSAdm}\]
The previously cited proof of \cite{L} can be adapted to show existence for admissible solutions in this case as well. This result is stated in \kcref{thm:teor:ExistLeraySol}. In this case as well, admissible solutions are also known as Leray solutions.
\nipar In \cite{CDLDR} and \cite{DR} it was proved that there exist infinitely many $\s{C}^\bg$ initial data for $\EXP<\bg<{}^1/_3$ which generate infinitely many $\s{C}^0_t\s{C}^\bg_x$ solutions which, by \cite{CD}, are in fact $\s{C}^\bg_{t,x}$. In fact, \cite{CDLDR} produces such data for $\EXP$ up to $^1/_2$, with solutions that can only be proved to be admissible (i.e. satisfy \eqref{eq:FNSAdm}) for $\EXP<{}^1/_5$.
\nipar Inspired by the results in \cite{DSz} and \cite{DRSz}, which prove the $L^2$-density of such initial data for the Euler equations \eqref{eq:Eul}, respectively for $\bg<{}^1/_5$ and for $\bg<{}^1/_3$, we investigate the existence of an $L^2$-dense class of $\s{C}^\bg$ wild initial data (namely data for which non-uniqueness holds) for admissible solutions to \eqref{eq:FNS} in $L^2_tH^\EXP_x$. The strategy proposed in \cite{DSz} and \cite{DRSz} provides a quantitative criterion for non-uniqueness based on the existence of approximate solutions called \emph{adapted subsolutions}.
\nipar Here, we explore and extend that strategy to the hypodissipative Navier-Stokes equations. The main issue with respect to the Euler setting is to control the dissipative term in the energy. Our main results are the following.
\xbegin{teor}[\texorpdfstring{$\s{C}^\bg$}{Cbg} weak solutions with data close to \texorpdfstring{$L^2$}{L2} functions and time of admissibility][thm:teor:TimeExist]
Let $\EXP<\bg<\xfr13,w\in L^2(\T^3)$. Then, for all $\hg>0$, there exist a time $T=T(\hg)>0$, an initial datum $w_\hg\in\s{C}^\bg(\T^3)$ such that $\|w_\hg-w\|_{L^2}<\hg$, and infinitely many weak solutions $v_\hg\in\s{C}^0([0,T],\s{C}^\bg(\T^3)\!)$, with initial datum $v_\hg|_{t=0}=w_\hg$, which satisfy \eqref{eq:EulAdm} on $[0,T]$, but can be proved to satisfy \eqref{eq:FNSAdm} (i.e to be admissible) only on $[0,T(\hg)]$. Moreover
\[\lim_{\hg\to0}T(\hg)=0.\]
\xend{teor}
The fact that the admissibility condition in \kcref{thm:teor:TimeExist} cannot be guaranteed to hold for $\s{C}^\bg$ solutions on a fixed set of times is due to the necessity of controlling the dissipation term in the energy.
\xbegin{defi}[Wild initial data][thm:defi:WID]
Let $X$ be a function space. A divergence-free vector field $w\in L^2(\T^3)$ is a \emph{$(\EXP,X,T)$-wild initial datum} for \eqref{eq:FNS} if there exist infinitely many weak solutions $v:[0,T]\x\T^3\to\R^3$ of \eqref{eq:FNS} such that $v\in X$, $v(x,0)=w(x)$ a.e. in $\T^3$, and the admissibility condition \eqref{eq:FNSAdm} holds on $[0,T]$. The set of such data is denoted by $W_{\EXP,X,T}$. If $X=L^\8_t\s{C}^\bg_x$, we will speak of $(\EXP,\bg,T)$-wild data, and of the set $W_{\EXP,\bg,T}$.
\xend{defi}
As a consequence of \kcref{thm:teor:TimeExist}, we obtain the following corollary.
\xbegin{cor}[Density of wild initial data -- Hölder solutions][thm:cor:DensWIDUnT]
The set $\bigcup_TW_{\EXP,\bg,T}$ is dense in the set of divergence-free $L^2$ vector fields, for all $\EXP<\bg<\xfr13$.
\xend{cor}
Moreover, by taking a solution $v_\hg$ as given by \kcref{thm:teor:TimeExist}, and continuing it with a Leray solution $\tilde v_\hg:[T(\hg),\8]\x\T^3\to\R^3$ with datum $\tilde v_\hg(T(\hg)\!)=v_\hg(T(\hg)\!)$, as provided by \kcref{thm:teor:ExistLeraySol}, we obtain the following.
\xbegin{teor}[Density of wild initial data -- Sobolev solutions][thm:teor:DensWIDNoHoeld]
$W_{\EXP,L^2_tH^\EXP_x,T}$ is dense in the sef ot divergence-free $L^2$ vector fields, for all $\EXP<\xfr13,T>0$.
\xend{teor}
The general strategy of the paper is to define suitable relaxations of the notion of solution (the so-called ``subsolutions'' of Section \ref{SubsolDefi}), and approximating one kind of subsolution with another one which is closer to the notion of solution. This is done constructing sequences of subsolutions that converge, in an appropriate sense, to a ``stronger'' subsolution. We will need two such approximations, and therefore two convex integration schemes: the first one will converge to a subsolution which is a solution at $t=0$, and has the $\s{C}^1$ norm of the velocity blowing up at a controlled rate at $t=0$. The second approximation will lead to a weak solution.
\nipar The paper is organized as follows. Section \ref{Prelims} establishes some notations and contains a few useful preliminary results. Section \ref{SubsolDefi} introduces three kinds of subsolutions, namely strict, strong, and adapted. Section \ref{Equestimates} contains some important definitions and relations used in the following sections. Section \ref{Propos} states the approximation results that lie at the heart of the proof. Section \ref{FinArg} deduces \kcref{thm:teor:TimeExist} from those approximation results. Section \ref{S2S} shows how one can approximate strict subsolutions with strong ones. Sections \ref{Gluing} and \ref{MPS} contain the two substeps of each convex integration step, respectively a gluing step and a perturbation step. Sections \ref{S2AProof} and \ref{A2SProof} then prove the other approximation results, namely the approximation of strong subsolutions with adapted ones, and that of adapted subsolutions with weak solutions.

\sect{Preliminaries}\label{Prelims}
Throughout the paper, we will use the following notations:
\bi
\item $\s{S}^{3\x3}$ are the symmetric 3-by-3 matrices; within this set, $\s{S}^{3\x3}_+$ are the positive definite ones, $\s{S}^{3\x3}_0$ are the traceless ones, and $\s{S}^{3\x3}_{\geq0}$ are the positive semidefinite ones.
\item If $R\in\s{S}^{3\x3}$, we decompose it as
\[R=\xfr13\opn{tr}R\opn{Id}+\pint R=\rg\opn{Id}+\pint R,\]
where $\pint R\in\s{S}^{3\x3}_0$ is the \emph{traceless part of $R$}.
\item For scalar functions $f$, we write $\grad f\coloneq(\pd_1f,\pd_2f,\pd_3f)\eqcolon\fD f$;
\item However, for vector fields $v$, we define $\fD v$ so that $(\fD v)_{ij}=\pd_jv_i$, whereas $\grad v=(\fD v)^T$; with these choices, $(v\per\grad)v=\fD v\per v=v\per\grad v$;
\item In a similar fashion, for tensor fields $S$, $\fD S$ is defined by $(\fD S)_{ijk}=\pd_kS_{ij},$ whereas $\grad S$ is defined by $(\grad S)_{ijk}=\pd_iS_{jk}$;
\item The Hölder norms are defined as follows
\begin{gather*}
\norm{f}_0\coloneq\sup\abs{f(x)}, \qquad [f]_k\coloneq\max_{\sstack{\bg\in*\N^3 \\ |\bg|=k}}\sup\abs{\pd_\bg f}, \qquad [f]_\ag\coloneq\sup_{x,y}\xfr{\abs{f(x)-f(y)}}{\abs{x-y}^\ag} \\
\norm{f}_k\coloneq\norm{f}_0+\sum_{i=1}^k[f]_i, \qquad\qquad \norm{f}_{k+\ag}\coloneq\norm{f}_k+\max_{\sstack{|\bg|=k \\ \bg\in*\N^3}}[\pd_\bg f]_\ag,
\end{gather*}
for $k\in\N,\ag\in(0,1)$.
\ei
Concerning the Hölder norms, we recall the following standard inequalities.
\xbegin{lemma}[Hölder norm inequalities][thm:lemma:HNI]
For $0\leq s\leq r$ and $f,g:\T^3\to\R^d$
\begin{align*}
[fg]_r\leq{}&C(r)([f]_r\|g\|_0+\|f\|_0[g]_r) \xtag{HProd}\label{eq:HProd}\\\relax
[f]_s\leq{}&C(r,s)\|f\|_0^{1-\xfr sr}[f]_r^{\xfr sr}. \xtag{HInterp}\label{eq:HInterp}
\end{align*}
Moreover, for $f:\T^3\to S\sbse\R^d$ and $\Yg:S\to\R$:
\begin{align*}
[\Yg\circ f]_m\leq{}&K(d,m)([\Yg]_1\|\fD f\|_{m-1}+\|\grad\Yg\|_{m-1}\|f\|_0^{m-1}\|f\|_m) \xtag{HCompos1}\label{eq:Compos1} \\\relax
[\Yg\circ f]_m\leq{}&K(d,m)([\Yg]_1\|\fD f\|_{m-1}+\|\grad\Yg\|_{m-1}[f]_1^m). \xtag{HCompos2}\label{eq:Compos2}
\end{align*}
Finally, for all $s,r\geq0$:
\[\begin{aligned}
\norm{f\ast\fg_\ell}_{r+s}\leq{}&C(r,s)\ell^{-s}\|f\|_r \\
\norm{f-f\ast\fg_\ell}_r\leq{}&C(r,s)\ell^1\|f\|_{r+1} \\
\norm{f-f\ast\fg_\ell}_r\leq{}&C(r,s)\ell^2\|f\|_{r+2} \\
\norm{(fg)\ast\fg_\ell-(f\ast\fg_\ell)(g\ast\fg_\ell)}_r\leq{}&C(r,s)\ell^{2-r}\|f\|_1\|g\|_1,
\end{aligned}
\xtag{Mollif}\label{eq:Mollif}\]
where $\fg$ is a standard mollification kernel, i.e. $\fg\in\Cinf_c(B_1;[0,1])$ and $\int\fg=1$, and $\fg_\ell\coloneq\xfr{1}{\ell^3}\fg(\xfr\per\ell)$.
\xend{lemma}
For functions $f,g:\T^3\x[0,T]\to\R^d$, we denote their time-slices by $f_t(x)\coloneq f(x),g_t(x)\coloneq g(x)$. The above lemma can then be applied to the time slices of time-dependent vector fields, e.g. the velocities of subsolutions, with the notation $\|f(t,\per)\|_{\s{C}^r},[f(t,\per)]_{\s{C}^r}$ for the (semi)norms of the slices. By taking supremum norms in time, the above inequalities can be formulated with $\s{C}^0_t\s{C}^r_x$ norms.
\nipar We now introduce Mikado flows, the basic building blocks of the perturbations, and the important Stationary Phase Lemma. The proofs of the following results can be found in \cite{DSz}.
\xbegin{lemma}[Mikado flows][thm:lemma:MikFlows]
For any compact subset $\s{N}\sbs\sbs\s{S}_+^{3\x3}$ there exists a smooth vector field $W:\s{N}\x\T^3\to\R^3$ such that, for every $R\in\s{N}$
\[\begin{sistema}
\opn{div}_\jg(W(R,\jg)\otimes W(R,\jg)\!)=0 \\
\opn{div}_\jg W(R,\jg)=0
\end{sistema}, \xtag{MikDiff}\label{eq:MikDiff}\]
and
\[\begin{sistema}
\fint_{\T^3}W(R,\jg)\diff\jg=0 \\
\fint_{\T^3}W(R,\jg)\otimes W(R,\jg)\diff\jg=R
\end{sistema}.  \xtag{MikInt}\label{eq:MikInt}\]
Using Fourier series in $\jg$ and the above integral and differential relations, we obtain that
\begin{align*}
W(R,\jg)={}&\sum_{\mcl{k\in*\Z^3\ssm\{0\}}}a_k(R)A_ke^{ik\per\jg} \xtag{MikFS}\label{eq:MikFS} \\
W(R,\jg)\otimes W(R,\jg)={}&R+\sum_{\mcl{k\in*\Z^3\ssm\{0\}}}C_k(R)e^{ik\per\jg}, \xtag{MikTPFS}\label{eq:MikTPFS}
\end{align*}
where the coefficients $a_k,C_k\in\Cinf$, the $A_k$ satisfy $A_k\per k=0,|A_k|=1$, the $C_k$ satisfy $C_kk=0$, and moreover
\begin{align*}
\sup_{R\in*\s{N}}\abs{\fD_R^Na_k(R)}={}&\norm{a_k}_{\s{C}^N(\s{N})}\leq\xfr{C(\s{N},N,m)}{|k|^m} \xtag{Mikak}\label{eq:Mikak} \\
\sup_{R\in*\s{N}}\abs{\fD_R^NC_k(R)}={}&\norm{C_k}_{\s{C}^N(\s{N})}\leq\xfr{C(\s{N},N,m)}{|k|^m}. \xtag{MikCk}\label{eq:MikCk}
\end{align*}
\xend{lemma}
In Section \ref{S2S}, we will need the fact that, if we set
\[U(R,\jg)\coloneq\sum_ka_k(R)\xfr{ik\x A_k}{|k|^2}e^{ik\per\jg}, \xtag{AntiCurlW}\label{eq:AntiCurlW}\]
then we have that $\opn{curl}_\jg U=W$. Indeed
\begin{align*}
\opn{curl}_\jg U(R,\jg)={}&\sum_{k\ell mn}\eg_{\ell mn}\pd_m\pa{a_k(R)\xfr{ik\x A_k}{|k|^2}e^{ik\per\jg}}_ne_\ell=\sum_{\mcl{k\ell mnpq}}\eg_{\ell mn}\eg_{npq}a_k(R)ik_p(A_k)_q|k|^{-2}\per ik_me^{ik\per\jg}e_\ell={} \\
{}={}&-\sum_{\mcl{k\ell mpq}}(\dg_{\ell p}\dg_{mq}-\dg_{\ell q}\dg_{mp})a_k(R)k_p|k|^{-2}(A_k)_qk_me^{ik\per\jg}e_\ell={} \\
{}={}&\sum_ka_k(R)A_ke^{ik\per\jg}-\sum_{kpq}a_k(R)k_pk_q|k|^{-2}(A_k)_qe^{ik\per\jg}e_p={} \\
{}={}&W(R,\jg)-\sum_ka_k(R)\underbrace{(k\per A_k)}_0\xfr{k}{|k|^2}e^{ik\per\jg}.
\end{align*}
We now introduce a certain ``anti-divergence operator'' which will be used to obtain the new Reynolds stress $R$ in the various approximation results.
\xbegin{defi}[Anti-divergence \texorpdfstring{$\s{R}$}{R}][thm:defi:AntiDivR]
Define the operator $\diamond$ so that
\[\begin{sistema}
\Dg\diamond v=v-\fint_{\T^3}v\diff x \\
\int_{\T^3}\diamond v=0
\end{sistema},\]
and then define
\[\s{R}v\coloneq\xfr14(\fD\s{P}\diamond v+(\fD\s{P}\diamond v)^T)+\xfr34(\fD\diamond v+(\fD\diamond v)^T)-\xfr12(\opn{div}\diamond v)\opn{Id}, \xtag{DefR}\label{eq:DefR}\]
$\s{P}$ being the Leray projection onto divergence-free fields with zero average.
\xend{defi}
This operator satisfies the following properties.
\xbegin{lemma}[Divergence and \texorpdfstring{$\s{R}$}{R}][thm:lemma:DivAndR]
For any $\Cinf$ vector field $v$, $\s{R}v\in\s{S}^{3\x3}_0$ is symmetric and trace-free, and moreover
\[\opn{div}\s{R}v=v-\fint_{\T^3}v\diff x, \xtag{divR}\label{eq:divR}.\]
\xend{lemma}
Moreover, we have the following statement, which we will use numerous times throughout the paper.
\xbegin{lemma}[Stationary Phase Lemma][thm:lemma:SPL]
Let $\ag\in(0,1),N\geq1$. Let $a\in\Cinf(\T^3),\Fg\in\Cinf(\R^3,\R^3)$ be smooth functions and assume
\[\xfr1K\leq\abs{\fD\Fg}\leq K\qquad\text{on }\T^3.\]
Then
\[\abs{\xints{\T^3}{}\2a(x)e^{ik\per\Fg}\diff x}\leq C(K,N)\xfr{\norm{a}_N+\norm{a}_0\norm{\Fg}_N}{\abs k^N}, \xtag{SPL-Int}\label{eq:SPL-Int}\]
and for the operator $\s{R}$ of \eqref{eq:DefR} above we have that
\[\norm{\s{R}(a(x)e^{ik\per\Fg})}_\ag\leq C(\ag,K,N)\pa{\xfr{\norm{a}_0}{\abs k^{1-\ag}}+\xfr{\norm{a}_{N+\ag}+\norm{a}_0\norm{\Fg}_{N+\ag}}{\abs k^{N-\ag}}}. \xtag{SPL-R}\label{eq:SPL-R}\]
\xend{lemma}
We now recall some classical estimates regarding fractional laplacians.
\xbegin{teor}[Fractional laplacian and Hölder norms][thm:teor:FrLHoeldNorms]
Let $\gg,\eg>0$ and $\bg\geq0$ such that $2\gg+\bg+\eg\leq1$, and let $f:\T^3\to\R^3$. If $f\in\s{C}^{0,2\gg+\bg+\eg}$, then $\frl\gg f\in\s{C}^\bg$, moreover there exists a constant $C=C(\eg)$ such that
\[\norm{\frl\gg f}_\bg\leq C(\eg)[f]_{2\gg+\bg+\eg}. \xtag{FrLH}\label{eq:FrLH}\]
Moreover, for every $\gg\in(0,1)$, $\eg>0$ such that $0<\gg+\eg\leq1$, and $f$ as above, then there exists $C=C(\eg)>0$ such that
\[\xints{\T^3}{}\2\abs{\frl*\gg f}^2(x)\diff x\leq C(\eg)[f]_{\gg+\eg}^2\qquad\VA f\in\s{C}^{\gg+\eg}(\T^3). \xtag{FrLInt}\label{eq:FrLInt}\]
\xend{teor}
For a proof, see \cite[Theorem B.1 and Corollary B.1]{DR}.
\nipar Continuing, we recall some elementary calculations for the reader's convenience. With the definitions we gave for $\grad,\fD$, setting $D_t\stp v\coloneq\pd_t+v\per\grad$, we have that
\begin{align*}
\grad e^{ik\per\Fg}={}&i\grad\Fg\per ke^{ik\per\Fg}=ie^{k\per\Fg}k\per\fD\Fg \xtag{GradExp}\label{eq:GradExp} \\
D_t\stp v(\fD\Fg)={}&\fD(D_t\stp v\Fg)-\fD\Fg\cdot\fD v. \xtag{MDGrad}\label{eq:MDGrad}
\end{align*}
Observing that $\fD\Fg\fD\Fg^{-1}=\opn{Id}$ and thus $0=D_t\stp v(\fD\Fg\fD\Fg^{-1})=D_t\stp v(\fD\Fg)\cdot\fD\Fg^{-1}+\fD\Fg\cdot D_t\stp v(\fD\Fg^{-1})$, we can see that
\[D_t\stp v\fD\Fg^{-1}=-\fD\Fg^{-1}D_t\stp v(\fD\Fg)\fD\Fg^{-1}=(\grad\Fg^{-1}\grad v)^T-\fD\Fg^{-1}\per\fD D_t\stp v\Fg\per\fD\Fg^{-1}. \xtag{MDInvGrad}\label{eq:MDInvGrad}\]
We continue by recalling some classical estimates on the transport and transport-diffusion equations, which can be found e.g. in \cite[Proposition B.1]{BDLSzV} (transport) and \cite[Proposition 3.3]{DR} (transport-diffusion).
\xbegin{propo}[Estimates on the transport equation][thm:propo:EstTrasp]
Assume $|t-t_0|\|v\|_1\leq1$. Then, any solution $f$ of
\[\begin{sistema}
(\pd_t+v\per\grad)f=g \\
f(\per,0)=f_0
\end{sistema}\]
satisfies
\begin{align*}
\norm{f(t)}_0\leq{}&\norm{f_0}_0+\xints{t_0}{t}\norm{g(\per,\tg)}\diff s \\
\norm{f(t)}_\ag\leq{}&e^\ag\pa{\norm{f_0}_\ag+\xints{t_0}{t}\norm{g(\per,s)}_\ag\diff s}
\end{align*}
for all $0\leq\ag\leq1$ and, more generally, for any $N\geq1$ and $0\leq\ag<1$
\[[f(t)]_{N+\ag}\lsim[f_0]_{N+\ag}+|t|[v]_{N+\ag}[f_0]_1+\xints{t_0}{t}\pa[b]{[g(s)]_{N+\ag}+(t-s)[v]_{N+\ag}[g(s)]_1}\diff s.\]
Define $\Fg(t,\per)$ to be the inverse of the flux $X$ of $v$ starting at time $t_0$ as the identity (i.e. ${}^d\2/\!_{dt}\,X=v(X,t)$ and $X(x,t_0)=x$). Under the same assumptions as above we have that
\begin{align*}
\norm{\grad\Fg(t)-\opn{Id}}_0\lsim{}&|t|[v]_1 \\\relax
[\Fg(t)]_N\lsim{}&|t|[v]_N\quad\VA N\geq2.
\end{align*}
\end{propo}
\xbegin{propo}[Estimates on the transport-diffusion equation][thm:propo:EstTraspDiff]
Assume $0\leq(t-t_0)[v]_1\leq1$. Then, any solution of
\[\case{
(\pd_t+v\per\grad+\xfrl)u=f & \text{in }\T^3\x(t_0,T) \\
u(\per,t_0)=u_0 & \text{in }\T^3
}\]
satisfies
\[\norm{u(t)}_\ag\leq e^\ag\pa{\norm{u_0}_\ag+\xints{t_0}{t}\norm{f(\per,s)}_\ag\diff s}\]
for all $0\leq\ag\leq1$ and, more generally, for any $N\geq1$ and $0\leq\ag<1$
\[[u(t)]_{N+\ag}\lsim[u_0]_{N+\ag}+(t-t_0)[v]_{N+\ag}[u_0]_1+\xints{t_0}{t}\pa[b]{[f(s)]_{N+\ag}+(t-s)[v]_{N+\ag}[f(s)]_1}\diff s,\]
where the implicit constants depends only on $N,\ag$.
\xend{propo}
To conclude this section, we recall classical Schauder estimates (see e.g. the book \cite{GT}), which will be used in several places in this paper.
\xbegin{lemma}[Schauder estimates][thm:lemma:Schauder]
For any $\ag\in(0,1)$ and any $m\in\N$, there exists a constant $C(\ag,m)$ with the following properties. If $\fg,\yg:\T^3\to\R$ are the unique solutions of
\[\begin{sistema}
\Dg\fg=f \\
\fint\fg=0
\end{sistema}
\qquad\qquad
\begin{sistema}
\Dg\yg=\opn{div}F \\
\fint\yg=0
\end{sistema},\]
then
\[\norm{\fg}_{m+2+\ag}\leq C(m,\ag)\norm{f}_{m+\ag}\qquad\qquad\norm{\yg}_{m+1+\ag}\leq C(m,\ag)\norm{F}_{m+\ag}.\]
\end{lemma}
\clearpage

\sect{Approximate solutions}\label{SubsolDefi}
For the proof of \kcref{thm:teor:TimeExist}, we begin by introducing the various notions of subsolutions needed to perform the convex integration schemes.
\nipar The first notion of subsolution is very similar to the one used in \cite{DRSz,DSz,DLSz}.
\xbegin{defi}[Subsolutions and strict subsolutions][thm:defi:Strict]
A \emph{subsolution} is a triple $(v,p,R):\T^3\x(0,T)\to\R^3\x\R\x\s{S}^{3\x3}_{\geq0}$ such that $v\in L^2_{loc}$, $R\in L^1_{loc}$, $p$ is a distribution, the equations
\[\begin{sistema}
\pd_tv+\opn{div}(v\otimes v)+\grad p+\xfrl v=-\opn{div}R \\
\opn{div}v=0
\end{sistema} \xtag{FNSR}\label{eq:FNSR}\]
hold in the sense of distributions in $\T^3\x(0,T)$, and moreover $R\geq0$ a.e., i.e. it is positive semidefinite a.e.. If $R\in\s{S}^{3\x3}_+$ a.e., then the subsolution is said to be \emph{strict}.
\xend{defi}
The next notion of subsolution extends the ones of \DSz\ and \DRSz. As in \DRSz, the Reynolds stress is controlled by a power of the trace. However, the exponent $\gg$ will only act on the ``reduced'' trace $\rg\Wg^{-1}$, where $\Wg>0$ is a constant whose role is explained in Section \ref{Equestimates}.
\xbegin{defi}[Strong subsolutions][thm:defi:Strong]
A \emph{strong subsolution with parameters $\gg,\Wg>0$} is a subsolution $(v,p,R)$ such that in addition $\opn{tr}R$ is a function of $t$ only and, if
\[\rg(t)\coloneq\xfr13(\opn{tr}R)(t)\qquad\qquad\varrho(t)\coloneq\xfr{\rg(t)}{\Wg},\]
then
\[\abs{\pint R(x,t)}\leq\Wg\varrho^{1+\gg}(t)\qquad\VA(x,t). \xtag{Strong}\label{eq:Strong}\]
\end{defi}
\xbegin{oss}[On strength and parameters][thm:oss:StrengthParams]
In our schemes $\varrho$ will be sufficiently small so that in particular $\varrho^\gg\leq r_0$, where $r_0$ is the geometric constant in \cite[Definition 3.2]{DSz}, thus leading to the conclusion that \eqref{eq:Strong} implies that our strong subsolutions are also strong in the sense of \DSz, provided $\Wg=O(1)$ (specifically $\Wg\varrho^\gg\leq r_0$). Note also that, if $(v,p,R)$ is a strong subsolution for some parameters $\gg,\Wg>0$ with $\varrho<1$, then it is also a strong subsolution for any $0<\gg'<\gg$ with the same $\Wg$.
\xend{oss}
The last notion of subsolution has vanishing Reynolds stress at time $t=0$ and the $\s{C}^1$-norms blow up at certain rates as the Reynolds stress goes to zero. Such adapted subsolutions have been introduced in \cite{DSz,DRSz}. The blow-up rate in this paper is analogous to the one of \cite{DRSz}. Differently from \cite{DRSz}, the blow-up is controlled by the ``reduced'' trace $\varrho$ rather than the ``full'' trace $\rg$, and the estimates include a power of $\Wg$.
\xbegin{defi}[Adapted subsolutions][thm:defi:Adapted]
Given $\gg,\Wg>0,0<\bg<\xfr13$, and $\ng$ satisfying
\[\ng>\xfr{1-3\bg}{2\bg}, \xtag{ANuBeta}\label{eq:ANuBeta}\]
we call a triple $(v,p,R)$ a \emph{$\s{C}^\bg$-adapted subsolution on $[0,T]$ with parameters $\gg,\Wg,\ng$} if $(v,p,R)\in\Cinf(\T^3\x\pasq{0,T})\cap\s{C}(\T^3\x[0,T])$ is a strong subsolution with parameters $\gg,\Wg$ with initial datum
\[v(\per,0)\in\s{C}^\bg(\T^3)\quad\text{and}\quad R(\per,0)\equiv0, \xtag{AIV}\label{eq:AIV}\]
and, setting $\rg(t)\coloneq\xfr13\opn{tr}R(x,t)$ and $\varrho\coloneq\rg\Wg^{-1}$, for all $t>0$ we have that $\rg(t)>0$ and there exist $\ag\in(0,1)$ and $C\geq1$ such that
\begin{align*}
\norm{v}_{1+\ag}\leq{}&C\Wg^{\xfr12}\varrho^{-(1+\ng)} \xtag{ASv}\label{eq:ASv} \\
\abs{\pd_t\varrho}\leq{}&C\Wg^{\xfr12}\varrho^{-\ng}. \xtag{ASTr}\label{eq:ASTr}
\end{align*}
\xend{defi}

\sect{Strategy of the proof}\label{Equestimates}
The remainder of this paper closely follows the convex integration strategy adopted by \DRSz\ in the Euler setting.
\nipar Section \ref{Propos} states results that allow us to approximate one kind of subsolution (as defined in the previous section) with another. One of those results uses the parameters we will introduce in this section in \eqref{eq:DefParams}.
\nipar Section \ref{FinArg} proves the main theorem starting from the results of Section \ref{Propos}.
\nipar Section \ref{S2S} shows how to obtain a strong subsolution from a strict one. By iterating Sections \ref{Gluing}-\ref{MPS}, we produce sequences $(v_q,p_q,R_q)$ of strong subsolutions which converge to a $\s{C}^\bg$-adapted subsolution in Section \ref{S2AProof}, and to a weak solution in Section \ref{A2SProof}.
\nipar In passing from one subsolution to another, the $\s{C}^0$ and $\s{C}^1$ norms of the various subsolutions are estimated in terms of parameters $(\dg_q,\lg_q)$, where $\dg_q^{\nicefrac12}$ is the amplitude (in space) of $w_q\coloneq v_q-v_{q-1}$, and $\lg_q$ is the oscillation frequency (in space) of $w_q$. The parameters, however, are partially different from those chosen in \cite{DRSz} and closer to the ones used in \cite{DSz}. More precisely, we define
\[\lg_q\coloneq2\pg\lceil a^{b^q}\rceil\qquad\dg_q\coloneq\Lg\zg_q=\dg\lg_1^{2\bg}\lg_q^{-2\bg}\qquad\zg_q\coloneq\lg_q^{-2\bg}\qquad\Lg\coloneq\dg\lg_1^{2\bg}, \xtag{DefParams}\label{eq:DefParams}\]
where
\bi
\item $\lceil x\rceil$ denotes the ceiling of $x$, i.e. the smallest integer $n\geq x$;
\item $\bg\in(0,\nicefrac13)$ and $b\in(1,\nicefrac32)$ control the Hölder exponent of the scheme and are required to satisfy
\[1<b<\xfr{1-\bg}{2\bg}. \xtag{bBeta}\label{eq:bBeta}\]
\item $a\dgreat1$ is sufficiently large to absorb various $q$-independent constants in the course of the proofs.
\ei
The parameter $\Lg$, and thus the distinction between $\dg_q$ and $\zg_q$, were absent in \cite{DRSz}. They are added here to make sure $\dg_1=\dg$, thus making \eqref{eq:U1TrEst} an $a$-independent estimate. Thus, in particular, we are allowed to bound $\Lg$ from below, since such a bound will be satisfied for $a$ large enough, but not to bound it from above, which would cause $\dg$ to depend on $a$.
\nipar With this choice of parameters, we must require the conditions
\begin{align*}
\Lg\geq{}&1 \xtag{EquestLg}\label{eq:EquestLg} \\
\xfr13>\bg>{}&\EXP+\eg', \xtag{EquestBgQgAgbar}\label{eq:EquestBgQgAgbar}
\end{align*}
for some positive $\eg'$. Condition \eqref{eq:EquestLg} merely requires $a$ to be sufficiently large.
\nipar The main convex integration step will consist in stating that, for a certain universal constant $M>1$, some sufficiently small $\ag,\gg>0$, and a sufficiently large $a\dgreat1$, if $(v_q,p_q,R_q)$ is a strong subsolution satisfying
\begin{align*}
\norm{\pint R_q}_0\leq{}&\Lg\varrho_q^{1+\gg} \xtag{GStrong}\label{eq:GStrong} \\
\norm{v_q}_{1+\ag}\leq{}&M\dg_q^{\xfr12}\lg_q^{1+\ag} \xtag{Gv}\label{eq:Gv} \\
\xfr34\dg_{q+2}\leq{}&\rg_q\leq\xfr72\dg_{q+1} \xtag{GTr}\label{eq:GTr} \\
\abs{\pd_t\rg_q}\leq{}&\rg_q\dg_q^{\xfr12}\lg_q \xtag{GTrDer}\label{eq:GTrDer} \\
\norm{v_q}_{\EXP+\eg}\leq{}& M\pa{1+\sum_{i=0}^q\lg_i^{\EXP+\eg-\bg}}, \xtag{GQg+Eg}\label{eq:GQg+Eg}
\end{align*}
where $\rg_q\coloneq\xfr13\opn{tr}R_q$, and $\varrho_q\coloneq\Lg^{-1}\rg_q$, then there exists a strong subsolution $(v_{q+1},p_{q+1},R_{q+1})$ satisfying the conditions \eqref{eq:GStrong}-\eqref{eq:GQg+Eg} with $q$ replaced by $q+1$ as well as the following additional estimate
\[\norm{v_{q+1}-v_q}_0+\lg_{q+1}\norm{v_{q+1}-v_q}_{H^{-1}}+\lg_{q+1}^{-1-\ag}\norm{v_{q+1}-v_q}_{1+\ag}\leq M\dg_{q+1}^{\xfr12}.\]
\clearpage\noindent
The proof consists of three steps:
\ben
\item A mollification step, moving from $(v_q,p_q,R_q)$ to $(v_{\ell_{q,i}},p_{\ell_{q,i}},R_{\ell_{q,i}})$, where the mollification parameter $\ell_{q,i}$ varies on suitably chosen subintervals, as required by the different orders of the upper and lower bounds on $\rg_q$ in \eqref{eq:GTr};
\item A gluing step, which goes from $(v_{\ell_{q,i}},p_{\ell_{q,i}},R_{\ell_{q,i}})$ to $(\lbar v_q,\lbar p_q,\lbar R_q)$;
\item A perturbation step going from $(\lbar v_q,\lbar p_q,\lbar R_q)$ to $(v_{q+1},p_{q+1},R_{q+1})$.
\een
The change in condition \eqref{eq:GStrong} with respect to \cite{DRSz} was made in order to prevent the new definition of $\dg_q$ from causing bounds of the form $\Lg^A\leq1$, with $A>0$, to appear in the proofs. Condition \eqref{eq:GQg+Eg} was added in order to control the new trace terms.
\nipar Section \ref{Gluing} proves the mollification and gluing steps, and Section \ref{MPS} addresses the perturbation step. The gluing step was introduced in \cite{I} to ensure $\pint{\lbar R}_q$ is supported in pairwise disjoint time intervals. This allows us to construct the perturbation as $w=\sum_i(w_{o,i}+w_{c,i})$, where the $w_{o,i}$ are Mikado flows with pairwise disjoint supports and $\opn{supp}w_{c,i}\sbse\opn{supp}w_{o,i}$, thus preventing $w\otimes w$ from containing ``mixed terms'' $w_{o,i}\otimes w_{o,j}$ with $i\neq j$, which are harder to deal with.
\nipar Fixing $\ag>0,\gg>0$, we also define
\begin{align*}
\ell_q\coloneq{}&\xfr{\zg_{q+2}^{\xfr{1+\gg}{2}}}{\zg_q^{\xfr12}\lg_q\lg_{q+1}^{2\ag}}=\xfr{\dg_{q+2}^{\xfr{1+\gg}{2}}\Lg^{-\xfr\gg2}}{\dg_q^{\xfr12}\lg_q\lg_{q+1}^{2\ag}} \xtag{DefEll}\label{eq:DefEll} \\
\tg_q\coloneq{}&\xfr{\ell_q^{4\ag}}{\dg_q^{\xfr12}\lg_q}. \xtag{DefTau}\label{eq:DefTau}
\end{align*}\vsp{-1em}
\xbegin{oss}[Homogeneity in \texorpdfstring{$\Lg$ of $\ell_q,\tg_q$}{Lg of ellq,tgq}][thm:oss:EllqTgqLg]
$\ell_q$, as well as the $\ell_{q,i}$ defined in Section \ref{Gluing}, are 0-homogeneous in $\Lg$, whereas $\tg_q$ is ${}^{1\!}/_2$-homogeneous. The last property allows us to cancel the $\Lg^{\nicefrac12}$ factors we will see appearing in the course of the proof.
\xend{oss}
We also assume
\[\xfr{\dg_{q+1}^{\xfr12}\dg_q^{\xfr12}\lg_q}{\lg_{q+1}^{1-15\ag-\bg\gg}}\leq\dg_{q+2}, \xtag{DLRel}\label{eq:DLRel}\]
which can be achieved if $a$ is sufficiently large assuming $(15\ag+\bg\gg)b<(b-1)(1-\bg-2b\bg)$. Moreover, we assume
\begin{align*}
\lg_{q+1}^{-1}\leq{}&\ell_q\leq\lg_q^{-1}. \xtag{LambdaEll}\label{eq:LambdaEll} \\
\intertext{The right inequality in \eqref{eq:LambdaEll} is evident from the definition. The left inequality can be reduced to $-b<\bg b^2(1+\gg)+\bg-1-2b\ag$, which can easily be verified for $\ag=0=\gg$, and thus also for $\ag,\gg$ sufficiently small. We will in fact need the following sharper bound:}
\lg_{q+1}^{1-\lbar N}\leq{}&\ell_q^{\lbar N+1},\xtag{LambdaEll.2}\label{eq:LambdaEll.2}
\end{align*}
which can be achieved by imposing the following condition:
\[\lbar N[(b-1)(1-\bg(b+1)\!)-\gg\bg b^2-2\ag b]>1+b+(1+\gg)\bg b^2+2\ag b-\bg. \xtag{LambdaEllq2log}\label{eq:LambdaEllq2log}\]
The above conditions can be obtained by choosing, in this order
\bi
\item $b,\bg$ as in \eqref{eq:bBeta}, so that in particular $\bg(1+b)<1$;
\item $0<\ag,\gg$ sufficiently small depending on $b,\bg$;
\item $\lbar N\in\N$ sufficiently large depending on $b,\bg,\ag,\gg$ so as to get \eqref{eq:LambdaEllq2log}.
\ei
One last notational remark: $A\lsim B$ (resp. $A\gsim B$) will mean $A\leq C(b,\bg,\ag,\gg,M)B$ (resp. $A\geq C(b,\bg,\ag,\gg,M)B$), or $C(N,b,\bg,\ag,\gg,M)$ if norms depending on $N$ are involved (e.g. $\s{C}^{N+1+\ag}$-norms). $A\sim B$ will mean $A\lsim B$ and $A\gsim B$. Note that $C$ does not depend on $a\dgreat1$.

\sect{Main iterative propositions}\label{Propos}
In this section, we state the main propositions which allow us to pass from one kind of subsolution to another one, which is closer to the notion of solution. Both the below statements use the parameters $\dg_q,\lg_q,\Lg$ defined in the previous section. The combination of these propositions leads to our main theorem, as illustrated in Section \ref{FinArg}.
\nipar In the first proposition it is shown that a smooth strict subsolution can be approximated with an adapted subsolution.
\xbegin{propo}[From strict to adapted subsolutions][thm:propo:S2A]
Let $\strsub$ be a smooth strict subsolution on $[0,T]$. Then, for any $\EXP<\adH<{}^1/_3$, $\ng>\xfr{1-3\adH}{2\adH}$, and $\dg,\sg>0$, there exist $\gg,\Wg>0$ and a $\s{C}^{\adH}$-adapted subsolution $\adsub$ with parameters $\gg,\Wg,\ng$ such that $\ad\rg\leq\xfr54\dg$ and, for all $t\in[0,T]$
\begin{align*}
\xints{\T^3}{}\2\pa[b]{\abs{\ad v}^2+\opn{tr}\ad R}\diff x={}&\xints{\T^3}{}\pa[b]{\abs{\str v}^2+\opn{tr}\str R}\diff x \xtag{S2AIC}\label{eq:S2AIC} \\
\norm{\str v-\ad v_{\s{C}^0}\lsim{}}&1+\dg^{\xfr12} \xtag{S2AvC0}\label{eq:S2AvC0} \\
\norm{\str v-\ad v}_{H^{-1}}<{}&\sg \xtag{S2Av}\label{eq:S2Av}
\end{align*}
Moreover, if we define
\[\ad{\s{T}}(t)\coloneq\xints0t\xints{\T^3}{}\pa[B]{\abs{\xfrl*\str v}^2-\abs{\xfrl*\ad v}^2}\diff x\diff s, \xtag{S2AE}\label{eq:S2AE}\]
we have the bound
\[\abs{\pd_t\ad{\s{T}}}\lsim\sum_q\Lg\lg_q^{\EXP+\eg-\bg}. \xtag{S2ADerT}\label{eq:S2ADerT}\]
The $q=0$ term of this sum is the largest, and is $\dg\lg_1^{2\bg}\lg_0^{\EXP+\eg-\bg}$, which is $a$-increasing. Since we can see that $a\to\8$ for $\dg\to0$, for any $\hg>0$, it can only be ensured that
\[\abs{\!\!\ad{\,\,\s{T}}}(t)\leq\hg\qquad t\in[0,\ad T(\hg,\dg,a)],\]
where $\ad T(\hg,\dg,a)\sim\hg\dg^{-1}\lg_1^{-2\bg}\lg_0^{\bg-\EXP-\eg}\to0$ if $a\to\8$ or $\hg\to0$.
\xend{propo}
The proof will be given in Section \ref{S2AProof}.
\nipar Next, we show that adapted subsolutions can be approximated by weak solutions with the same initial data.
\xbegin{propo}[From adapted subsolutions to weak solutions][thm:propo:A2S]
Let $\EXP<\solH<\adH<\xfr13$, $\gg>0$, and $\ng>0$ with
\[\xfr{1-3\adH}{2\adH}<\ng<\xfr{1-3\solH}{2\solH}. \xtag{A2SBetasNu}\label{eq:A2SBetasNu}\]
The following holds for all $\dg<1$.
\nipar If $\adsub$ is a $\s{C}^{\adH}$-adapted subsolution with parameters $\gg,\Wg,\ng$ and $\ad\rg\leq\xfr52\dg$, then, for all $\sg>0$, there exists a $\s{C}^\solH$ weak solution $v$ of \eqref{eq:FNS} with initial datum
\[\sol v(\per,0)=\ad v(\per,0) \xtag{A2SIV}\label{eq:A2SIV}\]
and such that, for all $t\in[0,T]$
\begin{align*}
\xints{\T^3}{}\2\abs{\sol v}^2\diff x={}&\xints{\T^3}{}\pa[b]{\abs{\ad v}^2+\opn{tr}\ad R}\diff x \xtag{A2SIC}\label{eq:A2SIC} \\
\norm{\sol v-\ad v}_{\s{C}^0}\lsim{}&\dg^{\xfr12} \xtag{A2SvC0}\label{eq:A2SvC0} \\
\norm{\sol v-\ad v}_{H^{-1}}<{}&\sg \xtag{A2Sv}\label{eq:A2Sv}
\end{align*}
Moreover, if we define
\[\sol{\s{T}}(t)\coloneq\xints0t\xints{\T^3}{}\pa[B]{\abs{\xfrl*\hat v}^2-\abs{\xfrl*v}^2}\diff x\diff s, \xtag{A2SE}\label{eq:A2SE}\]
once again we can see that
\[\abs{\pd_t\sol{\s{T}}}\lsim\sum\Lg\lg_q^{\EXP+\eg-\bg}, \xtag{A2SDerT}\label{eq:A2SDerT}\]
so that, like in the previous proposition, for any $\hg>0$, it can be ensured that
\[\abs{\sol{\s{T}}(t)}\leq\hg\qquad\VA t\in[0,\sol T(\hg,\dg,a)],\]
where $\sol T(\hg,\dg,a)\sim\hg\dg^{-1}\lg_1^{-2\bg}\lg_0^{\bg-\EXP-\eg}\to0$ if $\hg\to0$ or $a\to\8$.
\nipar Finally, consider the family of strong subsolutions $(\ad v,\ad p,\ad R+{}^e/_3\opn{Id})$, where $e:[0,T]\to\R$ satisfies the following conditions:
\begin{align*}
e(t)\leq{}&\xfr52\dg-\ad\rg(t) \xtag{CondE1}\label{eq:CondE1} \\
\abs{\pd_te}\leq{}&\rad{\dg_0}\lg_0e. \xtag{CondE2}\label{eq:CondE2} \\
e\geq{}&0. \xtag{CondE3}\label{eq:CondE3}
\end{align*}
This family can be used to yield infinitely many distinct weak solutions with the same initial data as $\solu$.
\xend{propo}
The proof will be given in Section \ref{A2SProof}.
\nipar \kcref{thm:propo:A2S} allows us to prove the following wildness criterion.
\xbegin{cor}[Wildness criterion][thm:cor:WildnessCrit]
Let $\EXP<\bg<\xfr13$ and $(\hat v,\hat p,\hat R)$ be a $\s{C}^\bg$-adapted subsolution such that $\rg\leq\xfr52\dg$ for some small $\dg>0$ and $\rg^{-1}|\pd_t\rg|\leq M\dg$ for some suitably large $M>0$. Assume that the following admissibility condition is satisfied for all $t\in[0,t_a]$ for some sufficiently small $t_a$:
\[\xfr12\xints{\T^3}{}\pa[b]{\abs{\hat v}^2(x,t)+\opn{tr}\hat R(x,t)}\diff x+\xints0t\xints{\T^3}{}\2\abs{\xfrl*\hat v}^2(x,s)\diff x\diff s\leq\xfr12\xints{\T^3}{}\pa[b]{\abs{\hat v}^2(x,0)+\opn{tr}\hat R(x,0)}\diff x, \xtag{WCAdm}\label{eq:WCAdm}\]
with a strict inequality for at least some $t\in[0,t_a]$. Then $\hat v(x,0)\in W_{\EXP,\bg-\eg,t_a}$ for any $\eg>0$.
\xend{cor}
The existence of infinitely many $\s{C}^{\bg-\eg}$ weak solutions with $\hat v(x,0)$ as their initial datum is a consequence of \kcref{thm:propo:A2S} above. The admissibility of those solutions follows from \eqref{eq:WCAdm} as shown in the next section, where it is also seen that the strictness of \eqref{eq:WCAdm} for at least some $t$ is vital to the admissibility of the solutions.
\nipar We shall henceforth adopt the following notational convention, already applied in the statements of the propositions:
\bi
\item $\strsub$ will always denote strict subsolutions;
\item In Sections \ref{FinArg}-\ref{S2S}, $\strngsub$ will always denote strong subsolutions; in Sections \ref{Gluing}-\ref{MPS}, all subsolutions will be strong, and in Sections \ref{S2AProof}-\ref{A2SProof}, the subscripts will mark strong subsolutions;
\item $\adsub$ will always denote adapted subsolutions, and $\ad\bg$ will be the Hölder regularity of adapted subsolutions;
\item $\solu$ will always denote (weak) solutions.
\ei
\clearpage

\sect{Proof of the existence theorem}\label{FinArg}
We start by recalling the following classical result.
\xbegin{teor}[Existence of Leray solutions][thm:teor:ExistLeraySol]
For any $w\in L^2(\T^3)$ with $\opn{div}w=0$ and every $\EXP\in(0,1)$ there is a weak solution $\sol v\in L^\8(\R^+,L^2(\T^3)\!)\cap L^2(\R^+,H^\EXP(\T^3)\!)$ of \eqref{eq:FNS} such that $\sol v(\per,0)=w$ and
\[\xfr12\xints{\T^3}{}\2\abs{\sol v}^2(x,t)\diff x+\xints0t\xints{\T^3}{}\2\abs{\xfrl*\sol v}^2(x,s)\diff x\diff s\leq\xfr12\xints{\T^3}{}\2\abs w^2(x)\diff x\qquad\VA t\geq0. \xtag{EI}\label{eq:energy}\]
In fact, the following form of energy inequality also holds:
\[\xfr12\xints{\T^3}{}\2\abs{\sol v}^2(x,t)\diff x+\xints st\xints{\T^3}{}\2\abs{\xfrl*\sol v}^2(x,\tg)\diff x\diff\tg\leq\xfr12\xints{\T^3}{}\2\abs{\sol v}^2(x,s)\diff x\qquad\text{a.e. }s,\VA t>s.\]
\xend{teor}
Recalling the \kcref{thm:defi:Strict} of subsolutions and strict subsolutions, one can prove the following existence result.
\xbegin{lemma}[Existence of strict subsolutions][thm:lemma:ExistStrictSubsol]
Let $w\in L^2(\T^3)$ with $\opn{div}w=0$. For any $\dg>0$ there exists a smooth strict subsolution $\strsub$ defined on $[0,T)$ such that
\begin{align*}
\norm{\str v|_{t=0}-w}_{L^2(\T^3)}\leq{}&\dg, \xtag{ExStrictSubsolIV}\label{eq:ExStrictSubsolIV} \\
\intertext{and for all $t\in[0,T]$}
\xfr12\xints{\T^3}{}\pa[b]{\abs{\str v}^2(x,t)+\opn{tr}\str R(x,t)}\diff x+\xints0t\xints{\T^3}{}\2\abs{\xfrl*\str v}^2(x,s)\diff x\diff s\leq{}&\xfr12\xints{\T^3}{}\abs{w}^2(x)\diff x+\dg. \xtag{ExStrictSubsolEI}\label{eq:ExStrictSubsolEI}
\end{align*}
\xend{lemma}
The proof is an adaptation of the one of \cite[Lemma 6.8, p. 38]{S}, and is reported in Appendix \ref{ExistStrictSubsol}.
\nipar The proof of the main result then follows the steps of \cite[Section 4]{DSz}, taking care of the additional dissipation by suitably increasing the energy of the starting strict subsolution.
\nipar\begin{qeddim*}[\kcref{thm:teor:TimeExist}]
We choose $\hg>0,\EXP<\bg<\bg',w\in L^2$ with $\opn{div}w=0$. Using the above result, we obtain a smooth strict subsolution $(\str v',\str p',\str R')$ on $[0,T]$ such that \eqref{eq:ExStrictSubsolIV}-\eqref{eq:ExStrictSubsolEI} hold for some $\dg>0$ which we will fix later. We now note that adding a smoothly time-dependent non-negative multiple of the identity to $\str R'$ does not change the fact that $(\str v',\str p',\str R')$ is a smooth strict subsolution. We may thus substitute our strict subsolution with $\strsub\coloneq(\str v',\str p',\str R'+2(3|\T^3|)^{-1}e_K(t)\opn{Id})$, where $K$ is a constant to be specified later in this proof, $0\leq e_K(t)\leq({}^\dg/_2-Kt)^+$, and $e_K(0)={}^\dg/_2$. Combining the choice of $e_K$ with \eqref{eq:ExStrictSubsolEI}, we obtain the following relations for ${}^\dg/_2-Kt>0$:
\begin{align*}
\xfr12\xints{\T^3}{}\pa[b]{\abs{\str v}^2(x,0)+\opn{tr}\str R(x,0)}\diff x={}&\xfr12\xints{\T^3}{}\abs{w}^2(x)\diff x+\xfr32\dg \xtag{ThmEI1}\label{eq:ThmEI1} \\
\xfr12\xints{\T^3}{}\pa[b]{\abs{\str v}^2(x,t)+\opn{tr}\str R(x,t)}\diff x+\xints0t\xints{\T^3}{}\2\abs{\xfrl*\str v}^2(x,s)\diff x\diff s\leq{}&\xfr12\xints{\T^3}{}\2\abs{w}^2(x)\diff x+\xfr32\dg-Kt. \xtag{ThmEI2}\label{eq:ThmEI2}
\end{align*}
Indeed, passing from $(\str v',\str p',\str R')$ to $\strsub$ adds a term $e_K$ to the left-hand side, since $\opn{tr}\str R=\opn{tr}\str R'+2|\T^3|^{-1}e_K$. Now let $\str v_0$ be the initial datum of $\str v$, and note that
\begin{align*}
\xints{\T^3}{}\2\opn{tr}\str R(x,0)\diff x={}&\norm{w}_{L^2}^2-\norm{\str v_0}_{L^2}^2+3\dg\leq\norm{w-\str v_0}_{L^2}(\norm{w}_{L^2}+\norm{\str v_0}_{L^2})+3\dg\leq{} \\
{}\leq{}&\dg(2\norm{w}_{L^2}+\dg)+3\dg\leq C(w)\dg. \xtag{ThmIntTrEst}\label{eq:ThmIntTrEst}
\end{align*}
Using \kcref{thm:propo:S2A} and \kcref{thm:propo:A2S}, we can produce a $\s{C}^{\adH}$-adapted subsolution $\adsub$ and a $\s{C}^{\solH}$ weak solution $\solu$, satisfying the integral equalities \eqref{eq:S2AIC} and \eqref{eq:A2SIC} and the $H^{-1}$ estimates \eqref{eq:S2Av} and \eqref{eq:A2Sv}, and the functions $\ad{\s{T}},\sol{\s{T}}$ of \eqref{eq:S2AE} and \eqref{eq:A2SE}. Recall that we have that
\[\xints{\T^3}{}\2(|\str v|^2+\opn{tr}\str R)(x,t)\diff x=\xints{\T^3}{}\2(|\ad v|^2+\opn{tr}\ad R)(x,t)\diff x=\xints{\T^3}{}\2\abs{\sol v(x,t)}^2\diff x, \xtag{ThmEnEq}\label{eq:ThmEnEq}\]
and thus
\[\norm{\sol v(t)}_2^2-\norm{\str v(t)}_2^2=\int\opn{tr}\str R(x,t)\diff x. \xtag{ThmDiffL2Norm}\label{eq:ThmDiffL2Norm}\]
Call $\sol v_0$ the initial datum of $\sol v$ and of $\ad v$, and note that, by \eqref{eq:S2Av}, \eqref{eq:ThmDiffL2Norm} and \eqref{eq:ThmIntTrEst}, we have that
\[\norm{\sol v_0-\str v_0}_2^2=\norm{\sol v_0}_2^2-\norm{\str v_0}_2^2-2\per\e\xints{\T^3}{}\2\str v_0\per(\sol v_0-\str v_0)\diff x\leq C(w)\dg+2\sg\norm{\str v_0}_{H^1}.\]
Thus, we first choose $\dg$ sufficiently small so that $C(w)\dg<\xfr{\hg^2}{2}$ and obtain $\str v$, then we fix $\sg<\xfr{\hg^2}{4\|\str v_0\|_{H^1}}$ and obtain $\ad v$, and finally we conclude that
\[\norm{\sol v_0-\str v_0}_2^2\leq\hg^2\implies\norm{\str v_0-\sol v_0}_{L^2}\leq\hg.\]
As for the admissibility condition, choosing $K$ so that $|\pd_t(\sol{\s{T}}+\ad{\s{T}})|\leq K-1$, as is made possible by \eqref{eq:S2ADerT} and \eqref{eq:A2SDerT}, we have that
\begin{multline*}
\xfr12\xints{\T^3}{}\2\abs{\sol v(x,t)}^2\diff x+\xints0t\xints{\T^3}{}\abs{\xfrl*\sol v}^2(x,s)\diff s\diff x \\
{}\overset{\text{\eqref{eq:ThmEnEq}}}{=}\xfr12\xints{\T^3}{}\pa[B]{\abs{\str v(x,t)}^2+\opn{tr}\str R(x,t)}\diff x-(\ad{\s{T}}+\sol{\s{T}})(t)+\xints0t\xints{\T^3}{}\2\abs{\xfrl*\str v}^2(x,s)\diff s\diff x \\
{}\overset{\text{\eqref{eq:ThmEI2}}}{\leq}\xints{\T^3}{}\2\xfr12\abs{w}^2(x)\diff x+\xfr32\dg-t\overset{\text{\eqref{eq:ThmEI1}}}{\leq}\xfr12\xints{\T^3}{}\pa[b]{\abs{\str v_0}^2+\opn{tr}\str R(x,0)}\diff x\overset{\text{\eqref{eq:ThmEnEq}}}{=}\xfr12\xints{\T^3}{}\2\abs{\sol v}^2(x,0)\diff x,
\end{multline*}
where the second-last inequality is strict for all $t\neq0$ where \eqref{eq:ThmEI2} is valid. This yields the energy inequality for $t$ sufficiently small. Since we can only estimate $|\pd_t\sol{\s{T}}+\pd_t\ad{\s{T}}|$ with a quantity which is potentially unbounded as $\dg\to0$ (as seen in \eqref{eq:S2ADerT} and \eqref{eq:A2SDerT}), and $C(w)\dg<\xfr{\hg^2}{2}$ implies $\dg\to0$ as $\hg\to0$, we conclude that our time $T(\hg)$ of guaranteed admissibility satisfies
\[\lim_{\hg\to0}T(\hg)=0.\]
So far, we have only obtained one solution for each $\hg$. Suppose that, from $\strsub$, we produced the adapted subsolution $\adsub$, and from there the solution $\solu$. As noted in \kcref{thm:propo:A2S}, considering
\[(\ad v',\ad p',\ad R')\coloneq\pa{\ad v,\ad p,\ad R+\xfr e3\opn{Id}}\]
with $e$ satisfying a suitable set of conditions, we can obtain more weak solutions and ensure these solutions are admissible up to $T(\hg)$. The required conditions are listed below.
\ben
\item The first condition ensures $\ad R'(\per,0)\equiv0$:
\[e(0)=0;\]
\item The second one ensures $\opn{tr}(\ad R')\geq0$:
\[e(t)\geq0;\]
it would be enough to require $\ad\rg'\geq|e|$, but we exclude $e<0$ for convenience (cfr. Step 3 of Section \ref{A2SProof});
\item The third one ensures the admissibility of the new solutions:
\begin{align*}
\xfr12|\T^3|e(t)\leq{}&\xfr12\xints{\T^3}{}\pa[b]{\abs{\ad v(x,0)}^2+\opn{tr}\ad R(x,0)-\abs{\ad v(x,t)}^2-\opn{tr}\ad R(x,t)}\diff x \\
&{}-\xints0t\xints{\T^3}{}\2\abs{\xfrl*\ad v(x,s)}^2\diff x\diff s-\ad{\s{T}}-Kt,
\end{align*}
$K$ being the same constant used to find $(v,p)$; since the right-hand side of the above inequality is strictly positive for all $t\neq0$ where $\solu$ is admissible, this condition is compatible with requiring that $e\geq0$ as done above;
\item The last conditions are \eqref{eq:CondE1}-\eqref{eq:CondE3}.
\een
This completes the proof.
\end{qeddim*}

\sect{From strict to strong subsolutions}\label{S2S}
\numberwithin{defi}{section}
\numberwithin{propo}{section}
We state here an analogue of \cite[Proposition 3.1]{DSz}.
\xbegin{propo}[][thm:propo:S2S]
Let $\strsub$ be a smooth solution of \eqref{eq:FNSR}, and $S\in\Cinf(\T^3\x[0,T];\s{S}^{3\x3}_+)$ be a smooth positive-definite matrix field. Fix $\lbar\ag\in(0,1)$ and $\eg>0$. Then for any $\lg>1$ there exists a smooth solution $\strngsub$ of \eqref{eq:FNSR} with
\begin{align*}
\strngsub={}&\strsub\quad\text{for }t\nin\opn{supp}\opn{tr}S, \xtag{S2SEq}\label{eq:S2SEq} \\
\xints{\T^3}{}\pa[b]{|\strng v|^2+\opn{tr}\strng R}(x,t)\diff x={}&\xints{\T^3}{}\pa[b]{|\str v|^2+\opn{tr}\str R}(x,t)\diff x\qquad\VA t\in[0,T], \xtag{S2SIC}\label{eq:S2SIC}
\end{align*}
and the following estimates hold
\begin{align*}
\norm{\strng v-\str v}_{H^{-1}}\leq{}&\xfr C\lg \xtag{S2SEst1}\label{eq:S2SEst1} \\
\norm{\strng v}_k\leq{}&C\lg^k\qquad k=1,2 \xtag{S2SEst2}\label{eq:S2SEst2} \\
\norm{\str R-\strng R-S}_N\leq{}&\xfr{C}{\lg^{1-2\EXP-\lbar\ag-N}} \xtag{S2SEst3}\label{eq:S2SEst3}
\end{align*}
Moreover, $\opn{tr}(\str R(x,t)-\strng R(x,t)-S(x,t)\!)\eqcolon\mg(t)$ is a function of $t$ only and satisfies
\[\abs{\mg'}(t)\leq C\lg^{\lbar\ag}. \xtag{S2STrEst}\label{eq:S2STrEst}\]
The constant $C\geq1$ above depends on $(v,p,R),S$ and $\lbar\ag$, but not on $\lg$. Finally, defining
\[\strng{\s{T}}(t)\coloneq\xints0t\xints{\T^3}{}\pa[B]{\abs{\xfrl*\strng v(x,s)}^2-\abs{\xfrl*\str v(x,s)}^2}\diff x\diff s,\]
we have that
\[\abs{\pd_t\strng{\s{T}}(t)}\leq C\lg^{2(\EXP+\eg)}. \xtag{S2SHypoTrace}\label{eq:S2SHypoTrace}\]
\xend{propo}
\begin{qeddim}
Define the inverse flow of $\str v$, $\Fg:\T^3\x[0,T]\to\T^3$, as the solution of
\[\begin{sistema}
\pd_t\Fg(x,t)+(\str v\per\grad)\Fg(x,t)=0 \\
\Fg(x,0)=x\qquad x\in\T^3
\end{sistema},\]
and set
\[\lbar R(x,t)=\fD\Fg(x,t)S(x,t)\fD^T\Fg(x,t).\]
Observe that $\lbar R$ is defined on the compact set $\T^3\x[0,T]$ and, being continuous, has a compact image $\s{N}_0\coloneq\lbar R(\T^3\x[0,T])\sbs\s{S}_+^{3\x3}$.
\nipar By \kcref{thm:lemma:MikFlows} there exists a smooth vector field $W:\s{N}_0\x\T^3\to\T^3$ satisfying the differential equations \eqref{eq:MikDiff} and the integral equations \eqref{eq:MikInt}. Define
\begin{align*}
w_o(x,t)={}&\fD\Fg^{-1}W(\lbar R,\lg\Fg(x,t)\!) \\
w_c(x,t)={}&\xfr1\lg\opn{curl}(\fD^T\Fg U(\lbar R,\lg\Fg(x,t)\!)\!)-w_o,
\end{align*}
where $U=U(R,\jg)$ is defined as in \eqref{eq:AntiCurlW} and thus satisfies $\opn{curl}_\jg U=W$. Moreover, set
\[\strng v\coloneq\str v+w_o+w_c\qquad\strng p=\str p+\lbar p\qquad\strng R\coloneq\str R-S-\pint{\s{E}}\stp1-\s{E}\stp2,\]
where
\begin{align*}
\lbar p\coloneq{}&-\xfr13(w_c\per\strng v+w_o\per w_c) \\
\pint{\s{E}}\stp1\coloneq{}&\s{R}(F)+(w_c\otimes\strng v+w_o\otimes w_c+\lbar p\opn{Id}) \\
F\coloneq{}&\opn{div}(w_o\otimes w_o-S)+(\pd_t+\str v\per\grad)w_o \\
&{}+[(w_o+w_c)\per\grad]\str v+\pd_tw_c+\xfrl(w_o+w_c) \\
\s{E}\stp2\coloneq{}&\xfr13\per\!\fint_{\T^3}\pa[B]{\abs{\strng v}^2-\abs{\str v}^2-\opn{tr}S}\diff x\per\opn{Id},
\end{align*}
with $\s{R}$ defined as in \eqref{eq:DefR}. By construction, the relation \eqref{eq:S2SIC} holds, $\pint{\s{E}}\stp1$ is traceless, $\s{E}\stp2$ is only $t$-dependent, and $\strngsub$ solves \eqref{eq:FNSR}. To verify this last claim, we can see that
\begin{align*}
\opn{div}\pint{\s{E}}\stp1={}&\opn{div}(\strng v\otimes\strng v-\str v\otimes\str v-S+\lbar p\opn{Id})+\pd_t(\strng v-\str v)+\xfrl(w_o+w_c) \\
{}={}&\pd_t\strng v+\opn{div}(\strng v\otimes\strng v-S+\pint R)+\grad\strng p+\xfrl\strng v.
\end{align*}
We call $w\coloneq w_o+w_c=\strng v-\str v$. Recall that
\begin{align*}
W(R,\jg)={}&\sum_{k\neq0}a_k(R)A_ke^{ik\per\jg} \\
U(R,\jg)={}&\sum_ka_k(R)\xfr{ik\x A_k}{|k|^2}e^{ik\per\jg},
\end{align*}
which are respectively \eqref{eq:MikFS} and \eqref{eq:AntiCurlW}, with the $a_k$ satisfying \eqref{eq:Mikak}. This allows us to decompose
\begin{align*}
w_o={}&\sum_{k\neq0}\fD\Fg^{-1}a_k(\fD\Fg S\fD^T\Fg)A_ke^{ik\per\lg\Fg}=\sum_{k\neq0}b_ke^{ik\per\lg\Fg} \xtag{Decwo}\label{eq:Decwo} \\
w_c={}&\xfr i\lg\sum_{k\neq0}\grad(a_k(\fD\Fg S\fD^T\Fg)\!)\x\xfr{\fD^T\Fg(k\x A_k)}{|k|^2}e^{ik\per\lg\Fg}=\sum_{k\neq0}\xfr{c_k}{\lg}e^{ik\per\lg\Fg}. \xtag{Decwc}\label{eq:Decwc}
\end{align*}
The estimate \eqref{eq:S2SEst3} is deduced by combining arguments from \cite{DSz} with estimates for $\s{E}\stp2$ and $\s{R}(\xfrl w)$. $\s{E}\stp2$ is estimated in a similar fashion to how we estimate $\pint{\s{E}}\stp1$ below. To estimate $\s{R}(\xfrl w)$, using the fact that $[\s{R},\xfrl]=0$ and \eqref{eq:FrLH}, we see that
\begin{align*}
\norm{\s{R}(\xfrl w)}_0\lsim\norm{\s{R}w}_{2\EXP+\lbar\ag}\lsim{}&\sum_k\bigg(\xfr{\norm{b_k+\lg^{-1}c_k}_0}{|\lg k|^{1-\lbar\ag-2\EXP}}+\xfr{\norm{b_k+\lg^{-1}c_k}_{N+2\EXP+\lbar\ag}}{|\lg k|^{N-2\EXP-\lbar\ag}} \\
&{\hsp{.5cm}}+\xfr{\norm{b_k+\lg^{-1}c_k}_0\norm{\Fg}_{N+2\EXP+\lbar\ag}}{|\lg k|^{N-2\EXP-\lbar\ag}}\bigg) \\
{}\lsim{}&\lg^{\lbar\ag+2\EXP-1}\per\!\sum_k\bigg(\xfr{1}{|k|^{7-\lbar\ag-2\EXP}}+\xfr{1+C_\Fg(N,\ag,\EXP)}{\lg^{N-1}|k|^{N+6-2\EXP-\lbar\ag}}\bigg),
\end{align*}
where we used \eqref{eq:Mikak} to get the extra $|k|^{-6}$ in each term, and the boundedness of $\Fg$ to get the $C_\Fg(N,\ag,\EXP)$.
\nipar Concerning \eqref{eq:S2SEst2}, the smoothness of $\Fg,S$ combined with \eqref{eq:Mikak} gives us
\[\max\{\|c_k\|_N,\|b_k\|_N\}\lsim|k|^{-m}, \sxtag{S2SEstCoeff}\]
for all integers $m>0$, where the $b_k$ and $c_k$ are as in the decompositions of $w_o,w_c$ above. This easily allows us to conclude that
\[\norm{w}_N\lsim\lg^N,\]
since differentiating the exponential gives us a factor of $\lg$ for each derivative. We then note that
\[\norm{\strng v}_N\leq\norm{\str v}_N+\norm{w}_N\lsim1+\lg^N\lsim\lg^N,\]
where the second step used the smoothness of $\str v$. For $N=1,2$ the above reduces to \eqref{eq:S2SEst2}.
\nipar Estimate \eqref{eq:S2SEst1} is proved separately for $w_c$ and $w_o$. The former is straightforward, since $w_c$ is already of order $\lg^{-1}$ in $L^2\hra H^{-1}$ thanks to \eqref{eq:Decwc}, and \eqref{eq:S2SEst1} is an $H^{-1}$ estimate. For the latter, the main idea is to write
\[e^{i\lg\Fg\per k}=\xfr{\pd_je^{i\lg\Fg\per k}}{i\lg\pd_j\Fg\per k},\]
and integrating by parts. One must then make sure that $j=j(x)$ is chosen in such a way that the denominator is bounded from below. This can be done by using the fact that $\Fg$ is a diffeomorphism, which implies $|\fD\Fg|$ is bounded below.
\nipar To continue, we note that $\mg=\opn{tr}\s{E}\stp2=\fint|\strng v|^2-|\str v|^2-\opn{tr}S\diff x$. $\str v$ and $\opn{tr}S$ are both smooth, so they are bounded. In order to estimate $\int|\strng v(t)|^2$, note that the following energy identity for $\strng v$ follows from \eqref{eq:FNSR}:
\[\pd_t\xfr12|\strng v|^2+\opn{div}\pa{\strng v\pa{\xfr{|\strng v|^2}{2}+\tilde p}\e}+\strng v\per\xfrl\str v=-\strng v\per\opn{div}(\pint{\strng R}-\s{R}(\!\xfrl w)\!).\]
Moreover, with the arguments used to estimate $R_{11}$ in \cite[pp. 18-20]{DSz}, we conclude that
\[\norm{\pint{\s{E}}\stp1-\s{R}(\!\xfrl w)}_0\lsim\lg^{\lbar\ag-1}.\]
These two bounds, by integrating in $x$ and using \eqref{eq:S2SEst2}, yield
\[\abs{\Der{t}\fint\xfr12|\strng v|^2\diff x}\lsim\fint|\strng v|\abs{\xfrl\str v}+|\grad\strng v|\abs{\pint{\strng R}-\s{R}(\!\xfrl w)}\diff x\leq C(1+\lg^{\lbar\ag}).\]
Thus, the estimate \eqref{eq:S2STrEst} is proved.
\nipar The last thing left is to estimate $|\pd_t\strng{\s{T}}|$. Combining some simple calculations with \eqref{eq:S2SEst2}, \eqref{eq:FrLH}, \eqref{eq:FrLInt}, and \eqref{eq:HInterp}, we obtain that
\begin{align*}
\abs{\pd_t\strng{\s{T}}}={}&\abs{\xints{\T^3}{}\2\abs{\xfrl*\strng v}^2-\abs{\xfrl*\str v}^2\diff x\diff s}=\abs{\xints{\T^3}{}\2[2\xfrl*\str v+\xfrl*w]\per\xfrl*w\diff x\diff s} \\
{}\lsim{}&\abs{\xints{\T^3}{}\22\xfrl*\str v\per\xfrl*w\diff x\diff s}+[w]_{\EXP+\eg}^2\diff s\lsim\underbrace{\norm{\str v}_{\EXP+\eg}\norm{w}_{\EXP+\eg}}_{\eqcolon I}+\lg^{2(\EXP+\eg)}.
\end{align*}
The velocity $\str v$ is bounded by smoothness, so $I\leq K(\str v,\EXP,\eg)\lg^{\EXP+\eg}$. Since $\lg>1$, this yields \eqref{eq:S2SHypoTrace}, thus concluding the proof.
\end{qeddim}
\kcref{thm:propo:S2S} will be applied in the situation described by the following corollary.
\xbegin{cor}[Strict to strong][thm:cor:S2Suse1]
Let $\strsub$ be a smooth strict subsolution on $[0,T]$. There exist $\tilde\dg,\gg>0$ such that the following holds.
\nipar For any $0<\dg<\tilde\dg$, $\ag,\gg>0$ and $0<\eg<\bg-\EXP$ sufficiently small, there exists a smooth strong subsolution $\strngsub$ with $\strng R(x,t)=\strng\rg(t)\opn{Id}+\pint{\strng R}(x,t)$, and a ``dissipative trace term'' as isolated in \kcref{thm:propo:S2S}, i.e.
\[\strng{\s{T}}(t)\coloneq\xints0t\xints{\T^3}{}\pa[B]{\abs{\xfrl*\str v(x,s)}^2-\abs{\xfrl*\strng v(x,s)}^2}\diff x\diff s,\]
such that, for all $t\in[0,T]$
\begin{align*}
\xints{\T^3}{}\pa[b]{|\str v(x,t)|^2+\opn{tr}\str R(x,t)}\diff x={}&\xints{\T^3}{}\pa[b]{|\strng v(x,t)|^2+\opn{tr}\strng R(x,t)}\diff x \xtag{U1EC}\label{eq:U1EC} \\
\xfr34\dg\leq\strng\rg\leq{}&\xfr54\dg \xtag{U1TrEst}\label{eq:U1TrEst} \\
\abs{\pint{\strng R\,\,}\!\!}\leq{}&\Lg\strng\varrho^{1+\gg} \xtag{U1Strong}\label{eq:U1Strong} \\
\norm{\strng v-\str v}_{H^{-1}}\leq{}&\dg\lg_0^{-1} \xtag{U1vH-1}\label{eq:U1vH-1} \\
\norm{\strng v}_{1+\ag}\leq{}&\dg_0\lg_0^{1+\ag} \xtag{U1v}\label{eq:U1v} \\
\abs{\pd_t\strng\rg(t)}\leq{}&\dg\dg_0^{\xfr12}\lg_0 \xtag{U1TrDer}\label{eq:U1TrDer} \\
\abs{\pd_t\strng{\s{T}}(t)}\leq{}&\Lg^{\xfr12}\dg_0^{\xfr12}\lg_0^{\EXP+\eg} \xtag{U1HypoTrace}\label{eq:U1HypoTrace} \\
\norm{\strng v}_{\EXP+\eg}\leq{}&K(1+\dg_0^{\xfr12}\lg_0^{\EXP+\eg}). \xtag{U1Fracv}\label{eq:U1Fracv}
\end{align*}
where the constant $K$ depends on $\strsub$ and $\eg$, the parameters $\dg_q,\lg_q,\zg_q,\Lg$ are defined as in \eqref{eq:DefParams} with sufficiently large $a$, and $\ag$ is the small parameter from Section \ref{Equestimates}.
\xend{cor}
\begin{qeddim}[The proof will proceed by first reducing all the claims in the corollary to a series of conditions on $\lg$, and then, at the end, proving that all those conditions can be satisfied simultaneously. This is necessary because some of them are upper bounds on $\lg$, and some are lower bounds.]
Let
\[\tilde\dg\coloneq\xfr12\inf\br{\str R(x,t)\jg\per\jg:|\jg|=1,x\in\T^3,t\in[0,T]}.\]
Since $\str R$ is a smooth positive definite tensor on a compact set, $\tilde\dg>0$. Then $S\coloneq\str R-\dg\opn{Id}$ is positive definite for any $\dg<\tilde\dg$. We may in addition assume without loss of generality that $\dg\leq1$. We apply \kcref{thm:propo:S2S} with $\strsub,S$, and $\lbar\ag\in(0,1),\eg>0$ to be chosen below. This yields a smooth solution $\strngsub$ of \eqref{eq:FNSR} with properties \eqref{eq:S2SIC}, \eqref{eq:S2SEst1}-\eqref{eq:S2SEst3}, and \eqref{eq:S2STrEst}. We first note that \eqref{eq:U1EC} coincides with \eqref{eq:S2SIC}. Next, we observe that $\strng R-\str R+S=\strng R-\dg\opn{Id}$, so that, since $\mg(t)=\opn{tr}(\strng R-\str R+S)$ is a function of time only, the function
\[\strng\rg=\xfr13\opn{tr}(\strng R-\str R+S)+\dg \xtag{U1TrEq}\label{eq:U1TrEq}\]
is independent of $x$.
\nipar Let us now prove \eqref{eq:U1TrEst}. By the above and \eqref{eq:S2SEst3} for $N=0$, we have that
\[\abs{\strng\rg-\dg}=\xfr13\abs{\opn{tr}(\strng R-\str R+S)}\leq\norm{\strng R-\str R+S}_0\leq C\lg^{2\EXP+\lbar\ag-1}, \xtag{U1PTrEst}\label{eq:U1PTrEst}\]
We require now the following condition on $\lg$:
\[C\lg^{2\EXP-1+\lbar\ag}\leq\dg_0^{\xfr12}\lg_0^{2\EXP-1+\lbar\ag}. \sxtag{U1lg1}\]
Then we notice that, for $\gg$ sufficiently small and $a$ sufficiently large, we have that
\[\dg_q^{\xfr12}\lg_q^{2\EXP+\lbar\ag-1}\leq\xfr14\dg_{q+1}\zg_{q+1}^\gg. \xtag{U1DLRel}\label{eq:U1DLRel}\]
Indeed, rewriting the above in terms of $\lg_q$, it reads
\[\Lg^{\xfr12}\lg_q^{-\bg+2\EXP+\lbar\ag-1}\leq\Lg\lg_q^{-2b\bg(1+\gg)}.\]
Since $\Lg\geq1$ by \eqref{eq:EquestLg}, this reduces to showing that
\[-\bg+2\EXP+\lbar\ag-1<-2b\bg(1+\gg). \xtag{U1bbgqg}\label{eq:U1bbgqg}\]
and taking $a$ sufficiently large. In turn, \eqref{eq:U1bbgqg} can be proved using that, by assumption, $\EXP<\bg,2b\bg<1-\bg$ (see \eqref{eq:bBeta}), and taking $\lbar\ag,\gg$ sufficiently small. Thus \eqref{eq:U1DLRel} is proved.
\nipar Now from \eqref{eq:U1PTrEst}, \eqref{eq:U1lg1}, and \eqref{eq:U1DLRel} for $q=0$, it follows that
\[\abs{\strng\rg-\dg}\leq\xfr14\dg_1\zg_1^\gg=\xfr14\Lg^{-\gg}\dg^{1+\gg}\leq\xfr14\dg, \xtag{U1Blabber}\label{eq:U1Blabber}\]
where in the last inequality we used the fact that $\dg<1<\Lg$. We have thus proved \eqref{eq:U1TrEst}.
\nipar From this estimate we can in turn deduce \eqref{eq:U1Strong}. Indeed, since $\,\,\pint{\!\!\strng R}=\,\,\pint{\!\!\strng R}-\,\,\pint{\!\!\str R}+\pint S$, by chaining the inequalities \eqref{eq:S2SEst3} for $N=0$, \eqref{eq:U1lg1}, and \eqref{eq:U1DLRel} for $q=0$, we analogously deduce that:
\[\abs{\,\,\pint{\!\!\strng R}}\leq\xfr14\dg_1\zg_1^\gg=\xfr14\Lg^{-\gg}\dg^{1+\gg}\leq\pa{\xfr34}^{1+\gg}\dg^{1+\gg}\Lg^{-\gg}\leq\strng\rg^{1+\gg}\Lg^{-\gg}\leq\Lg\strng\varrho^{1+\gg}.\]
The bound \eqref{eq:U1vH-1} follows from \eqref{eq:S2SEst1} together with the following condition on $\lg$:
\[C\lg^{-1}\leq\dg\lg_0^{-1}. \sxtag{U1lg3}\]
To obtain \eqref{eq:U1v}, we first use standard interpolation estimates together with \eqref{eq:S2SEst2} to obtain that
\[\norm{\strng v}_{1+\ag}\leq C_I\norm{\strng v}_1^{1-\ag}\norm{\strng v}_2^\ag\leq C_IC\lg^{1+\ag}.\]
Therefore, \eqref{eq:U1v} reduces to the following condition on $\lg$:
\[CC_I\lg^{1+\ag}\leq\dg_0^{\xfr12}\lg_0^{1+\ag}, \sxtag{U1lg4}\]
The estimate \eqref{eq:U1TrDer} follows from \eqref{eq:S2STrEst} and \eqref{eq:U1TrEq}), giving
\[\abs{\pd_t\strng\rg}=\xfr13\abs{\pd_t\opn{tr}(\strng R-\str R+S)}\leq\xfr C3\lg^{\lbar\ag}.\]
Therefore, \eqref{eq:U1TrDer} amounts to
\[\xfr C3\lg^{\lbar\ag}\leq\dg\Lg^{\xfr12}\lg_0^{1-\bg}. \sxtag{U1lg5}\]
Since by \eqref{eq:S2SHypoTrace} one has that $|\pd_t\strng{\s{T}}|\leq C\lg^{2(\EXP+\eg)}$, to obtain \eqref{eq:U1HypoTrace} we require
\[C\lg^{2(\EXP+\eg)}\leq\Lg^{\xfr12}\dg_0^{\xfr12}\lg_0^{\EXP+\eg}. \sxtag{U1lg6}\]
Finally, to obtain \eqref{eq:U1Fracv}, we note that $\strng v$ is smooth and thus bounded by a constant $C_0$, so that, by interpolation and \eqref{eq:S2SEst2}, we have that
\[\norm{\strng v}_{\EXP+\eg}\leq C_I\norm{\strng v}_0^{1-\EXP-\eg}\norm{\strng v}_1^{\EXP+\eg}\leq C_IC_0^{1-\EXP-\eg}(C\lg)^{\EXP+\eg}.\]
Therefore, we will require
\[C_IC_0^{1-\EXP-\eg}(C\lg)^{\EXP+\eg}\leq\dg_0^{\xfr12}\lg_0^{\EXP+\eg}. \sxtag{U1lg7}\]
To conclude the proof of the corollary, we now show that, for suitable choices of $\dg,\gg,\lbar\ag$, there exists a $\lg$ satisfying conditions \eqref{eq:U1lg1}, \eqref{eq:U1lg3}, \eqref{eq:U1lg4}, \eqref{eq:U1lg5}, \eqref{eq:U1lg6}, and \eqref{eq:U1lg7}.
\nipar In particular, for fixed constants $\lbar C,\lbar K$ independent of the parameters $a,\dg,b,\lbar\ag,\bg$, the following conditions must be satisfied by $\lg$:
\begin{align*}
\lg\geq{}&\lbar K\dg_0^{\xfr12(2\EXP-1+\lbar\ag)^{-1}}\lg_0 \sxtag{U1lg1a}\\
\lg\geq{}&\lbar K\dg^{-1}\lg_0 \sxtag{U1lg3a}\\
\lg\leq{}&\lbar C\lg_0\dg_0^{\xfr{1}{2+2\ag}} \sxtag{U1lg4a}\\
\lg\leq{}&\lbar C\dg^{\xfr{1}{\lbar\ag}}\Lg^{\xfr{1}{2\lbar\ag}}\lg_0^{\xfr{1-\bg}{\lbar\ag}} \sxtag{U1lg5a}\\
\lg\leq{}&\lbar C\Lg^{\xfr{1}{4(\EXP+\eg)}}\dg_0^{\xfr{1}{4(\EXP+\eg)}}\lg_0^{\xfr12} \sxtag{U1lg6a}\\
\lg\leq{}&\lbar C\dg_0^{\xfr12(\EXP+\eg)^{-1}}\lg_0 \sxtag{U1lg7a}
\end{align*}
First we choose $\dg<1$, and $\lbar\ag,\gg$ sufficiently small, and then show that, for $a$ sufficiently large there exists a $\lg$ satisfying all the above inequalities.
\nipar First of all, notice that, since $\dg_0=\dg\lg_0^{2\bg(b-1)}\dgreat1$ if $\dg$ is fixed and $a$ is sufficiently large, then \eqref{eq:U1lg3a} implies \eqref{eq:U1lg1a}, and \eqref{eq:U1lg4a} implies \eqref{eq:U1lg7a} independently of the choice of $\ag>0$, since $\EXP+\eg<\bg<\xfr13<1+\ag$.
\nipar Hence, we are left with showing that \eqref{eq:U1lg3a} is compatible with \eqref{eq:U1lg4a}-\eqref{eq:U1lg6a}.
\nipar The compatibility of \eqref{eq:U1lg3a} and \eqref{eq:U1lg4a}, independently of $\ag>0$, is straightforward, since $\dg_0\dgreat1$ when $a$ is sufficiently large.
\nipar Inequality \eqref{eq:U1lg5a} does not contradict \eqref{eq:U1lg3a} provided we choose $\lbar\ag$ so small that $\xfr{1-\bg}{\lbar\ag}>1$, and then $a$ sufficiently large.
\nipar The compatibility of \eqref{eq:U1lg3a} with \eqref{eq:U1lg6a} rewrites as
\[\lg_0^{\xfr12}\leq\xfr{\lbar C\dg}{\lbar K}\Lg^{\xfr{1}{4(\EXP+\eg)}}\dg_0^{\xfr{1}{4(\EXP+\eg)}}\]
and, inserting the definitions of $\dg_0,\Lg,\lg_1$, as
\[\lg_0^{\xfr12}\leq\xfr{\lbar C}{\lbar K}\dg^{\xfr{1}{2(\EXP+\eg)}+1}\lg_0^{(2b-1)\pa{\xfr{\bg}{2(\EXP+\eg)}}}.\]
Hence the above reduces to showing that
\[\xfr12\leq\xfr{\bg(2b-1)}{2(\EXP+\eg)},\]
which holds since $b>1$ and $\EXP+\eg<\bg$.
\nipar The proof is thus complete.
\end{qeddim}

\sect{Localized gluing step}\label{Gluing}
\xbegin{defi}[Decomposing the time interval][thm:defi:Subintervs+]
Let $0\leq T_1<T_2\leq T$ such that $T_2-T_1>4\tg_q$. We define sequences of intervals $\{I_i\},\{J_i\}$ as follows. Let
\[t_i\coloneq i\tg_q\qquad I_i\coloneq\sq{t_i+\xfr13\tg_q,t_i+\xfr23\tg_q}\cap[0,T], \xtag{$t_i,I_i$i}\label{eq:DeftiIi}\]
and let
\[\ubar n\coloneq\case{
\min\br{i:t_i-\xfr23\tg_q\geq T_1} & T_1>0 \\
0 & T_1=0
} \qquad\qquad \lbar n\coloneq\max\br{i:t_i+\xfr23\tg_q\leq T_2}. \xtag{$\lbar n,\ubar n$}\label{eq:lbarnubarn}\]
Moreover, define
\[\begin{aligned}
J_i\coloneq{}&\pa{t_i-\xfr13\tg_q,t_i+\xfr13\tg_q}\cap[0,T]\qquad\ubar n\leq i\leq\lbar n \\
J_{\ubar n-1}\coloneq{}&\sqpa{0,t_{\ubar n}-\xfr23\tg_q} \qquad\qquad J_{\lbar n+1}\coloneq\pasq{t_{\lbar n}+\xfr23\tg_q,T}.
\end{aligned} \xtag{Ji}\label{eq:Ji}\]
\xend{defi}
These form a pairwise disjoint decomposition of $[0,T]$:
\[[0,T]=J_{\ubar n-1}\cup I_{\ubar n-1}\cup[J_{\ubar n}\cup\dotso\cup J_{\lbar n}]\cup I_{\lbar n}\cup J_{\lbar n+1}, \xtag{Decomp[0,T]}\label{eq:Decomp0T}\]
and
\[t_{\ubar n}<T_1+\xfr53\tg_q<T_2-\xfr53\tg_q<t_{\lbar n}. \xtag*\label{eq:6.4}\]
Moreover, if $T_1>0$, $\ubar n\geq1$, otherwise we have both that $\ubar n=0$ and that $J_{\ubar n-1}\cup I_{\ubar n-1}=\0$.
\nipar Given a subsolution $(v_q,p_q,R_q)$ with $\rg_q\coloneq\xfr13\opn{tr}R_q$ and $\varrho_q\coloneq\Lg^{-1}\rg_q$, we will define:
\[(\rg_{q,i},\varrho_{q,i},\ell_{q,i})\coloneq\pa{\rg_q(t_i),\varrho_q(t_i),\xfr{\varrho_{q,i}^{\xfr{1+\gg}{2}}}{\zg_q^{\xfr12}\lg_q^{1+\ag}\ell_q^{-\ag}}}, \sxtag{rgqiellqi}\]
where $\ag,\gg$ are the parameters of Section \ref{Equestimates}. Using \eqref{eq:GUTr} and assuming $a\dgreat1$ is sufficiently large (as in \eqref{eq:LambdaEll}, depending on $\ag,\gg,b$), we may ensure that
\[\lg_{q+1}^{-1}\leq\ell_q\leq\ell_{q,i}\leq\lg_q^{-1}. \xtag{$\lg,\ell_q,\ell_{q,i}$}\label{eq:lgellqellqi}\]
Since we will always be working with $\EXP<\bg$, recalling $\ell_q^{-1}\leq\lg_{q+1}$, and assuming $\eg\leq\ag$, we observe that
\[\tg_q\ell_{q,i}^{-2\EXP-\eg}\leq\tg_q\ell_q^{-2\EXP-\eg}\leq\ell_q^{3\ag}\Lg^{-\xfr12}\pa{\zg_q^{\xfr12}\lg_q}^{-1}\lg_{q+1}^{2\EXP}\leq\lg_q^{\bg-1+2b\EXP}<1, \xtag{Gellqtgq}\label{eq:Gellqtgq}\]
for $b$ sufficiently close to 1.
\xbegin{propo}[Gluing step][thm:propo:Gluing]
Let $b,\bg,\ag,\gg$ and $(\dg_q,\lg_q,\Lg,\zg_q,\ell_q,\tg_q)$ be as in Section \ref{Equestimates}, with
\begin{align*}
\ag b<{}&\bg\gg \xtag{GParRel1}\label{eq:GParRel1} \\
b^2(1+\gg)<{}&\xfr{1-\bg}{2\bg}, \xtag{GParRel2}\label{eq:GParRel2} \\
\intertext{Let $[T_1,T_2]\sbs[0,T]$ with $|T_2-T_1|>4\tg_q$. Let $(v_q,p_q,R_q)$ be a strong subsolution on $[0,T]$ which on $[T_1,T_2]$ satisfies the estimates}
\xfr34\dg_{q+2}\leq{}&\rg_q\leq\xfr72\dg_{q+1} \xtag{GUTr}\label{eq:GUTr} \\
\norm{\pint R_q}_0\leq{}&\Lg\varrho_q^{1+\gg} \xtag{GUStrong}\label{eq:GUStrong} \\
\norm{v_q}_{1+\ag}\leq{}&M\dg_q^{\xfr12}\lg_q^{1+\ag} \xtag{GUv}\label{eq:GUv} \\
\norm{v_q}_{\EXP+\eg}\leq{}&M\pa{1+\!\sum_{i=0}^q\dg_i^{\xfr12}\lg_i^{\EXP+\eg}} \xtag{GUFracv}\label{eq:GUFracv} \\
\abs{\pd_t\rg_q}\lsim{}&\rg_q\dg_q^{\xfr12}\lg_q \xtag{GUDer}\label{eq:GUDer}
\end{align*}
with some constant $M>0$, where
\[\rg_q\coloneq\xfr13\opn{tr}R_q\qquad\qquad\qquad\varrho_q\coloneq\xfr{\rg_q}{\Lg}.\]
Define $\rg_{q,i},\varrho_{q,i},\ell_{q,i}$ as in \eqref{eq:rgqiellqi}.
\nipar Then, provided $a\dgreat1$ is sufficiently large, there exists $(\lbar v_q,\lbar p_q,\lbar R_q)$ solution of \eqref{eq:FNSR} on $[0,T]$ such that
\[(\lbar v_q,\lbar p_q,\lbar R_q)=(v_q,p_q,R_q)\qquad\text{on }[0,T]\ssm[T_1,T_2], \xtag{GG-GUEq}\label{eq:GG-GUEq}\]
and on $[T_1,T_2]$ the following estimates hold:
\begin{align*}
\norm{\lbar v_q-v_q}_\ag\lsim{}&\Lg^{\xfr12}\lbar\varrho_q^{\xfr{1+\gg}{2}}\ell_q^\ag \xtag{GGDist}\label{eq:GGDist} \\
\norm{\lbar v_q}_{1+\ag}\lsim{}&\dg_q^{\xfr12}\lg_q^{1+\ag} \xtag{GGv}\label{eq:GGv} \\
\norm{\lbar v_q}_{\EXP+\eg}\leq{}&M\pa{1+\sum_{i=0}^{q+1}\dg_i^{\xfr12}\lg_i^{\EXP+\eg}} \xtag{GGFracv}\label{eq:GGFracv} \\
\norm{\pint{\lbar R}_q}_0\lsim{}&\Lg\lbar\varrho_q^{1+\gg}\ell_q^{-2\ag} \xtag{GGStrong}\label{eq:GGStrong} \\
\xfr78\rg_q\leq{}&\lbar\rg_q\leq\xfr98\rg_q \xtag{GGTr}\label{eq:GGTr} \\
\abs{\pd_t\lbar\rg_q}\lsim{}&\lbar\rg_q\dg_q^{\xfr12}\lg_q \xtag{GGDer}\label{eq:GGDer} \\
\norm{\rg_q-\lbar\rg_q}_0={}&\xfr13\abs{\xints{\T^3}{}\pa[b]{\abs{v_q}^2-\abs{\lbar v_q}^2}\diff x}\lsim\Lg\lbar\varrho_q^{1+\gg}\lg_q^{-\ag}\ell_q^\ag. \xtag{GGTrDist}\label{eq:GGTrDist}
\end{align*}
Moreover, on $[t_{\ubar n},t_{\lbar n}]$ the following additional estimates hold for $t\in I_{i-1}\cup J_i\cup I_i$:
\begin{align*}
\norm{\lbar v_q}_{N+1+\ag}\lsim{}&\dg_q^{\xfr12}\lg_q^{1+\ag}\ell_{q,i}^{-N} \xtag{GGv2}\label{eq:GGv2} \\
\norm{\pint{\lbar R}_q}_{N+\ag}\lsim{}&\Lg\lbar\varrho_q^{1+\gg}\ell_{q,i}^{-N}\ell_q^{-2\ag}\lg_q^{-\ag} \xtag{GGStrong2}\label{eq:GGStrong2} \\
\norm{(\pd_t+\lbar v_q\per\grad)\pint{\lbar R}_q}_{N+\ag}\lsim{}&\Lg\lbar\varrho_q^{1+\gg}\tg_q^{-1}\ell_{q,i}^{-N}\ell_q^{-2\ag}\lg_q^{-\ag}. \xtag{GGMatDer}\label{eq:GGMatDer}
\end{align*}
Regarding the support of the Reynolds stress, we have that
\[\pint{\lbar R}_q(\per,t)\equiv0\qquad\VA t\in\bigcup_{i=\ubar n}^{\lbar n}J_i. \xtag{GGR0}\label{eq:GGR0}\]
In terms of energy, we have that
\[\xints{\T^3}{}\pa[b]{\abs{\lbar v_q}^2(x,t)+\opn{tr}\lbar R_q(x,t)}\diff x=\xints{\T^3}{}\pa[b]{\abs{v_q}^2(x,t)+\opn{tr}R_q(x,t)}\diff x, \xtag{GIC}\label{eq:GIC}\]
and the function
\[\s{T}_g\coloneq\xints{\T^3}{}\2\xints0t\pa[B]{\abs{\xfrl*\lbar v_q}^2-\abs{\xfrl*v_q}^2}\diff s\diff x,\]
satisfies
\[\abs{\pd_t\s{T}_g}\lsim\Lg^{\xfr12}\dg_{q+1}^{\xfr12}\lg_{q+1}^{\EXP+\eg}, \xtag{GGRhoD}\label{eq:GGRhoD}\]
and therefore
\[\abs{\s{T}_g(t)}\lsim t\Lg\lbar\varrho_q^{1+\gg}\zg_q^{\xfr12}\ell_{q,i}^{-N}\ell_q^{-6\ag}\leq T\Lg\lbar\varrho_q^{1+\gg}\zg_q^{\xfr12}\ell_{q,i}^{-N}\ell_q^{-6\ag}.\]
Finally, if $2\ag<\bg\gg$, then
\[\norm{\lbar v_q-v_q}_\ag\lsim\dg_{q+1}^{\xfr12}\ell_q^{(\xfr2b+\xfr12)\ag}. \xtag{GGDistA2S}\label{eq:GGDistA2S}\]
\xend{propo}
The proof closely follows the gluing procedure of \cite[Section 6]{DRSz}, which in turn draws heavily from \cite[Sections 3-4]{BDLSzV}. Recall that the solution is left unchanged outside $[T_1,T_2]$ and the gluing only happens in that interval. More precisely, recalling the decomposition \eqref{eq:Decomp0T}:
\bi
\item The gluing procedure is carried out in the interval
\[J_{\ubar n}\cup\dotso\cup J_{\lbar n}=\pa{t_{\ubar n}-\xfr13\tg_q,t_{\lbar n}+\xfr13\tg_q}; \xtag{GlInt}\label{eq:GlInt}\]
\item The subsolution is left unchanged in $J_{\ubar n-1}\cup J_{\lbar n+1}$;
\item The intervals $I_{\ubar n-1}$ and $I_{\lbar n}$ are used as cutoff regions between the glued and unglued subsolutions.
\ei
Recall also that, since the trace $\rg_q={}^1/_3\opn{tr}R_q$ has different lower and upper bounds on $[T_1,T_2]$ (respectively of order $\dg_{q+2}$ and $\dg_{q+1}$), mollification with different parameters $\ell_{q,i}$ depending on $\rg_q(t_i)$ on intervals of size $\tg_q$ around the points $t_i$ is necessary.
\nipar We will also make use of the following estimates.
\xbegin{lemma}[Material derivative estimates for subsolutions and potentials][thm:lemma:MatDerEstSubsPot]
Let $(v,p,R)$, $(v',p',R')$ be two solutions of \eqref{eq:FNSR}, and let $z\coloneq(-\Dg)^{-1}\opn{curl}v,z'\coloneq(-\Dg)^{-1}\opn{curl}v$. Then the following estimates hold for every $N\in\N,\ag\in(0,1)$:
\begin{align*}
\norm{(\pd_t+v\per\grad+\xfrl)(v-v')}_{N+\ag}\lsim{}&\norm{v-v'}_{N+\ag}(\norm{v}_\ag+\norm{v'}_\ag) \\
&{}+\norm{v-v'}_\ag(\norm{v}_{N+\ag}+\norm{v'}_{N+\ag})+\norm{R-R'}_{N+1+\ag} \xtag{MatDerEst1}\label{eq:MatDerEst1} \\
\norm{(\pd_t+v\per\grad+\xfrl)(z'-z)}_{N+\ag}\lsim{}&\norm{z'-z}_{N+\ag}\norm{v}_{1+\ag}+\norm{z'-z}_\ag\norm{v}_{N+1+\ag} \\
&{}+\norm{R'-R}_{N+\ag}. \xtag{MatDerEst2}\label{eq:MatDerEst2}
\end{align*}
\xend{lemma}
The proof of \eqref{eq:MatDerEst1} can be found in the proof of \cite[Proposition 5.3]{DR}, whereas that of \eqref{eq:MatDerEst2} can be found within the proof of \cite[Proposition 5.4]{DR}. Note that, since the proofs require the Schauder estimates of \kcref{thm:lemma:Schauder}, this result does not hold for $\ag=0$.
\nipar\begin{qeddim*}[\kcref{thm:propo:Gluing}]
\begin{center}
\tbold{Step 1: Mollification}
\end{center}
Let $\fg$ be a standard mollification kernel in space and define
\begin{align*}
v_{\ell_{q,i}}\coloneq{}&v_q\ast\fg_{\ell_{q,i}} \\
p_{\ell_{q,i}}\coloneq{}&p_q\ast\fg_{\ell_{q,i}}+\xfr13(|v_q|^2\ast\fg_\ell-|v_{\ell_{q,i}}|^2) \\
\pint R_{\ell_{q,i}}\coloneq{}&\pint R_q\ast\fg_{\ell_{q,i}}+(v_q\potimes v_q)\ast\fg_{\ell_{q,i}}-v_{\ell_{q,i}}\potimes v_{\ell_{q,i}}.
\end{align*}
With this definition, \eqref{eq:FNSR} holds for the triple $(v_{\ell_{q,i}},p_{\ell_{q,i}}. R_{\ell_{q,i}})$. Using the estimates \eqref{eq:GUStrong} and \eqref{eq:GUv}, together with \eqref{eq:Mollif}, we deduce that
\begin{align*}
\norm{v_{\ell_{q,i}}-v_q}_\ag\lsim{}&\dg_q^{\xfr12}\lg_q^{1+\ag}\ell_{q,i}=\Lg^{\xfr12}\varrho_{q,i}^{\xfr{1+\gg}{2}}\ell_q^\ag \xtag{GMDist}\label{eq:GMDist} \\
\norm{v_{\ell_{q,i}}}_{N+1+\ag}\lsim{}&\dg_q^{\xfr12}\lg_q^{1+\ag}\ell_{q,i}^{-N} \xtag{GMv}\label{eq:GMv} \\
\norm{v_{\ell_{q,i}}}_{\EXP+\eg}\lsim{}&\Lg^{\xfr12} \xtag{GMFracv}\label{eq:GMFracv} \\
\norm{\pint R_{\ell_{q,i}}}_{N+\ag}\lsim{}&\Lg\varrho_q^{1+\gg}\ell_{q,i}^{-N-\ag}+\dg_q\lg_q^{2+2\ag}\ell_{q,i}^{2-N-\ag} \xtag{GMR}\label{eq:GMR} \\
\lsim{}&\Lg\varrho_q^{1+\gg}\ell_{q,i}^{-N-\ag}+\Lg\varrho_{q,i}^{1+\gg}\ell_{q,i}^{-N}\ell_q^\ag \\
\abs{\xints{\T^3}{}\2\abs{v_q}^2-\abs{v_{\ell_{q,i}}}^2\diff x}\lsim{}&\dg_q\lg_q^{2+2\ag}\ell_{q,i}^2=\Lg\varrho_{q,i}^{1+\gg}\ell_q^{2\ag}. \xtag{GMEk}\label{eq:GMEk}
\end{align*}
To obtain \eqref{eq:GMEk}, we also use the trivial identity $\int f\ast\fg_\ell=\int f$ for $f=|v_q|^2$.
\begin{center}
\tbold{Step 2: Gluing procedure}
\end{center}
Let $\{I_i\}_{\ubar n\leq i\leq\lbar n}$ be the sequence of intervals corresponding to $[T_1,T_2]$ according to \kcref{thm:defi:Subintervs+} above. We now fix a partition of unity on $[0,T]$
\[\sum_{\mcl{i=\ubar n-1}}^{\lbar n+1}\xg_i\equiv1\]
subordinate to the decomposition \eqref{eq:Decomp0T}, i.e. $[0,T]=J_{\ubar n-1}\cup I_{\ubar n-1}\cup[J_{\ubar n}\cup\dotso\cup J_{\lbar n}]\cup I_{\lbar n}\cup J_{\lbar n+1}$. More precisely, for each $\ubar n-1\leq i\leq\lbar n+1$, the function $\xg_i\geq0$ satisfies
\[\opn{supp}\xg_i\sbs I_{i-1}\cup J_i\cup I_i, \qquad\qquad \xg_i\eval{J_i}\equiv1 \qquad\qquad \abs{\pd_t^N\xg_i}\lsim\tg_q^{-N}\quad\VA N\geq0.\]
We define
\[\lbar v_q\coloneq\2\sum_{\mcl{i=\ubar n-1}}^{\lbar n+1}\xg_iv_i,\qquad\qquad\lbar p_q\stp1\coloneq\2\sum_{\mcl{i=\ubar n-1}}^{\lbar n+1}\xg_ip_i, \e[10] \xtag{GGDefvp1}\label{eq:GGDefvp1}\]
where $(v_i,p_i)$ is defined as follows. For $\ubar n\leq i\leq\lbar n$ we define $(v_i,p_i)$ as the solution of
\[\begin{sistema}
\pd_tv_i+\opn{div}(v_i\otimes v_i)+\grad p_i+\xfrl v_i=0 \\
\opn{div}v_i=0 \\
v_i(\per,t_i)=v_{\ell_{q,i}}(\per,t_i)
\end{sistema}, \xtag{GIVP}\label{eq:GIVP}\]
and set $(v_i,p_i)=(v_q,p_q)$ for $i=\ubar n-1$ and $i=\lbar n+1$. Thus, we note first of all that $\opn{div}\lbar v_q=0$, and moreover
\[(\lbar v_q,\lbar p_q\stp1)=(v_q,p_q), \qquad t\in[0,T]\ssm[T_1,T_2].\]
Next, we define $\lbar R_q$. We have that $\xg_i+\xg_{i+1}=1$ for $t\in J_i\cup I_i\cup J_{i+1}$, and therefore
\begin{align*}
\pd_t\lbar v_q+\opn{div}(\lbar v_q\otimes\lbar v_q)+\grad\lbar p_q\stp1+\xfrl\lbar v_q={}&\pd_t\xg_i\per(v_i-v_{i+1}) \\
&{}-\xg_i\per(1-\xg_i)\opn{div}(\!(v_i-v_{i+1})\otimes(v_i-v_{i+1})\!) \\
&{}-\opn{div}(\xg_iR_i+(1-\xg_i)R_{i+1})
\end{align*}
for all $\ubar n-1\leq i\leq\lbar n$, where we wrote $R_i=0$ for $\ubar n\leq i\leq\lbar n$ and $R_i=R_q$ otherwise. Thus, recalling the operator $\s{R}$ from \kcref{thm:defi:AntiDivR}, set
\[\begin{aligned}
\pint{\lbar R}_q\stp1\coloneq{}&\case{
\mat{c} \e[80]-\pd_t\xg_i\s{R}(v_i-v_{i+1}) \\ {\e[-8]}+\xg_i\per(1-\xg_i)(v_i-v_{i+1})\mathbin{\pint\otimes}(v_i-v_{i+1}) \emat & t\in I_i \\
0 & t\in J_i\cup J_{i+1}
} \\
\pint{\lbar R}_q\stp2\coloneq{}&\sum_{\mcl{i=\ubar n-1}}^{\lbar n+1}\xg_i\pint R_i=(\xg_{\ubar n-1}+\xg_{\lbar n+1})\pint R_q,
\end{aligned}
\xtag{GGDefR12}\label{eq:GGDefR12}\]
and
\[\lbar p_q\stp2\coloneq\sum_{\mcl{i=\ubar n-1}}^{\lbar n+1}\xg_i\per(1-\xg_i)\pa{\abs{v_i-v_{i+1}}^2-\fint_{\T^3}\abs{v_i-v_{i+1}}^2\diff x}. \xtag{GGDefp2}\label{eq:GGDefp2}\]
Finally, we define
\[\lbar R_q=\pint{\lbar R}_q\stp1+\pint{\lbar R}_q\stp2+\lbar\rg_q\opn{Id},\qquad\lbar p_q\coloneq\lbar p_q\stp1+\lbar p_q\stp2, \xtag{GGDefRp}\label{eq:GGDefRp}\]
where
\[\lbar\rg_q\coloneq\rg_q+\xfr13\fint_{\T^3}\pa[b]{\abs{v_q}^2-\abs{\lbar v_q}^2}\diff x. \xtag{GGDefrho}\label{eq:GGDefrho}\]
Define also
\[\s{T}_g(t)\coloneq\xfr13\fint_{\T^3}\pa[B]{\abs{\xfrl*v_q}^2-\abs{\xfrl*\lbar v_q}^2}\diff x.\]
By construction, we have that
\[\pd_t\lbar v_q+\opn{div}(\lbar v_q\otimes\lbar v_q)+\grad\lbar p_q=-\opn{div}\lbar R_q,\]
and \eqref{eq:GG-GUEq} and \eqref{eq:GIC} hold. Moreover
\[\pint{\lbar R}_q=0\qquad\VA t\in\bigcup_{i=\ubar n}^{\lbar n}J_i.\]
\begin{center}
\tbold{Step 3: Stability estimates on classical solutions}
\end{center}
Throughout this step and the next, we will assume estimate \eqref{eq:GGTr}, which will be proved in Step 5 below. This estimate will allow us to replace $\lbar\rg_q$ with $\rg_q$ and viceversa whenever we need to do so in our estimates, since the two are of the same order.
\nipar Let us consider for the moment $\ubar n\leq i\leq\lbar n$. We recall the classical existence result for solutions of \eqref{eq:GIVP} found in \cite[Proposition 3.5]{DR}, by which $(v_i,p_i)$ in \eqref{eq:GGDefvp1} above is defined at least on an interval of length $\sim\|v_{\ell_{q,i}}\|_{1+\ag}^{-1}$. By \eqref{eq:GMv} and \eqref{eq:DefTau}, we have that
\[\norm{v_{\ell_{q,i}}}_{1+\ag}\lsim\dg_q^{\xfr12}\lg_q^{1+\ag}=\tg_q^{-1}\ell_q^{4\ag}\lg_q^\ag\leq\tg_q^{-1}.\]
Therefore, provided $a\dgreat1$ is sufficiently large, $v_i$ is defined on $I_{i-1}\cup J_i\cup I_i$, so that $\lbar v_q$ in \eqref{eq:GGDefvp1} is well-defined.
\nipar Next, we deduce from \eqref{eq:GUDer} that $|\pd_t\log\rg_q|\leq\dg_q^{{}^1\!/_2}\lg_q=\tg_q^{-1}\ell_q^{4\ag}$, so that, by assuming $a\dgreat1$ is sufficiently large, we may ensure that
\[\rg_q(t_1)\leq4\rg_q(t_2)\qquad\VA t_1,t_2\in I_{i-1}\cup J_i\cup I_i, \xtag{GRhoRatio}\label{eq:GRhoRatio}\]
for any $i$. In particular $\rg_q\sim\rg_{q,i}$ and $\varrho_q\sim\varrho_{q,i}$ in $I_{i-1}\cup J_i\cup I_i$. We apply \kcref{thm:lemma:MatDerEstSubsPot} to $(v_{\ell_{q,i}},p_{\ell_{q,i}},R_{\ell_{q,i}})$ and $(v_i,p_i,0)$, which, using \eqref{eq:GMDist}, \eqref{eq:GMv}, \eqref{eq:GMR}, and \eqref{eq:GRhoRatio}, immediately yields
\begin{align*}
\norm{(\pd_t+v_{\ell_{q,i}}\per\grad+\xfrl)(v_i-v_{\ell_{q,i}})}_{N+\ag}\leq{}&\norm{v_i-v_{\ell_{q,i}}}_{N+\ag}(\norm{v_i}_{1+\ag}+\norm{v_{\ell_{q,i}}}_{1+\ag}) \\
&{}+\norm{v_i-v_{\ell_{q,i}}}_\ag(\norm{v_i}_{N+1+\ag}+\norm{v_{\ell_{q,i}}}_{N+1+\ag}) \\
&{}+\norm{\pint R_{\ell_{q,i}}}_{N+1+\ag} \\
{}\lsim{}&\dg_q^{\xfr12}\lg_q^{1+\ag}\ell_{q,i}^{-N}\norm{v_i-v_{\ell_{q,i}}}_\ag+\dg_q^{\xfr12}\lg_q^{1+\ag}\norm{v_i-v_{\ell_{q,i}}}_{N+\ag} \\
&{}+\Lg\varrho_{q,i}^{1+\gg}\ell_{q,i}^{-1-N-\ag}. \xtag{GGMDvi-vl.Lemma}\label{eq:GGMDvi-vl.Lemma}
\end{align*}
Combining this with \kcref{thm:propo:EstTraspDiff}, working first for the case $N=0$, then for $N=1$, and finally for the general case, yields
\[\norm{v_i-v_{\ell_{q,i}}}_{N+\ag}\lsim\xints{t_i}{t}\tg_q^{-1}\norm{v_i-v_{\ell_{q,i}}}_{N+\ag}\diff s+\Lg\tg_q\varrho_{q,i}^{1+\gg}\ell_{q,i}^{-N-1-\ag}.\]
The Grönwall inequality then implies that
\[\norm{v_i-v_{\ell_{q,i}}}_{N+\ag}\lsim\Lg\tg_q\varrho_{q,i}^{1+\gg}\ell_{q,i}^{-N-1-\ag}\lsim\Lg^{\xfr12}\varrho_{q,i}^{\xfr{1+\gg}{2}}\ell_{q,i}^{-N}\ell_q^\ag. \xtag{GGDistN}\label{eq:GGDistN}\]
The case $N=0$ of \eqref{eq:GGDistN}, together with \eqref{eq:GMDist}, leads to
\[\norm{\lbar v_q-v_q}_\ag\leq\sum_{\mcl{j=i-1}}^{i+1}(\norm{v_j-v_{\ell_{q,j}}}_\ag+\norm{v_{\ell_{q,j}}-v_q}_\ag)\lsim\Lg^{\xfr12}\varrho_{q,i}^{\xfr{1+\gg}{2}}\ell_q^\ag. \xtag{GGDist2}\label{eq:GGDist2}\]
By \eqref{eq:GGTr}, this is equivalent to \eqref{eq:GGDist}.
\nipar The case $N=1$ of \eqref{eq:GGDistN} leads to
\[\norm{v_i-v_{\ell_{q,i}}}_{1+\ag}\lsim\Lg^{\xfr12}\varrho_{q,i}^{\xfr{1+\gg}{2}}\ell_{q,i}^{-1}\ell_q^\ag=\Lg^{\xfr12}\zg_q^{\xfr12}\lg_q^{1+\ag}\ell_q^{-\ag}\ell_q^\ag=\dg_q^{\xfr12}\lg_q^{1+\ag}.\]
Combining the above estimate with \eqref{eq:GUv} outside the gluing region and in $I_{\ubar n-1}\cup I_{\lbar n}$, and with \eqref{eq:GMv} in $J_{\ubar n}\cup I_{\ubar n}\cup\dotso\cup I_{\lbar n-1}\cup J_{\lbar n}$, we deduce that \eqref{eq:GGv} is verified. More generally, as we did above for $N=0$, we deduce from \eqref{eq:GMv} and \eqref{eq:GGDistN} that
\[\norm{\lbar v_q}_{1+N+\ag}\leq\sum_{j=i-1}^{i+1}\xg_j\pa{\norm{v_j-v_{\ell_{q,j}}}_{1+N+\ag}+\norm{v_{\ell_{q,j}}}_{1+N+\ag}}\lsim\dg_q^{\xfr12}\lg_q^{1+\ag}\ell_{q,i}^{-N}\qquad\qquad\VA t\in\opn{supp}\xg_i.\]
We have used the fact that $\ell_{q,i}\sim\ell_{q,i+1}\sim\ell_{q,i-1}$. The above inequality coincides with \eqref{eq:GGv2}.
\nipar We also remark the following simple interpolation of the $N=0$ and $N=1$ cases of \eqref{eq:GGDistN}, which will be used in Step 5 below.
\[\norm{v_i-v_{\ell_{q,i}}}_{\EXP+\eg}\lsim\norm{v_i-v_{\ell_{q,i}}}_\ag^{1-\EXP-\eg}\norm{v_i-v_{\ell_{q,i}}}_{1+\ag}^{\EXP+\eg}\lsim\Lg\tg_q\varrho_{q,i}^{1+\gg}\ell_{q,i}^{-\EXP-\eg-1-\ag}. \xtag{GGFracDist}\label{eq:GGFracDist}\]
\nipar Further in the proof, we will need estimates for $\|v_i-v_{i+1}\|_{N+\ag}$ and $\|(\pd_t+v_{\ell_{q,i}}\per\grad)(v_i-v_{i+1})\|_{N+\ag}$. Concerning the former, by applying the triangle inequality, we see that
\[\norm{v_i-v_{i+1}}_{N+\ag}\leq\norm{v_i-v_{\ell_{q,i}}}_{N+\ag}+\norm{v_{\ell_{q,i}}-v_{\ell_{q,i+1}}}_{N+\ag}+\norm{v_{\ell_{q,i+1}}-v_{i+1}}_{N+\ag}.\]
The first and third term are estimated by \eqref{eq:GGDistN}. The second one is readily shown to satisfy the same estimate, so that
\[\norm{v_i-v_{i+1}}_{N+\ag}\lsim\Lg^{\xfr12}\varrho_{q,i}^{\xfr{1+\gg}{2}}\ell_{q,i}^{-N}\ell_q^\ag. \xtag{GGvi-vi+1}\label{eq:GGvi-vi+1}\]
As for the material derivative, we note that
\begin{align*}
\norm{(\pd_t+v_{\ell_{q,i}}\per\grad)(v_i-v_{i+1})}_{N+\ag}\leq{}&\norm{(\pd_t+v_{\ell_{q,i}}\per\grad)(v_i-v_{\ell_{q,i}})}_{N+\ag} \\
&{}+\norm{(\pd_t+v_{\ell_{q,i}}\per\grad)(v_{\ell_{q,i}}-v_{\ell_{q,i+1}})}_{N+\ag} \\
&{}+\norm{(v_{\ell_{q,i}}-v_{\ell_{q,i+1}})\per\grad(v_{\ell_{q,i+1}}-v_{i+1})}_{N+\ag} \\
&{}+\norm{(\pd_t+v_{\ell_{q,i+1}}\per\grad)(v_{\ell_{q,i+1}}-v_{i+1})}_{N+\ag} \\
{}\eqcolon{}&I+II+III+IV.
\end{align*}
We start by estimating $I$. Inserting \eqref{eq:GGDistN} into \eqref{eq:GGMDvi-vl.Lemma} gives us
\[\norm{(\pd_t+v_{\ell_{q,i}}\per\grad+\xfrl)(v_i-v_{\ell_{q,i}})}_{N+\ag}\lsim\Lg\varrho_{q,i}^{1+\gg}\ell_{q,i}^{-1-N-\ag}. \xtag{GGMD+vi-vl}\label{eq:GGMD+vi-vl}\]
Using \eqref{eq:Gellqtgq}, we can easily conclude, by interpolating in \eqref{eq:GGDistN}, that
\begin{align*}
\norm{\xfrl(v_i-v_{\ell_{q,i}})}_{N+\ag}\lsim{}&\norm{v_i-v_{\ell_{q,i}}}_{N+\ag+2\EXP+\eg}\lsim\norm{v_i-v_{\ell_{q,i}}}_{N+\ag}^{1-2\EXP-\eg}\norm{v_i-v_{\ell_{q,i}}}_{N+\ag+1}^{2\EXP+\eg} \\
{}\lsim{}&\Lg\tg_q\ell_{q,i}^{-2\EXP-\eg}\varrho_{q,i}^{1+\gg}\ell_{q,i}^{-N-1-\ag}\lsim\Lg\varrho_{q,i}^{1+\gg}\ell_{q,i}^{-1-N-\ag},
\end{align*}
so that
\[I\leq\norm{(\pd_t+v_{\ell_{q,i}}\per\grad+\xfrl)(v_i-v_{\ell_{q,i}})}_{N+\ag}+\norm{\xfrl(v_i-v_{\ell_{q,i}})}_{N+\ag}\lsim\Lg\varrho_{q,i}^{1+\gg}\ell_{q,i}^{-1-N-\ag}.\]
The term $IV$ is handled similarly. Since $v_{\ell_{q,i}}-v_{\ell_{q,i+1}}$ obeys the bound \eqref{eq:GGvi-vi+1}, using \kcref{thm:lemma:MatDerEstSubsPot} we conclude that $II$ also obeys the above bound. We now consider $III$, which can be estimated as follows:
\[III=\norm{(v_{\ell_{q,i}}-v_{\ell_{q,i+1}})\per\grad(v_{i+1}-v_{\ell_{q,i+1}})}_{N+\ag}\lsim\Lg^{\xfr12}\varrho_{q,i}^{\xfr{1+\gg}{2}}\ell_{q,i}^\ag\per\Lg^{\xfr12}\varrho_{q,i}^{\xfr{1+\gg}{2}}\ell_{q,i}^\ag\per\ell_{q,i}^{-N-1}=\Lg\varrho_{q,i}^{1+\gg}\ell_{q,i}^{2\ag-N-1},\]
in particular satisfying the bound \eqref{eq:GGMD+vi-vl}. We thus conclude that
\[\norm{(\pd_t+v_{\ell_{q,i}}\per\grad)(v_i-v_{i+1})}_{N+\ag}\lsim\Lg\varrho_{q,i}^{1+\gg}\ell_{q,i}^{-1-N-\ag}=\Lg^{\xfr12}\varrho_{q,i}^{\xfr{1+\gg}{2}}\tg_q^{-1}\ell_q^{3\ag}\lg_q^\ag\ell_{q,i}^{-N-\ag}. \xtag{GGMDvi-vi+1}\label{eq:GGMDvi-vi+1}\]
\begin{center}
\tbold{Step 4: Estimates on the new Reynolds stress}
\end{center}
As is done in \cite[Section 3.3]{BDLSzV}, we define the vector potentials
\[z_i\coloneq(-\Dg)^{-1}\opn{curl}v_i\qquad z_{\ell_{q,i}}\coloneq(-\Dg)^{-1}\opn{curl}v_{\ell_{q,i}},\qquad z_q\coloneq(-\Dg)^{-1}\opn{curl}v_q\]
and, by a combination of \kcref{thm:lemma:MatDerEstSubsPot}, \kcref{thm:propo:EstTraspDiff}, and the Grönwall inequality similar to the one used to obtain \eqref{eq:GGDistN} in Step 3 above, obtain that
\[\norm{z_i-z_{\ell_{q,i}}}_{N+\ag}\lsim\Lg\tg_q\varrho_{q,i}^{1+\gg}\ell_{q,i}^{-N}\ell_q^{-\ag}. \xtag{GGzi-zl}\label{eq:GGzi-zl}\]
Combining this with \kcref{thm:lemma:MatDerEstSubsPot} yields
\[\norm{(\pd_t+v_{\ell_{q,i}}\per\grad+\xfrl)(z_i-z_{\ell_{q,i}})}_{N+\ag}\lsim\Lg\varrho_{q,i}^{1+\gg}\ell_{q,i}^{-N}\ell_q^{-\ag}. \xtag{GGStep4MatDer+}\label{eq:GGStep4MatDer+}\]
By \eqref{eq:Gellqtgq} and the Schauder estimates of \kcref{thm:lemma:Schauder}, we deduce that
\begin{align*}
\norm{\xfrl(z_i-z_{\ell_{q,i}})}_{N+\ag}\lsim{}&\norm{z_i-z_{\ell_{q,i}}}_{N+\ag+2\EXP+\eg}\lsim\norm{v_i-v_{\ell_{q,i}}}_{N+\ag-1+2\EXP+\eg} \\
{}\lsim{}&\Lg\tg_q\varrho_{q,i}^{1+\gg}\ell_{q,i}^{-N-2\EXP-\eg}\ell_q^{-\ag}\lsim\Lg\varrho_{q,i}^{1+\gg}\ell_{q,i}^{-N}\ell_q^{-\ag}. \xtag{GStep4Lapl}\label{eq:GStep4Lapla}
\end{align*}
By the triangle inequality, \eqref{eq:GStep4Lapla}, and \eqref{eq:GGStep4MatDer+}, we thus conclude that
\[\norm{(\pd_t+v_{\ell_{q,i}}\per\grad)(z_i-z_{\ell_{q,i}})}_{N+\ag}\lsim\Lg\varrho_{q,i}^{1+\gg}\ell_{q,i}^{-N}\ell_q^{-\ag}. \xtag{GGMDzi-zl}\label{eq:GGMDzi-zl}\]
Both \eqref{eq:GGzi-zl} and \eqref{eq:GGMDzi-zl} are valid in $I_{i-1}\cup J_i\cup I_i$ for any $\ubar n\leq i\leq\lbar n$.
\nipar The sequel of this proof will require estimates on $z_i-z_{i+1}$, which means we must now bound $z_{\ell_{q,i}}-z_{\ell_{q,i+1}}$. We note that $z_{\ell_{q,i}}=z_q\ast\fg_{\ell_{q,i}}$ and $z_{\ell_{q,i+1}}=z_q\ast\fg_{\ell_{q,i+1}}$, so that, using \eqref{eq:Mollif} and Schauder estimates (\kcref{thm:lemma:Schauder}), we get
\[\norm{z_{\ell_{q,i}}-z_{\ell_{q,i+1}}}_{N+\ag}\lsim\norm{z_q}_{2+\ag}(\ell_{q,i}^{2-N}+\ell_{q,i+1}^{2-N})\lsim\Lg\tg_q\varrho_{q,i}^{1+\gg}\ell_q^{-2\ag}\lg_q^{-\ag}\ell_{q,i}^{-N}. \xtag{GGzlqi-zlqi+1}\label{eq:GGzlqi-zlqi+1}\]
The final estimate for $z_i-z_{i+1}$ is thus
\[\norm{z_i-z_{i+1}}_{N+\ag}\leq\norm{z_i-z_{\ell_{q,i}}}_{N+\ag}+\norm{z_{\ell_{q,i}}-z_{\ell_{q,i+1}}}_{N+\ag}+\norm{z_{\ell_{q,i+1}}-z_{i+1}}_{N+\ag}\lsim\Lg\varrho_{q,i}^{1+\gg}\tg_q\ell_{q,i}^{-N}\ell_q^{-2\ag}\lg_q^{-\ag}, \xtag{GGzi-zi+1}\label{eq:GGzi-zi+1}\]
slightly coarser than \eqref{eq:GGzi-zl} above. As for the material derivatives, we must estimate
\begin{align*}
\norm{(\pd_t+v_{\ell_{q,i}}\per\grad)(z_i-z_{i+1})}_{N+\ag}\leq{}&\norm{(\pd_t+v_{\ell_{q,i}}\per\grad)(z_i-z_{\ell_{q,i}})}_{N+\ag} \\
&{}+\norm{(\pd_t+v_{\ell_{q,i}}\per\grad)(z_{\ell_{q,i}}-z_{\ell_{q,i+1}})}_{N+\ag} \\
&{}+\norm{(v_{\ell_{q,i}}-v_{\ell_{q,i+1}})\per\grad(z_{\ell_{q,i+1}}-z_{i+1})}_{N+\ag} \\
&{}+\norm{(\pd_t+v_{\ell_{q,i+1}}\per\grad)(z_{\ell_{q,i+1}}-z_{i+1})}_{N+\ag} \\
{}\eqcolon{}&I+II+III+IV.
\end{align*}
The terms I and IV are estimated by \eqref{eq:GGMDzi-zl}. To estimate $II$, we apply \kcref{thm:lemma:MatDerEstSubsPot} to $(v_{\ell_{q,i}},p_{\ell_{q,i}},R_{\ell_{q,i}})$ and $(v_{\ell_{q,i+1}},p_{\ell_{q,i+1}},R_{\ell_{q,i+1}})$ and use \eqref{eq:GMv} and \eqref{eq:GGzlqi-zlqi+1}, obtaining that
\[\norm{(\pd_t+v_{\ell_{q,i}}\per\grad+\xfrl)(z_{\ell_{q,i}}-z_{\ell_{q,i+1}})}_{N+\ag}\lsim\rg_{q,i}^{1+\gg}\ell_{q,i}^{-N-\ag}(\ell_{q,i}^\ag\ell_q^{2\ag}+1).\]
We then note that, by interpolation
\[\norm{\xfrl(z_{\ell_{q,i}}-z_{\ell_{q,i+1}})}_{N+\ag}\lsim\Lg\tg_q\varrho_{q,i}^{1+\gg}\ell_{q,i}^{-N-2\EXP-\eg}\ell_q^{-2\ag}\lg_q^{-\ag}\lsim\Lg\varrho_{q,i}^{1+\gg}\ell_{q,i}^{-N}\ell_q^{-2\ag}\lg_q^{-\ag}.\]
The above two bounds combine to yield
\[II\lsim\Lg\varrho_{q,i}^{1+\gg}\ell_{q,i}^{-N}\ell_q^{-2\ag}\lg_q^{-\ag}. \xtag{GGMDzconvol}\label{eq:GGMDzconvol}\]
Coming to $III$, we estimate it by combining \eqref{eq:GGzi-zl} with \eqref{eq:GGvi-vi+1}:
\begin{align*}
\norm{(v_{\ell_{q,i}}-v_{\ell_{q,i+1}})\per\grad(z_{i+1}-z_{\ell_{q,i+1}})}_{N+\ag}\lsim{}&\Lg^{\xfr12}\varrho_{q,i}^{\xfr{1+\gg}{2}}\ell_{q,i}^\ag\per\Lg\tg_q\varrho_{q,i}^{1+\gg}\ell_q^{-\ag}\per\ell_{q,i}^{-N-1} \\
{}\leq{}&\Lg\varrho_{q,i}^{1+\gg}\ell_{q,i}^{\ag-N}. \xtag{GGvChange}\label{eq:GGvChange}
\end{align*}
Combining \eqref{eq:GGMDzi-zl}, \eqref{eq:GGMDzconvol}, and \eqref{eq:GGvChange}, we thus obtain that
\[\norm{(\pd_t+v_{\ell_{q,i}}\per\grad)(z_i-z_{I+1})}_{N+\ag}\lsim\Lg\varrho_{q,i}^{1+\gg}\ell_{q,i}^{-N}\ell_q^{-2\ag}\lg_q^{-\ag}. \xtag{GGMDzi-zi}\label{eq:GGMDzi-zi}\]
Recalling the expression for $\pint{\lbar R}_q$ in \eqref{eq:GGDefRp}, using \eqref{eq:GRhoRatio}, \eqref{eq:GGzi-zi+1}, and \eqref{eq:GGvi-vi+1}, we obtain, as in the proof of \cite[Proposition 4.3]{BDLSzV}, that
\begin{align*}
\norm{\pint{\lbar R}_q}_{N+\ag}\lsim{}&\tg_q^{-1}\norm{z_i-z_{i+1}}_{N+\ag}+\norm{v_i-v_{i+1}}_{N+\ag}\norm{v_i-v_{i+1}}_\ag \\
{}\lsim{}&\Lg\varrho_{q,i}^{1+\gg}\ell_{q,i}^{-N}\ell_q^{-2\ag}\lg_q^{-\ag}+\Lg^{\xfr12}\varrho_{q,i}^{\xfr{1+\gg}{2}}\ell_{q,i}^{-N+\ag}\per\Lg^{\xfr12}\varrho_{q,i}^{\xfr{1+\gg}{2}}\ell_{q,i}^\ag \\
{}\lsim{}&\Lg\varrho_{q,i}^{1+\gg}\ell_{q,i}^{-N}\ell_q^{-\ag}(\ell_q^{-\ag}\lg_q^{-\ag}+\ell_{q,i}^{3\ag}).
\end{align*}
This, together with \eqref{eq:GGTr}, gives us \eqref{eq:GGStrong2}.
\nipar As for \eqref{eq:GGMatDer}, we note that
\begin{align*}
\norm{(\pd_t+\lbar v_q\per\grad)\pint{\lbar R}_q}_{N+\ag}\lsim{}&\tg_q^{-2}\norm{z_i-z_{i+1}}_{N+\ag}+\tg_q^{-1}\norm{(\pd_t+v_{\ell_{q,i}}\per\grad)(z_i-z_{i+1})}_{N+\ag}+{} \\
&{}+\tg_q^{-1}\norm{v_{\ell_{q,i}}}_{1+\ag}\norm{z_i-z_{i+1}}_{N+\ag}+\tg_q^{-1}\norm{v_{\ell_{q,i}}}_{N+1+\ag}\norm{z_i-z_{i+1}}_\ag \\
&{}+\tg_q^{-1}\norm{v_i-v_{i+1}}_{N+\ag}\norm{v_i-v_{i+1}}_\ag \\
&{}+\norm{(\pd_t+v_{\ell_{q,i}}\per\grad)(v_i-v_{i+1})}_{N+\ag}\norm{v_i-v_{i+1}}_\ag \\
&{}+\norm{(\pd_t+v_{\ell_{q,i}}\per\grad)(v_i-v_{i+1})}_\ag\norm{v_i-v_{i+1}}_{N+\ag} \\
&{}+\norm{(v_{\ell_{q,i}}-\lbar v_q)\per\grad\pint{\lbar R}_q}_{N+\ag}.
\end{align*}
Combining the above bound on $\pint{\lbar R}_q$ with \eqref{eq:GGzi-zi+1}, \eqref{eq:GGMDzi-zi}, \eqref{eq:GMv}, \eqref{eq:GGvi-vi+1}, \eqref{eq:GGMDvi-vi+1}, and the bound \eqref{eq:GGDist2} applied to $v_{\ell_{q,i}}-\lbar v_q$, we obtain that
\[\norm{(\pd_t+\lbar v_q\per\grad)\pint{\lbar R}_q}_{N+\ag}\lsim\Lg\tg_q^{-1}\varrho_{q,i}^{1+\gg}\ell_{q,i}^{-N}\ell_q^{-2\ag}\lg_q^{-\ag},\xtag{GGPreMatDer}\label{eq:GGPreMatDer}\]
which yields \eqref{eq:GGMatDer} once combined with \eqref{eq:GGTr}. \\\begin{center}
\tbold{Step 5: $\lbar\rg_q$, $\s{T}_g$, and \eqref{eq:GGDistA2S}}
\end{center}
Next, we estimate $\lbar\rg_q$, recalling its definition in \eqref{eq:GGDefrho}. We wish to estimate $\lbar\rg_q-\rg_q$. We note that
\[\norm{\lbar\rg_q-\rg_q}_0=\xfr13\abs{\xints{\T^3}{}\2\abs{\lbar v_q}^2-\abs{v_q}^2\diff x}\leq\abs{\xints{\T^3}{}\2\abs{\lbar v_q}^2-\abs{v_{\ell_{q,i}}}^2\diff x}+\abs{\xints{\T^3}{}\2\abs{v_{\ell_{q,i}}}^2-\abs{v_q}^2\diff x}.\]
The second term above is already estimated by \eqref{eq:GMEk}, so we proceed to estimate the first term.
\nipar As in \cite[Proposition 4.4]{BDLSzV}, one has that
\begin{align*}
\abs{\lbar v_q}^2-\abs{v_{\ell_{q,i}}}^2={}&\xg_i(|v_i|^2-|v_{\ell_{q,i}}|^2)+(1-\xg_i)(|v_{i+1}|^2-|v_{\ell_{q,i+1}}|^2) \\
&{}+(1-\xg_i)(|v_{\ell_{q,i+1}}|^2-|v_{\ell_{q,i}}|^2)-\xg_i(1-\xg_i)|v_i-v_{i+1}|^2.
\end{align*}
Therefore
\begin{align*}
\abs{\xints{\T^3}{}\2\abs{\lbar v_q}^2-\abs{v_{\ell_{q,i}}}^2\diff x}\leq{}&\abs{\xints{\T^3}{}\2\abs{v_i}^2-\abs{v_{\ell_{q,i}}}^2\diff x}+\abs{\xints{\T^3}{}\2\abs{v_{i+1}}^2-\abs{v_{\ell_{q,i+1}}}^2\diff x} \\
&{}+\abs{\xints{\T^3}{}\2\abs{v_{\ell_{q,i}}}^2-\abs{v_{\ell_{q,i+1}}}^2\diff x}+\abs{\xints{\T^3}{}\2\abs{v_i-v_{i+1}}^2\diff x}. \xtag{GGDecTr}\label{eq:GGDecTr}
\end{align*}
We start by estimating the fourth term as follows by using \eqref{eq:GGvi-vi+1}:
\[\abs{\xints{\T^3}{}\2\abs{v_i-v_{i+1}}^2\diff x}\leq\norm{v_i-v_{i+1}}_\ag^2\lsim\Lg\varrho_{q,i}^{1+\gg}\ell_q^{2\ag}. \xtag{GGSqDiff}\label{eq:GGSqDiff}\]
We then proceed to estimate the third term in \eqref{eq:GGDecTr} by using the triangle inequality, \eqref{eq:GMEk}, and the fact $\varrho_{q,i}\sim\varrho_{q,i+1}$:
\begin{align*}
\abs{\xints{\T^3}{}\pa[b]{\abs{v_{\ell_{q,i}}}^2-\abs{v_{\ell_{q,i+1}}}^2}\diff x}\leq{}&\abs{\xints{\T^3}{}\pa[b]{\abs{v_{\ell_{q,i}}}^2-\abs{v_q}^2}\diff x}+\abs{\xints{\T^3}{}\pa[b]{\abs{v_q}^2-\abs{v_{\ell_{q,i+1}}}^2}\diff x}\lsim{} \\
{}\lsim{}&\Lg\varrho_{q,i}^{1+\gg}\ell_q^{2\ag}. \xtag{GGMollDiff}\label{eq:GGMollDiff}\
\end{align*}
The first and second terms in \eqref{eq:GGDecTr} are estimated in similar ways, so we only estimate the former. To that end, we proceed in a way similar to \cite[Proposition 5.5]{DR}. We start by using the fact that $(v_i,p_i,0)$ and $(v_{\ell_{q,i}},p_{\ell_{q,i}},\pint R_{\ell_{q,i}})$ are subsolutions, \kcref{thm:teor:FrLHoeldNorms}, \eqref{eq:GMv} and \eqref{eq:GMR}, and \eqref{eq:GMFracv} and \eqref{eq:GGFracDist}, to obtain that
\begin{align*}
\abs{\Der{t}\xints{\T^3}{}\2\abs{v_i}^2-\abs{v_{\ell_{q,i}}}^2\diff x}\leq{}&2\pa{\abs{\xints{\T^3}{}\2\fD v_{\ell_{q,i}}\mathbin:\pint R_{\ell_{q,i}}\diff x}+\abs{\xints{\T^3}{}\2\abs{\xfrl*v_i}^2-\abs{\xfrl*v_{\ell_{q,i}}}^2\diff x}} \\
{}\lsim{}&\norm{v_{\ell_{q,i}}}_{1+\ag}\norm{\pint R_{\ell_{q,i}}}_\ag+\norm{v_i+v_{\ell_{q,i}}}_{\EXP+\eg}\norm{v_i-v_{\ell_{q,i}}}_{\EXP+\eg} \\
{}\lsim{}&\dg_q^{\xfr12}\lg_q^{1+\ag}\per\Lg\varrho_{q,i}^{1+\gg}\ell_{q,i}^{-\ag}+\Lg^{\xfr12}\per\Lg\varrho_{q,i}^{1+\gg}\tg_q\ell_{q,i}^{-\EXP-\eg-1-\ag} \\
={}&\Lg\varrho_{q,i}^{1+\gg}\tg_q^{-1}\ell_q^{2\ag}(F_1+F_2),
\end{align*}
where we write
\begin{align*}
F_1\coloneq{}&\dg_q^{\xfr12}\lg_q^{1+\ag}\ell_{q,i}^{-\ag}\tg_q\ell_q^{-2\ag} \\
F_2\coloneq{}&\Lg^{\xfr12}\tg_q^2\ell_q^{-2\ag}\ell_{q,i}^{-1-\EXP}\ell_{q,i}^{-\eg-\ag}.
\end{align*}
It is easy to see that
\[F_1=\ell_q^{4\ag}\lg_q^\ag\ell_{q,i}^{-\ag}\ell_q^{-2\ag}\lsim1.\]
To prove that $F_2$ is also bounded above by a constant, we note that
\[F_2=f(\ag,\eg,\gg)\xfr{\Lg^{\xfr12}}{\dg_q\lg_q^2}\xfr{(\zg_q^{\xfr12}\lg_q)^{1+\EXP}}{\varrho_{q,i}^{\xfr{1+\EXP}{2}}},\]
where $f(\ag,\eg,\gg)\sim1$ for $\ag,\gg,\eg\ll1$. We then observe that \eqref{eq:GParRel2}, together with the fact that $\EXP<\xfr13,b>1$, implies that
\[\xfr{\Lg^{\xfr12}}{\dg_q\lg_q^2}\xfr{(\zg_q^{\xfr12}\lg_q)^{1+\EXP}}{\varrho_{q,i}^{\xfr{1+\EXP}{2}}}\leq\varrho_{q,i}^{-\xfr23}(\zg_q^{\xfr12}\lg_q)^{-\xfr23}\lsim\zg_{q+2}^{-\xfr23}\lg_q^{-\xfr23(1-\bg)}=\lg_q^{\xfr43\bg b^2-\xfr23(1-\bg)}=\lg_q^{-\xfr23(1-\bg-2b^2\bg)}\leq1.\]
We have thus proved that
\[F_2\lsim1,\]
provided $\ag,\gg,\eg\ll1$ are sufficiently small and $a$ is sufficiently large. This means that
\[\abs{\Der{t}\xints{\T^3}{}\2\abs{v_i}^2-\abs{v_{\ell_{q,i}}}^2\diff x}\lsim\Lg\varrho_{q,i}^{1+\gg}\tg_q^{-1}\ell_q^{2\ag}. \xtag{GGKinDer1a}\label{eq:GGKinDer1a}\]
Combining \eqref{eq:GGDecTr}-\eqref{eq:GGKinDer1a}, we conclude that
\[\abs{\xints{\T^3}{}\2\abs{\lbar v_q}^2-\abs{v_{\ell_{q,i}}}^2\diff x}\lsim\Lg\varrho_{q,i}^{1+\gg}\ell_q^{2\ag}. \xtag{GGEK1}\label{eq:GGEK1}\]
Estimates \eqref{eq:GGEK1} and \eqref{eq:GMEk} imply
\[\abs{\xints{\T^3}{}\2\abs{\lbar v_q}^2-\abs{v_q}^2\diff x}\lsim\Lg\varrho_{q,i}^{1+\gg}\ell_q^{2\ag}.\]
This proves in particular that $\lbar\rg_q\sim\rg_q$ and \eqref{eq:GGTr}, as well as \eqref{eq:GGTrDist}.
\nipar Similarly, using \eqref{eq:FNSR} for $(v_q,p_q,R_q)$ and $(v_{\ell_{q,i}},p_{\ell_{q,i}},\pint R_{\ell_{q,i}})$ first, and for $(\lbar v_q,\lbar p_q,\lbar R_q)$ and $(v_{\ell_{q,i}},p_{\ell_{q,i}},R_{\ell_{q,i}})$ afterwards, we also deduce that
\begin{align*}
\abs{\Der{t}\xints{\T^3}{}\2\abs{v_{\ell_{q,i}}}^2-\abs{v_q}^2\diff x}\lsim{}&\Lg\varrho_{q,i}^{1+\gg}\tg_q^{-1}\ell_q^{2\ag} \xtag{GGKinDer2}\label{eq:GGKinDer2} \\
\abs{\Der{t}\xints{\T^3}{}\2\abs{\lbar v_q}^2-\abs{v_{\ell_{q,i}}}^2\diff x}\lsim{}&\Lg\varrho_{q,i}^{1+\gg}\tg_q^{-1}\ell_q^{2\ag}. \xtag{GGKinDer3}\label{eq:GGKinDer3}
\end{align*}
Combining \eqref{eq:GGKinDer2} and \eqref{eq:GGKinDer3}, we get
\[\abs{\pd_t\lbar\rg_q-\pd_t\rg_q}\lsim\Lg\varrho_{q,i}^{1+\gg}\tg_q^{-1}\ell_q^{2\ag}. \xtag{GGPreTrDist}\label{eq:GGPreTrDist}\]
To prove \eqref{eq:GGDer} note that
\[\Lg\varrho_{q,i}^{1+\gg}\tg_q^{-1}\ell_q^{2\ag}\lsim\rg_q\dg_q^{\xfr12}\lg_q\per\lg_{q+1}^{2\ag-2\bg\gg}\lsim\rg_q\dg_q^{\xfr12}\lg_q,\]
where we used the definitions of $\tg_q,\zg_{q+1}$, \eqref{eq:GUTr}, the relations $\ell_q^{-1}\leq\lg_{q+1}$ and $\Lg\varrho_{q,i}=\rg_{q,i}\sim\rg_q$, and the fact that $\ag<\ag b<\bg\gg$, which follows from \eqref{eq:GParRel1} since $b>1$. Therefore, since we showed above that $\lbar\rg_q\sim\rg_q$, we have \eqref{eq:GGDer}. \vsp{-3pt}
\nipar It remains to estimate $\norm[b]{\pint{\lbar R}_q}_0$ on $[T_1,T_2]$ in order to verify \eqref{eq:GGStrong} for the Reynolds stress. We already obtained \eqref{eq:GGStrong2} on $J_{\ubar n}\cup\dotso\cup J_{\lbar n}$ (recall the decomposition \eqref{eq:Decomp0T}). Moreover, on $J_{\ubar n-1}\cup J_{\lbar n+1}$ the subsolution remains unchanged, so there is nothing to prove. We are then left with the task of proving \eqref{eq:GGStrong} on the cut-off regions $I_{\ubar n-1}$ and $I_{\lbar n}$.
\nipar To do so, we need to estimate $\|v_i-v_q\|_\ag$ and $\|z_i-z_q\|_\ag$. For the former, we combine estimates of $v_i-v_{\ell_{q,i}}$ and of $v_{\ell_{q,i}}-v_q$. For the latter, we only need to estimate $\|z_{\ell_{q,i}}-z_q\|_\ag$, since we already handled $z_i-z_{\ell_{q,i}}$ above. One has that, by \eqref{eq:Mollif}, \kcref{thm:lemma:Schauder} (Schauder estimates), and \eqref{eq:GUv}
\[\norm{z_{\ell_{q,i}}-z_q}_\ag\lsim\norm{z_q}_{2+\ag}\ell_{q,i}^2\lsim\norm{\opn{curl}v_q}_\ag\ell_{q,i}^2\lsim\dg_q^{\xfr12}\lg_q^{1+\ag}\ell_{q,i}^2=\Lg\tg_q\varrho_{q,i}^{1+\gg}\ell_q^{-2\ag}\lg_q^{-\ag},\]
which gives us \eqref{eq:GGStrong} as desired.
\nipar We then have to verify \eqref{eq:GGFracv} and \eqref{eq:GGRhoD}. To that end, we observe that
\[\norm{\lbar v_q}_{\EXP+\eg}\leq\norm{\lbar v_q-v_q}_{\EXP+\eg}+\norm{v_q}_{\EXP+\eg}.\]
The second term is estimated by \eqref{eq:GUFracv}. As for the first one, we note that
\[\norm{\lbar v_q-v_q}_\ag\lsim\Lg^{\xfr12}\lbar\varrho_q^{\xfr{1+\gg}{2}}\ell_q^\ag\lsim\dg_{q+1}^{\xfr12}\zg_{q+1}^{\xfr\gg2}\ell_q^\ag\leq\dg_{q+1}^{\xfr12}\ell_q^\ag,\]
where the first step is due to \eqref{eq:GGDist2}. We can thus estimate $\|\lbar v_q-v_q\|_{\EXP+\eg}$ by interpolation:
\begin{align*}
\norm{\lbar v_q-v_q}_{\EXP+\eg}\lsim{}&\norm{\lbar v_q-v_q}_\ag^{1-\EXP-\eg}\norm{\lbar v_q-v_q}_{1+\ag}^{\EXP+\eg}\lsim(\dg_{q+1}^{\xfr12}\ell_q^\ag)^{1-\EXP-\eg}\per\pa{\dg_q^{\xfr12}\lg_q^{1+\ag}}^{\EXP+\eg} \\
{}\leq{}&\dg_{q+1}^{\xfr12}\lg_{q+1}^{\EXP+\eg}\per\ell_q^{\ag(1-2\EXP-2\eg)}.
\end{align*}
Lastly
\[\abs{\pd_t\s{T}_g}=\abs{\xints{\T^3}{}\2\abs{\xfrl*\lbar v_q}^2-\abs{\xfrl*v_q}^2\diff x}\lsim\norm{\lbar v_q+v_q}_{\EXP+\eg}\norm{\lbar v_q-v_q}_{\EXP+\eg}\lsim\Lg^{\xfr12}\dg_{q+1}^{\xfr12}\lg_{q+1}^{\EXP+\eg},\]
thus proving \eqref{eq:GGRhoD}.
\nipar We conclude the proof by obtaining \eqref{eq:GGDistA2S}. To that end, we first note that \eqref{eq:GGTr} and \eqref{eq:GUTr} combine to give us $\lbar\varrho_q\lsim\zg_{q+1}$. Combining this with \eqref{eq:GGDist}, we conclude that
\[\norm{\lbar v_q-v_q}_\ag\lsim\dg_{q+1}^{\xfr12}\zg_{q+1}^{\xfr\gg2}\ell_q^\ag,\]
thus reducing \eqref{eq:GGDistA2S} to
\[\lg_{q+1}^{-\bg\gg}\ell_q^{(\xfr12-\xfr2b)\ag}\lsim1.\]
Since $\ell_q^{-1}\leq\lg_{q+1}$, the above follows from
\[-\bg\gg+\xfr2b\ag-\xfr12\ag<0\iff\ag<\xfr{2b\bg\gg}{4-b}.\]
We recall that we wish to obtain \eqref{eq:GGDistA2S} only under the assumption $2\ag<\bg\gg$. This means the above relation follows from
\[\xfr{2b}{4-b}>\xfr12\iff5b>4,\]
which in turn follows from $b>1$. The proof is thus complete.
\end{qeddim*}

\xbegin{oss}[Multi-gluing][thm:oss:MultiGluing]
\kcref{thm:propo:Gluing} can easily be extended to a pairwise disjoint union of intervals $[T_1\stp i,T_2\stp i]\sbs[0,T]$ with $T_2\stp i-T_1\stp i\geq4\tg_q$ and $T_2\stp i<T_1\stp{i+1}$.
\xend{oss}%

\sect{Perturbation step}\label{MPS}
\xbegin{propo}[Main Perturbation Step][thm:propo:MPS]
Let $b,\bg,\ag,\gg,(\dg_q,\lg_q,\Lg,\zg_q,\ell_q,\tg_q)$ be as in Section \ref{Equestimates} with
\begin{align*}
\ag<{}&\bg\gg. \xtag{PParam1}\label{eq:PParam1}
\end{align*}
Let $[T_1,T_2]\sbs[0,T]$ and let $t_i,I_i,J_i$ be as in \eqref{eq:DeftiIi}. Let $(v,p,R)$ be a smooth strong subsolution on $[T_1,T_2]$ satisfying
\begin{align*}
\norm{R}_1\lsim{}&\Lg\varrho^{1+\gg}\ell_q^{-2\ag}\ell_{q,i}^{-1} \xtag{PRStrong}\label{eq:PRStrong} \\
\dg_{q+2}\lsim{}&\rg\lsim\dg_{q+1} \xtag{PRho}\label{eq:PRho} \\
\norm{v}_0\leq{}&C_P\xtag{Pv00}\label{eq:Pv00}
\end{align*}
on $K_i\coloneq[(i-1+\xfr13)\tg_q,(i+\xfr23)\tg_q]$ for $\ubar n-1\leq i\leq\lbar n+1$, where the $\ell_{q,i}$ are defined as in \eqref{eq:rgqiellqi}, and $C_P$ is a geometric constant. Further, let $\yg:[0,T]\to[0,1]$ be a cutoff function and $S_\yg\in\Cinf(\T^3\x[T_1,T_2];\s{S}^{3\x3})$ be a smooth matrix field with
\[S_\yg(x,t)=\sg_\yg(t)\opn{Id}+\pint S_\yg(x,t)=\Lg\varsigma_\yg(t)\opn{Id}+\pint S_\yg(x,t), \xtag{DecompS}\label{eq:DecompS}\]
where $\pint S_\yg=\yg^2\pint S$ is traceless, $\sg_\yg=\yg^2\sg$, and $(\varsigma,\varsigma_\yg)\coloneq\Lg^{-1}(\sg,\sg_\yg)$. Suppose $\yg$ satisfies
\begin{align*}
\abs{\yg'}\lsim{}&\dg_q^{\xfr12}\lg_q, \xtag{PCutDer}\label{eq:PCutDer} \\
\intertext{and $\sg$ satisfies}
0\leq{}&\sg(t)\leq4\dg_{q+1} \xtag{PSTr}\label{eq:PSTr} \\
\sg|_{K_i}\lsim{}&\rg(t_i) \xtag{PSigmaRho}\label{eq:PSigmaRho} \\
\abs{\pd_t\sg}\lsim{}&\sg\dg_q^{\xfr12}\lg_q, \xtag{PSTrDer}\label{eq:PSTrDer} \\
\intertext{Moreover, assume that for any $N\geq0,t\in I_{i-1}\cup J_i\cup I_i,\ubar n\leq i\leq\lbar n$}
\norm{\pint S}_{N+\ag}\lsim{}&\Lg\varsigma^{1+\gg}\ell_{q,i}^{-N}\ell_q^{-2\ag} \xtag{PSStrong}\label{eq:PSStrong} \\
\norm{v}_{N+1+\ag}\lsim{}&\dg_q^{\xfr12}\lg_q^{1+\ag}\ell_{q,i}^{-N} \xtag{Pv0}\label{eq:Pv0} \\
\norm{v}_{\EXP+\eg}\leq{}&M\pa{1+\sum_{i=0}^{q+1}\Lg^{\xfr12}\lg_i^{\EXP+\eg-\bg}} \xtag{PFracv0}\label{eq:PFracv0} \\
\norm{(\pd_t+v\per\grad)\pint S}_{N+\ag}\lsim{}&\Lg\varsigma^{1+\gg}\ell_{q,i}^{-N}\ell_q^{-6\ag}\dg_q^{\xfr12}\lg_q, \xtag{PSMatDer}\label{eq:PSMatDer}
\end{align*}
Finally, assume that
\[\opn{supp}\pint S_\yg\sbse\T^3\x\bigcup_iI_i. \xtag{PSupp}\label{eq:PSupp}\]
Then, provided $a\dgreat1$ is sufficiently large (depending on the implicit constants in \eqref{eq:PSTrDer}, \eqref{eq:PSStrong}, \eqref{eq:Pv0}, and \eqref{eq:PSMatDer}), there exist smooth $(\tilde v,\tilde p)\in\Cinf(\T^3\x[T_1,T_2];\R^3\x\R)$ and a smooth matrix field $\s{E}\in\Cinf(\T^3\x[T_1,T_2];\s{S}^{3\x3})$ with $\opn{supp}\s{E}\sbs\opn{supp}S$ such that, setting $\tilde R_1\coloneq R-S_\yg-\s{E}$, the triple $(\tilde v,\tilde p,\tilde R)$ is a strong subsolution with
\[\xints{\T^3}{}\2|\tilde v|^2+\opn{tr}\tilde R\diff x=\xints{\R^3}{}\2|v|^2+\opn{tr}R\diff x\qquad\VA t. \xtag{PPIC}\label{eq:PPIC}\]
Moreover, we have the estimates
\begin{align*}
\norm{\tilde v-v}_{H^{-1}}\leq&\xfr M2\dg_{q+1}^{\xfr12}\ell_{q,i}^{-1}\lg_{q+1}^{-1} \xtag{PDistH-1}\label{eq:PDistH-1} \\
\norm{\tilde v-v}_0\leq{}&\xfr M2\dg_{q+1}^{\xfr12} \xtag{PDist0}\label{eq:PDist0} \\
\norm{\tilde v-v}_{1+\ag}\leq{}&\xfr M2\dg_{q+1}^{\xfr12}\lg_{q+1}^{1+\ag} \xtag{PDist1}\label{eq:PDist1} \\
\norm{\tilde v-v}_{\EXP+\eg}\leq{}&M\dg_{q+1}^{\xfr12}\lg_{q+1}^{\EXP+\eg}, \xtag{PFracDist}\label{eq:PFracDist} \\
\intertext{and the error $\s{E}$ satisfies the estimates}
\norm{\s{E}}_0\leq{}&\dg_{q+2}\lg_{q+1}^{-6\ag} \xtag{PErr}\label{eq:PErr} \\
\abs{\pd_t\opn{tr}\s{E}}\leq{}&\dg_{q+2}\dg_{q+1}^{\xfr12}\lg_{q+1}^{1-6\ag}. \xtag{PErrDer}\label{eq:PErrDer}
\end{align*}
Finally, setting
\begin{align*}
\s{T}_p(t)\coloneq{}&\xfr13\xints0t\xints{\T^3}{}\pa[B]{\abs{\xfrl*\tilde v}^2-\abs{\xfrl*v}^2}\diff x\diff s \\
{}={}&\xfr13\xints0t\xints{\T^3}{}\pa[B]{\xfrl*(\tilde v+v)\per\xfrl*(\tilde v-v)}\diff x\diff s,
\end{align*}
we have that
\[\abs{\pd_t\s{T}_p}\lsim\Lg\lg_q^{\EXP+\eg-\bg}, \xtag{PHypoTrace}\label{eq:PHypoTrace}\]
for any $\eg>0$. Thus, $\s{T}_p\opn{Id}$ satisfies \eqref{eq:PErrDer} for all $t\in[0,T]$, and \eqref{eq:PErr} only for small times.
\xend{propo}
\begin{qeddim}[The proof extends \cite[Section 7]{DRSz}, which is a localization of the argument carried out in \cite[Section 5]{BDLSzV}. The difference between \cite{DRSz} and \cite{BDLSzV} is that the latter absorbs the whole $R$ with the perturbation flow, whereas the former, as well as the proof below, aims to only absorb $S$.]
\begin{center}
\tbold{Step 1: Squiggling Stripes and the Stress Tensors $\tilde S_i$}
\end{center}
As in \cite[Lemma 5.3]{BDLSzV}, we choose a family of smooth non-negative $\hg_i=\hg_i(x,t)$ with the following properties:
\begin{align*}
\hg_i\in{}&\relax\Cinf(\T^3\x[T_1,T_2];[0,1]) \xtag{$\hg_i$-i}\label{eq:hgi1} \\
\opn{supp}\hg_i\cap\opn{supp}\hg_j={}&\0\qquad\qquad\VA i\neq j \xtag{$\hg_i$-ii:DisjSupp}\label{eq:hgi2:DisjSupp} \\
\T^3\x I_i\sbs{}&\br{(x.t):\hg_i(x,t)=1} \xtag{$\hg_i$-iii:Restr$I_i$}\label{eq:hgi3:RestrIi} \\
\opn{supp}\hg_i\sbse{}&\T^3\x(J_i\cup I_i\cup J_{i+1}) \xtag{$\hg_i$-iv:Supp}\label{eq:hi4:Supp} \\
{}={}&\T^3\x\br{\pa{t_i-\xfr13\tg_q,t_i+\xfr13\tg_q}\cap[0,T]} \\
\3c_0>0:\quad\sum_i\xints{\T^3}{}\2\hg_i^2(x,t)\diff x\geq{}&c_0\quad\VA t\in[0,T] \xtag{$\hg_i$-v:$c_0$}\label{eq:hgi5:c0} \\
\norm{\pd_t^N\hg_i}_m\leq{}&C(N,m)\tg_q^{-N}\qquad N,m\geq0, \xtag{$\hg_i$-vi:NormDer}\label{eq:hgi6:NormDer}
\end{align*}
where the $c_0$ in \eqref{eq:hgi5:c0} is a geometric constant. Define
\[\sg_i(x,t)\coloneq\abs{\T^3}\xfr{\hg_i^2(x,t)}{\sum_j\int\hg_j^2(y,t)\diff y}\sg_\yg(t),\]
so that $\sum_i\int_{\T^3}\sg_i\diff x=\int_{\T^3}\sg_\yg\diff x$. Using the inverse flow $\Fg_i$ starting at time $t_i$
\[\begin{sistema}
(\pd_t+v\per\grad)\Fg_i=0 \\
\Fg_i(x,t_i)=x
\end{sistema},\]
set
\begin{align*}
S_i\coloneq{}&\sg_i\opn{Id}+\hg_i^2\pint S_\yg \\
\tilde S_i\coloneq{}&\xfr{\fD\Fg_iS_i\fD^T\Fg_i}{\sg_i}=\fD\Fg_i\pa{\opn{Id}+\xfr{\sum_j\int\hg_j^2}{|\T^3|\sg}\pint S}\fD^T\Fg_i.
\end{align*}
One can check from \eqref{eq:hgi5:c0}-\eqref{eq:hgi6:NormDer} and from \eqref{eq:PSTr} and \eqref{eq:PSigmaRho} that
\begin{align*}
\norm{\sg_i}_0\leq{}&4|\T^3|c_0^{-1}\dg_{q+1} \xtag{SigmaI0}\label{eq:SigmaI0} \\
\norm{\sg_i}_N\lsim{}&\rg_i\coloneq\rg(t_i)\lsim\dg_{q+1}, \xtag{SigmaI}\label{eq:SigmaI}
\end{align*}
and moreover, since by \eqref{eq:PSupp} $\opn{supp}\pint S_\yg\sbse\br{\sum_i\hg_i^2=1}$
\[\xfr13\opn{tr}\sum_i\xints{\T^3}{}\2S_i\diff x=\xfr13\opn{tr}S_\yg. \xtag{P.$\sg_i$-$\sg$}\qquad\qquad\sum_i\pint S_i=\pint S_\yg\label{eq:Psgisg}\]
We next claim that for all $(x,t)$
\[\tilde S_i(x,t)\in B_{\xfr12}(\opn{Id})\sbs\s{S}_+^{3\x3}, \xtag{P\~SBall}\label{eq:PStildeBall}\]
where $B_{\xfr12}(\opn{Id})$ is the ball of radius $\xfr12$ centered at the identity $\opn{Id}$ in $\s{S}^{3\x3}$. Indeed, by the classical estimates on transport equations reported in \kcref{thm:propo:EstTrasp}
\[\norm{\grad\Fg_i-\opn{Id}}_0\lsim\tg_q\dg_q^{\xfr12}\lg_q^{1+\ag}=\ell_q^{4\ag}\lg_q^\ag\leq\ell_q^{3\ag} \xtag{PDGrad}\label{eq:PDGrad}\]
for $t\in J_i\cup I_i\cup J_{i+1}$, since this is an interval of length $\abs{J_i\cup I_i\cup J_{i+1}}\sim\tg_q$. Using \eqref{eq:PSTr}, \eqref{eq:PSStrong} and \eqref{eq:LambdaEll}, we also have that, for any $N\geq0$
\[\norm{\xfr{\hg_i^2\pint S_\yg}{\sg_i}}_N\lsim\norm{\xfr{\pint S}{\sg}}_N\lsim\varsigma^\gg\ell_{q,i}^{-N}\ell_q^{-2\ag}\lsim\zg_{q+1}^\gg\lg_{q+1}^{2\ag}\ell_{q,i}^{-N}=\lg_{q+1}^{2\ag-2\bg\gg}\ell_{q,i}^{-N}. \xtag{PcutN}\label{eq:PcutN}\]
Then, using the decomposition
\[\tilde S_i-\opn{Id}=\fD\Fg_i\xfr{\hg_i^2\pint S_\yg}{\sg_i}\fD\Fg_i^T+\fD\Fg_i(\fD\Fg_i^T-\opn{Id})+\fD\Fg_i-\opn{Id},\]
we deduce from \eqref{eq:PDGrad}-\eqref{eq:PcutN} that
\[\abs{\tilde S_i-\opn{Id}}\lsim(1+\ell_q^{3\ag})\lg_{q+1}^{2\ag-2\bg\gg}(1+\ell_q^{3\ag})+2\ell_q^{3\ag}\leq\xfr12,\]
provided $a\dgreat1$ is sufficiently large, since we assumed $\ag<\bg\gg$ in \eqref{eq:PParam1}. This verifies \eqref{eq:PStildeBall}.
\begin{center}
\tbold{Step 2: The perturbation $w$.}
\end{center}
Now we can define the perturbation term as
\[w_o\coloneq\sum_i\rad{\sg_i}(\fD\Fg_i)^{-1}W(\tilde S_i,\lg_{q+1}\Fg_i)=\sum_iw_{oi},\]
where $W$ are the Mikado flows on the compact set $B_{\xfr12}(\opn{Id})$ as defined in \kcref{thm:lemma:MikFlows}. Notice that the supports of the $w_{oi}$ are disjoint and, using the Fourier series representation of the Mikado flows
\[w_{oi}\coloneq\sum_{k\neq0}(\fD\Fg_i)^{-1}b_{i,k}A_ke^{i\lg_{q+1}k\per\Fg_i}, \xtag{Pwoi}\label{eq:Pwoi}\]
where we write
\[b_{i,k}(x,t)\coloneq\rad{\sg_i(x,t)}a_k(\tilde S_i(x,t)\!).\]
We define $w_c$ so that $w\coloneq w_o+w_c$ is divergence-free:
\[w_c\coloneq\xfr{i}{\lg_{q+1}}\per\e\sum_{i,k\neq0}\fD(b_{i,k})\x\xfr{\fD\Fg_i^T\per(k\x A_k)}{|k|^2}e^{i\lg_{q+1}k\per\Fg_i}=\sum_{i,k\neq0}\xfr{c_{i,k}}{\lg_{q+1}}e^{i\lg_{q+1}k\per\Fg_i},\]
where we write
\[c_{i,k}= \fD(b_{i,k})\x\xfr{\fD\Fg_i^T\per(k\x A_k)}{|k|^2}. \xtag{PDefcik}\label{eq:PDefcik}\]
Define then
\begin{align*}
w\coloneq{}&w_o+w_c \\
\tilde v\coloneq{}&v+w \\
\tilde p\coloneq{}&p-\sum_i\sg_i \\
\s{E}(x,t)\coloneq{}&\pint{\s{E}}\stp1(x,t)+\s{E}\stp2,
\end{align*}
where
\[\pint{\s{E}}\stp1\coloneq\s{R}(\pd_t\tilde v+\opn{div}(\tilde v\otimes\tilde v)+\grad\tilde p+\xfrl\tilde v+\opn{div}(R-S_\yg)\!), \xtag{PDefE1}\label{eq:PDefE1}\]
with $\s{R}$ being the anti-divergence operator defined in \kcref{thm:defi:AntiDivR}, and
\[\s{E}\stp2(t)\coloneq\xfr{\opn{Id}}{3}\xints{\T^3}{}\pa[b]{|\tilde v|^2-|v|^2-\opn{tr}S_\yg}\diff x. \xtag{PDefE2}\label{eq:PDefE2}\]
Equations \eqref{eq:PPIC} and \eqref{eq:FNSR} follow by construction.
\begin{center}
\tbold{Step 3: Estimates on the perturbation}\label{PStep3}
\end{center}
The estimates on $\tilde v$ follow similarly to the ones for $v_{q+1}$ in \cite[Section 5-6]{BDLSzV}. Obtaining those requires estimates on the coeffcients $b_{i,k},c_{i,k}$, which in turn require estimates of $\tilde S_i$ and estimates of $\fD\Fg_i$. The latter read as follows:
\begin{align*}
\norm{\fD\Fg_i-\opn{Id}}_N+\norm{(\fD\Fg_i)^{-1}-\opn{Id}}_N\lsim{}&\ell_q^{3\ag}\ell_{q,i}^{-N} \xtag{PDistGrad}\label{eq:PDistGrad} \\
\norm{\fD\Fg_i}_N+\norm{(\fD\Fg_i)^{-1}}_N\lsim{}&\ell_q^{3\ag\one_{N\neq0}}\ell_{q,i}^{-N} \xtag{PGrad}\label{eq:PGrad} \\
\norm{(\pd_t+v\per\grad)\fD\Fg_i}_N\lsim{}&\dg_q^{\xfr12}\lg_q^{1+\ag}\ell_q^{3\ag\one_{N\neq0}}\ell_{q,i}^{-N}. \xtag{PMDGrad}\label{eq:PMDGrad}
\end{align*}
To obtain these, we first observe that $\Fg_i$ is a diffeomorphism, which implies both $\fD\Fg_i$ and $(\fD\Fg_i)^{-1}$ are bounded, thus yielding the $N=0$ case of \eqref{eq:PGrad}. To obtain \eqref{eq:PDistGrad}, we start by combining \eqref{eq:PDGrad} with the $N=0$ case of \eqref{eq:PGrad}, thus obtaining that
\[\norm{(\fD\Fg_i)^{-1}-\opn{Id}}_0\leq\norm{(\fD\Fg_i)^{-1}}_0\norm{\opn{Id}-\fD\Fg_i}_0\lsim\ell_q^{3\ag}.\]
This yields \eqref{eq:PDistGrad} for $N=0$. For $N\geq1$, we note that
\[\norm{\fD\Fg_i-\opn{Id}}_N+\norm{(\fD\Fg_i)^{-1}-\opn{Id}}_N\lsim\norm{\fD\Fg_i-\opn{Id}}_0+\norm{(\fD\Fg_i)^{-1}-\opn{Id}}_0+\norm{\fD^2\Fg_i}_{N-1}+\norm{\fD(\fD\Fg_i)^{-1}}_{N-1}.\]
The other cases of \eqref{eq:PDistGrad} follow by combining its $N=0$ case with the $N\geq1$ cases of \eqref{eq:PGrad}.
\nipar The estimates for $\|\fD\Fg_i\|_N$ for $N\geq1$ follow from \kcref{thm:propo:EstTrasp}. By combining \eqref{eq:MDGrad} with \kcref{thm:lemma:HNI}, we obtain that
\[\norm{(\pd_t+v\per\grad)\fD\Fg_i}_N\lsim\norm{\fD\Fg_i}_N\norm{\grad v}_0+\norm{\grad v}_N\norm{\fD\Fg_i}_0.\]
Estimates \eqref{eq:PGrad} and \eqref{eq:Pv0} then yield \eqref{eq:PMDGrad}.
\nipar To complete the proof of \eqref{eq:PGrad}, we are left with estimating $\|(\fD\Fg_i)^{-1}\|_N$. We note that $\fD^N(\Fg_i\circ\Fg_i^{-1})=0$ for $N\geq1$. We then use the Leibniz rule and the chain rule to write
\begin{align*}
\fD^2(\Fg_i\circ\Fg_i^{-1})={}&\fD(\!(\fD\Fg_i\circ\Fg_i^{-1})\fD\Fg_i^{-1})=(\fD^2\Fg_i\circ\Fg_i^{-1})(\fD\Fg_i^{-1})^2+(\fD\Fg_i\circ\Fg_i^{-1})\fD^2\Fg_i^{-1} \\
\fD^3(\Fg_i\circ\Fg_i^{-1})={}&\fD(\fD^2(\Fg_i\circ\Fg_i^{-1})\!) \\
{}={}&(\fD^3\Fg_i\circ\Fg_i^{-1})(\fD\Fg_i^{-1})^3+3(\fD^2\Fg_i\circ\Fg_i^{-1})(\fD\Fg_i^{-1})(\fD^2\Fg_i^{-1})+(\fD\Fg_i\circ\Fg_i^{-1})\fD^3\Fg_i^{-1} \\
\fD^4(\Fg_i\circ\Fg_i^{-1})={}&\fD(\fD^3(\Fg_i\circ\Fg_i^{-1})\!) \\
{}={}&(\fD^4\Fg_i\circ\Fg_i^{-1})(\fD\Fg_i^{-1})^4+6(\fD^3\Fg_i\circ\Fg_i^{-1})(\fD\Fg_i^{-1})^2\fD^2\Fg_i^{-1}+3(\fD^2\Fg_i\circ\Fg_i^{-1})(\fD^2\Fg_i^{-1})^2 \\
&{}+4(\fD^2\Fg_i\circ\Fg_i^{-1})(\fD\Fg_i^{-1})(\fD^3\Fg_i^{-1})+(\fD\Fg_i\circ\Fg_i^{-1})\fD^4\Fg_i^{-1}.
\end{align*}
From these, we can see that
\begin{align*}
\norm{\fD^2\Fg_i^{-1}}_0\leq{}&\norm{\fD\Fg_i^{-1}}_0\norm{\fD^2\Fg_i}_0\norm{\fD\Fg_i^{-1}}_0^2\lsim\ell_q^{3\ag}\ell_{q,i}^{-1} \\
\norm{\fD^3\Fg_i^{-1}}_0\leq{}&\norm{\fD\Fg_i^{-1}}_0\pa{\norm{\fD^3\Fg_i}_0\norm{\fD\Fg_i^{-1}}_0^3+3\norm{\fD^2\Fg_i}_0\norm{\fD\Fg_i^{-1}}_0\norm{\fD^2\Fg_i^{-1}}_0}\lsim\ell_q^{3\ag}\ell_{q,i}^{-2} \\
\norm{\fD^4\Fg_i^{-1}}_0\leq{}&\norm{\fD\Fg_i^{-1}}_0\left(\norm{\fD^4\Fg_i}_0\norm{\fD\Fg_i^{-1}}_0^4+6\norm{\fD^3\Fg_i}_0\norm{\fD\Fg_i^{-1}}_0^2\norm{\fD^2\Fg_i^{-1}}_0\right. \\
&{}+\left.3\norm{\fD^2\Fg_i}_0\norm{\fD^2\Fg_i^{-1}}_0^2+4\norm{\fD^2\Fg_i}_0\norm{\fD\Fg_i^{-1}}_0\norm{\fD^3\Fg_i^{-1}}_0\right)\lsim\ell_q^{3\ag}\ell_{q,i}^{-3}.
\end{align*}
These two examples show us that
\[\fD^M(\Fg_i\circ\Fg_i^{-1})=(\fD\Fg_i\circ\Fg_i^{-1})\fD^M\Fg_i^{-1}+\text{other terms},\]
where the other terms are of the form $(\fD^k\Fg_i\circ\Fg_i^{-1})(\fD\Fg_i^{-1})^\ell(\fD^m\Fg_i^{-1})^n$, where $k-1+(m-1)n=M-1$ and $m<M$. If we assume \eqref{eq:PGrad} for $N<M-1$, we see that such terms are estimated as $\ell_{q,i}^{-(k-1)}\ell_{q,i}^{-(m-1)n}=\ell_{q,i}^{-(1-M)}$, thus so is $\fD^M\Fg_i^{-1}=\fD(\fD^{M-1}\Fg_i^{-1})$, which proves \eqref{eq:PGrad} for $N=M-1$. Thus, by induction, the estimate \eqref{eq:PGrad} is proved.
\nipar The following estimates then follow precisely as in \cite[Propositions 5.7 and 5.9]{BDLSzV}:
\begin{align*}
\norm{\tilde S_i}_N\lsim{}&\ell_{q,i}^{-N} \xtag{PStildeN}\label{eq:PStildeN} \\
\norm{b_{i,k}}_N\lsim{}&\rg_i^{\xfr12}|k|^{-6}\ell_{q,i}^{-N} \xtag{Pbik}\label{eq:Pbik} \\
\norm{c_{i,k}}_N\lsim{}&\rg_i^{\xfr12}|k|^{-6}\ell_{q,i}^{-N-1} \xtag{Pcik}\label{eq:Pcik} \\
\norm{D_t\tilde S_i}_N\lsim{}&\tg_q^{-1}\ell_{q,i}^{-N} \xtag{PMDStilde}\label{eq:PMDStilde} \\
\norm{D_tb_{i,k}}_N\lsim{}&\dg_{q+1}^{\xfr12}\tg_q^{-1}\ell_{q,i}^{-N}|k|^{-6} \xtag{PMDbik}\label{eq:PMDbik} \\
\norm{D_tc_{i,k}}_N\lsim{}&\dg_{q+1}^{\xfr12}\tg_q^{-1}\ell_{q,i}^{-N-1}|k|^{-6}. \xtag{PMDcik}\label{eq:PMDcik}
\end{align*}
To obtain \eqref{eq:PStildeN}, observe first that, by its definition, we have that
\[\tilde S_i=\fD\Fg_i\fD^T\Fg_i+\fD\Fg_i\xfr{\sum_j\xints{\T^3}{}\2\hg_j^2\diff x}{|\T^3|\sg}\pint S\fD^T\Fg_i,\]
and therefore
\[\norm{\tilde S_i}_N\lsim\norm{\fD\Fg_i}_N+\norm{\fD\Fg_i}_N\norm{\xfr{\pint S}{\sg}}_0+\norm{\xfr{\pint S}{\sg}}_N,\]
where we used that $\|\fD\Fg_i\|_0\lsim1$. By \eqref{eq:PGrad}, the first term above obeys \eqref{eq:PStildeN}. To estimate the remaining two terms, we use \eqref{eq:PSStrong} and \eqref{eq:PSTr} to obtain that
\[\norm{\xfr{\pint S}{\sg}}_N\lsim\xfr{\norm{\pint S}_{N+\ag}}{\sg}\lsim\varsigma^\gg\ell_{q,i}^{-N}\ell_q^{-2\ag}\lsim\ell_{q,i}^{-N}\lg_{q+1}^{2\ag-2\bg\gg}. \sxtag{PS/sg}\]
Estimate \eqref{eq:PStildeN} then follows from \eqref{eq:PGrad} and the assumption \eqref{eq:PParam1}, i.e. that $\ag<\bg\gg$. The proof of \eqref{eq:PMDStilde} follows a similar strategy. First, we decompose $S_i\sg_i^{-1}$ and its material derivative as
\begin{align*}
\xfr{S_i}{\sg_i}=\opn{Id}+\xfr{\sum_j\xints{\T^3}{}\2\hg_j^2\diff x}{|\T^3|}\xfr{\pint S}{\sg}\eqcolon\opn{Id}+f(t)\xfr{\pint S}{\sg} \\
D_t\xfr{S_i}{\sg_i}=f'(t)\xfr{\pint S}{\sg}+f(t)\xfr{D_t\pint S}{\sg}-f(t)\xfr{\sg'}{\sg}\xfr{\pint S}{\sg}. \sxtag{PDecompMDSi/sgi}
\end{align*}
We then decompose $D_t\tilde S_i$ as follows:
\[D_t\tilde S_i=D_t\fD\Fg_i\xfr{S_i}{\sg_i}\fD^T\Fg_i+\fD\Fg_iD_t\xfr{S_i}{\sg_i}\fD^T\Fg_i+\fD\Fg_i\xfr{S_i}{\sg_i}D_t\fD^T\Fg_i\eqcolon I+II+III.\]
The terms $I$ and $III$ are estimated in similar ways, so we only estimate $I$. To that end, we note that, using \eqref{eq:PS/sg} and the fact that $f\lsim1$, we can obtain that
\[\norm{\xfr{S_i}{\sg_i}}_N\leq1+\norm{\xfr{\pint S}{\sg}}_N\lsim\ell_{q,i}^{-N}.\]
Using this bound together with \eqref{eq:PGrad} and \eqref{eq:PMDGrad}, we obtain that
\begin{align*}
\norm{I}_N\lsim{}&\norm{D_t\fD\Fg_i}_N\norm{\xfr{S_i}{\sg_i}}_0\norm{\fD\Fg_i}_0+\norm{D_t\fD\Fg_i}_0\norm{\xfr{S_i}{\sg_i}}_N\norm{\fD\Fg_i}_0+\norm{D_t\fD\Fg_i}_0\norm{\xfr{S_i}{\sg_i}}_0\norm{\fD\Fg_i}_N \\
{}\lsim{}&\dg_q^{\xfr12}\lg_q^{1+\ag}\ell_q^{3\ag\one_{N\neq0}}\ell_{q,i}^{-\ag}\ell_q^{3\ag\one_{N\neq0}}\ell_{q,i}^{-N}.
\end{align*}
Combining the above with
\[\dg_q^{\xfr12}\lg_q^{1+\ag}=\tg_q^{-1}\ell_q^{3\ag}\leq\tg_q^{-1}, \sxtag{PDgTg}\]
we see that $I$ satisfies \eqref{eq:PMDStilde}.
\nipar Coming to $II$, we first note that
\[\norm{II}_N\lsim\norm{\fD\Fg_i}_N\norm{D_t\xfr{S_i}{\sg_i}}_0\norm{\fD\Fg_i}_0+\norm{\fD\Fg_i}_0^2\norm{D_t\xfr{S_i}{\sg_i}}_N. \sxtag{PMDSi/sgi.II}\]
Therefore, to prove that $II$ satisfies \eqref{eq:PMDStilde}, we need to estimate $D_t(\sg_i^{-1}S_i)$. Recalling \eqref{eq:PDecompMDSi/sgi}, we note that
\[\norm{D_t\xfr{S_i}{\sg_i}}_N\leq\abs{f'(t)}\norm{\xfr{\pint S}{\sg}}_N+\xfr{f(t)}{\sg(t)}\norm{D_t\pint S}_N+\xfr{f(t)\sg'}{\sg^2}\norm{\pint S}_N\eqcolon P_1+P_2+P_3. \sxtag{PEstMDSi/sgi}\]
To estimate $P_1$, we first note that, by \eqref{eq:hgi6:NormDer}
\[f'(t)=\xfr{\sum_i\xints{\T^3}{}\2\pd_t\hg_i(x,t)\hg_i(x,t)\diff x}{|\T^3|}\lsim\tg_q^{-1}.\]
Combining this with \eqref{eq:PS/sg}, we conclude that
\[\norm{P_1}_N\lsim\tg_q^{-1}\ell_{q,i}^{-N}. \sxtag{PP1}\]
To bound $P_2$ from above, we recall that $f\lsim1$, and apply \eqref{eq:PSMatDer} to conclude that
\[\norm{P_2}_N\lsim\Lg\varsigma^\gg\ell_{q,i}^{-N}\ell_q^{-6\ag}\dg_q^{\xfr12}\lg_q.\]
We then note that
\[\varsigma^\gg\ell_q^{-6\ag}\dg_q^{\xfr12}\lg_q\lsim\tg_q^{-1}\lg_{q+1}^{2\ag-2\bg\gg}, \sxtag{PSorrata}\]
so that $P_2$ also satisfies \eqref{eq:PP1}.
\nipar Coming, finally, to $P_3$, we use \eqref{eq:PSTrDer}, \eqref{eq:PS/sg}, and \eqref{eq:PDgTg} to conclude that $P_3$ also satisfies \eqref{eq:PP1}. Combining the bounds on $P_1,P_2,P_3$ we have just obtained, we conclude that $D_t(\sg_i^{-1}S_i)$ also satisfies \eqref{eq:PP1}, i.e.
\[\norm{D_t\xfr{S_i}{\sg_i}}_N\lsim\tg_q^{-1}\ell_{q,i}^{-N}.\]
Combining this estimate with \eqref{eq:PGrad}, we conclude \eqref{eq:PMDStilde}.
\nipar To prove \eqref{eq:Pbik} and \eqref{eq:PMDbik}, we first prove some estimates on $\rad{\sg_i}$. Firstly, we note that, thanks to \eqref{eq:hgi5:c0} and \eqref{eq:PSigmaRho}
\[\norm{\rad{\sg_i}}_N\leq\xfr{|\T^3|^{\xfr12}}{\rad{C_0}}\rad{\sg_\yg}\norm{\hg_i}_N\leq\xfr{|\T^3|^{\xfr12}\yg}{C_0^{\xfr12}}\rg_i^{\xfr12}\lsim\rg_i^{\xfr12}. \xtag{P$\rad{\sg_i}$}\label{eq:PRadSgi}\]
As for the material derivative, setting
\[h(t)\coloneq\sum_j\int\hg_j^2(x,t)\diff x,\]
we have that
\[\sg_i=\xfr{|\T^3|\yg^2\hg_i^2}{h}\sg.\]
We can thus write
\begin{align*}
D_t\rad{\sg_i}={}&\xfr12\sg_i^{-\xfr12}D_t\sg_i=\xfr{\rad{h}}{2|\T^3|^{\xfr12}\yg\rad\sg\hg_i}\pa{\xfr{2|\T^3|\yg^2\sg}{h}\hg_iD_t\hg_i+|\T^3|\pd_t\pa{\xfr{\yg^2\sg}{h}}\hg_i^2} \\
{}={}&\xfr{|\T^3|^{\xfr12}\yg\rad{\sg}}{\rad{h}}D_t\hg_i+\xfr{|\T^3|^{\xfr12}\rad\sg\hg_i}{\rad{h}}\pd_t\yg+|\T^3|^{\xfr12}\xfr{1}{2\rad{h\sg}}\pd_t\sg\yg\hg_i-|\T^3|^{\xfr12}\xfr{\rad\sg}{2h^{\xfr32}}\pd_th\yg\hg_i \\
{}\eqcolon{}&I+II+III+IV.
\end{align*}
We then note that, by \eqref{eq:hgi6:NormDer} and \eqref{eq:Pv0}, we have that
\[\norm{D_t\hg_i}_N\leq\norm{\pd_t\hg_i}_N+\norm{v}_N\norm{\hg_i}_0+\norm{v}_0\norm{\hg_i}_N\lsim\tg_q^{-1}+\dg_q^{\xfr12}\lg_q^{1+\ag}\ell_{q,i}^{-N}\lsim\tg_q^{-1}\ell_{q,i}^{-N}.\]
Combining this with \eqref{eq:PSigmaRho}, \eqref{eq:hgi5:c0}, and the fact $\yg\leq1$, we easily see that
\[\norm{I}_N\leq\xfr{|\T^3|^{\xfr12}\rg_i^{\xfr12}}{c_0^{\xfr12}}\norm{D_t\hg_i}_N\lsim\rg_i^{\xfr12}\tg_q^{-1}\ell_{q,i}^{-N}.\]
Coming to $II$, we estimate it by using \eqref{eq:PSigmaRho}, \eqref{eq:hgi5:c0}, \eqref{eq:hgi6:NormDer}, and \eqref{eq:PCutDer}:
\[\norm{II}_N\leq\xfr{|\T^3|^{\xfr12}\rg_i^{\xfr12}\dg_q^{\xfr12}\lg_q}{c_0^{\xfr12}}\norm{\hg_i}_N\lsim\dg_q^{\xfr12}\lg_q\rg_i^{\xfr12}\leq\rg_i^{\xfr12}\tg_q^{-1}.\]
As for $III$, we use \eqref{eq:hgi5:c0}, \eqref{eq:PSTrDer}, \eqref{eq:hgi6:NormDer}, and the fact that $\yg\leq1$:
\[\norm{III}_N\leq\xfr{|\T^3|^{\xfr12}\dg_q^{\xfr12}\lg_q\sg^{\xfr12}\varsigma^\gg}{2c_0^{\xfr12}}\norm{\hg_i}_N\lsim\rad\sg\dg_q^{\xfr12}\lg_q\lsim\rg_i^{\xfr12}\tg_q^{-1}.\]
Finally, to estimate $IV$, we first note that
\[\norm{h'}_N=\norm{\sum_j\xints{\T^3}{}\22\hg_j\pd_t\hg_j\diff y}_N\lsim\sum_j\xints{\T^3}{}\22\pa{\norm{\hg_j}_N\norm{\pd_t\hg_j}_0+\norm{\hg_j}_0\norm{\pd_t\hg_j}_N}\diff y\leq K\tg_q^{-1}, \sxtag{PDerh}\]
where $K>0$ is a constant. We used \eqref{eq:hgi6:NormDer}. The term $IV$ is thus estimated by combining this bound with \eqref{eq:hgi5:c0}, \eqref{eq:PSigmaRho}, \eqref{eq:hgi6:NormDer}, and the fact $\yg\leq1$:
\[\norm{IV}_N\leq\xfr{K|\T^3|^{\xfr12}\rg_i^{\xfr12}\tg_q^{-1}}{2c_0^{\xfr32}}\norm{\hg_i}_N\lsim\rg_i^{\xfr12}\tg_q^{-1}.\]
By combining the above bounds on $I,II,III,IV$, we conclude that
\[\norm{D_t\rad{\sg_i}}_N\lsim\dg_{q+1}^{\xfr12}\tg_q^{-1}\ell_{q,i}^{-N}. \xtag{PMD$\rad{\sg_i}$}\label{eq:PMDRadSgi}\]
To prove \eqref{eq:Pbik} and \eqref{eq:PMDbik}, we note that
\begin{align*}
\norm{b_{i,k}}_N\lsim{}&\norm{\rad{\sg_i}}_N\norm{a_k(\tilde S_i)}_0+\norm{\rad{\sg_i}}_0\norm{a_k(\tilde S_i)}_N \\
\norm{D_tb_{i,k}}\lsim{}&\norm{D_t\rad{\sg_i}}_N\norm{a_k(\tilde S_i)}_0+\norm{D_t\rad{\sg_i}}_n\norm{a_k(\tilde S_i)}_N+\norm{\rad{\sg_i}}_N\norm{D_t[a_k(\tilde S_i)]}_0+\norm{\rad{\sg_i}}_0\norm{D_t[a_k(\tilde S_i)]}_N.
\end{align*}
The bounds \eqref{eq:Pbik} and \eqref{eq:PMDbik} then readily follow by combining \eqref{eq:PRadSgi}, \eqref{eq:PMDRadSgi}, and the following applications of \eqref{eq:Compos2}:
\begin{align*}
\norm{a_k(\tilde S_i)}_N\lsim{}&\norm{\fD a_k}_0\norm{\fD\tilde S_i}_{N-1}+\norm{\fD a_k}_{N-1}\norm{\fD\tilde S_i}_0^N\lsim\norm{a_k}_N(\norm{\tilde S_i}_N+\norm{\tilde S_i}_1^N), \\
\norm{D_t(a_k(\tilde S_i)\!)}_N\leq{}&\norm{(\fD_2a_k)(\tilde S_i)}_N\norm{D_t\tilde S_i}_0+\norm{(\fD_2a_k)(\tilde S_i)}_0\norm{D_t\tilde S_i}_N \\
{}\lsim{}&(\norm{a_k}_{N+1}\norm{\tilde S_i}_1^N+\norm{a_k}_2\norm{\tilde S_i}_N)\norm{D_t\tilde S_i}_0+\norm{a_k}_1\norm{D_t\tilde S_i}_N,
\end{align*}
where $(\fD_2a_k)_{ij}=\pd_{S_{ij}}a_k$ being the matrix of the first derivatives of $a_k$ w.r.t. the components of its argument.
\nipar To prove \eqref{eq:Pcik}, we note that, by Leibniz rule
\[\norm{c_{i,k}}_N\lsim\sum_{i=0}^N\norm{b_{i,k}}_{i+1}\norm{\fD^T\Fg_i}_{N-i},\]
from which \eqref{eq:Pcik} follows by \eqref{eq:Pbik} and \eqref{eq:PGrad}. Coming finally to \eqref{eq:PMDcik}, we start by noting that
\[D_t\grad(b_{i,k})=\grad D_t(b_{i,k})+[v\per\grad,\grad](b_{i,k}).\]
This means that
\[\norm{D_tc_{i,k}}_N\lsim\sum_{i=0}^N\pa{\norm{\grad D_t(b_{i,k})}_i+\norm{[v\per\grad,\grad](b_{i,k})}_i}\norm{\fD^T\Fg_i}_{N-i}+\sum_{i=0}^N\norm{\fD b_{i,k}}_i\norm{D_t\fD^T\Fg_i}_{N-i}.\]
Since, by the estimates \eqref{eq:PMDGrad}, \eqref{eq:PGrad}, \eqref{eq:Pbik}, \eqref{eq:PMDbik}, and \eqref{eq:Pv0} on the factors here involved, we see that this scales like $\ell_{q,i}^{-N}$, we will only need to prove the case $N=0$. In that case, we obtain that
\[\norm{D_tc_{i,k}}_0\lsim\norm{D_t(b_{i,k})}_1\norm{\fD\Fg_i}_0+\norm{[v\per\grad,\grad](b_{i,k})}_0\norm{\fD\Fg_i}_0+\norm{b_{i,k}}_1\norm{D_t\fD\Fg_i}_0\eqcolon I+II+III.\]
By \eqref{eq:PMDbik} and \eqref{eq:PGrad}, we see that
\[I\lsim\dg_{q+1}^{\xfr12}\tg_q^{-1}\ell_{q,i}^{-1}|k|^{-6}.\]
By \eqref{eq:Pbik}, \eqref{eq:PMDGrad}, and \eqref{eq:PDgTg}, we estimate $III$ as
\[III\lsim\rg_i^{\xfr12}|k|^{-6}\ell_{q,i}^{-1}\per\dg_q^{\xfr12}\lg_q^{1+\ag}\lsim\dg_{q+1}^{\xfr12}\tg_q^{-1}\ell_{q,i}^{-1}|k|^{-6}.\]
We are then left with proving that $II$ also satisfies this bound. To this end, we rewrite the commutator as
\[[v\per\grad,\grad](b_{i,k})=\sum_{j\ell}(v_j\pd_j(\pd_\ell b_{i,k})-\pd_\ell(v_j\pd_jb_{i,k})\!)e_\ell=-\sum_{j\ell}\pd_\ell v_j\pd_jb_{i,k}e_\ell=\grad v\grad b_{i,k}.\]
It then follows from \eqref{eq:PGrad}, \eqref{eq:Pbik}, and \eqref{eq:PDgTg} that $II$ satisfies the same bound as $I$ and $III$, thus proving \eqref{eq:PMDcik}.
\nipar In turn, the estimates on $\tilde v$ in \eqref{eq:PDist0}-\eqref{eq:PDist1} follow from the ones just given precisely as in \cite[Corollary 5.8, pp. 23-24]{BDLSzV}. Indeed, once we note that
\[\norm{\grad(e^{i\lg_{q+1}k\per\Fg_i})}_0\leq\lg_{q+1}|k|\norm{\fD\Fg_i}_0\leq2\lg_{q+1}|k|, \xtag{BDLSzV5.33}\label{eq:BDLSzV5.33}\]
we can deduce from the estimates above that
\begin{align*}
\norm{w_{o,i}}_N\lsim{}&\sum_{i,k}\norm{\fD\Fg_i^{-1}}_N\norm{b_{i,k}}_0\norm{e^{i\lg_{q+1}k\per\Fg_i}}_0+\sum_{i,k}\norm{\fD\Fg_i^{-1}}_0\norm{b_{i,k}}_N\norm{e^{i\lg_{q+1}k\per\Fg_i}}_0 \\
&{}+\sum_{i,k}\norm{\fD\Fg_i^{-1}}_0\norm{b_{i,k}}_0\norm{e^{i\lg_{q+1}k\per\Fg_i}}_N\lsim\dg_{q+1}^{\xfr12}\lg_{q+1}^N.
\end{align*}
The $w_{o,i}$ have pairwise disjoint supports, so the sum over $i$ always consists of a single term, which yields that the desired estimates hold for $w_o$. The estimates on $c_{i,k}\lg_{q+1}^{-1}$ are always better than those on $b_{i,k}$, meaning that any estimate that holds for $w_o$ holds for $w_c$ as well. Thus, \eqref{eq:PDist0} follows directly, and \eqref{eq:PDist1} and \eqref{eq:PFracDist} follow by interpolation.
\nipar Coming to \eqref{eq:PDistH-1}, the fact that $w_c$ satisfies this bound can be easily deduced from \eqref{eq:Pcik}, which tells us that $\lg_{q+1}^{-1}\|c_{i,k}\|_0\lsim\dg_{q+1}^{\xfr12}|k|^{-6}\ell_{q,i}^{-1}\lg_{q+1}^{-1}$. To estimate $w_o$, we use a procedure similar to the one employed in Section \ref{S2S} to prove \eqref{eq:S2SEst1}, replacing \eqref{eq:S2SEstCoeff} with \eqref{eq:Pbik}.
\begin{center}
\tbold{Step 4: Estimates on the new Reynolds term $\pint{\s{E}}\stp1$.}
\end{center}
The aim of this section is to prove $\pint{\s{E}}\stp1$ satisfies \eqref{eq:PErr}, namely
\[\norm{\pint{\s{E}}\stp1}_0\leq\dg_{q+2}\lg_{q+1}^{-6\ag}.\]
Drawing from \BDLSzV, we decompose $\pint{\s{E}}\stp1$ as
\begin{align*}
\pint{\s{E}}\stp1={}&\s{R}\pa[x]{\pd_tw+\opn{div}(v\otimes w+w\otimes v+w\otimes w)-\sum_i\grad\sg_i+\xfrl w-\opn{div}S_\yg} \\
{}={}&\s{R}\pa[x]{\pd_tw+w\per\grad v+v\per\grad w+\opn{div}(w\otimes w)-\sum_i\grad\sg_i+\xfrl w-\opn{div}S_\yg} \\
{}={}&\underbrace{\s{R}(w\per\grad v)}_{\mcl{\text{Nash error}{}\eqcolon\s{E}_N}}+\underbrace{\s{R}(\!(\pd_t+v\per\grad)w)}_{\text{Transport error}{}\eqcolon\s{E}_T}+\underbrace{\s{R}\sq[x]{\opn{div}(w\otimes w-\pint S_\yg)-\sum_i\grad\sg_i}}_{\text{Oscillation error}{}\eqcolon\s{E}_O}+\underbrace{\s{R}(\!\xfrl w)}_{\mcl{\text{Dissipation error}{}\eqcolon\s{E}_D}}.
\end{align*}
We then note that, since the $w_{o,i}$ have disjoint supports and $w_o=\sum_iw_{o,i}$, by \eqref{eq:Psgisg}, we have that
\[\opn{div}(w\otimes w-\pint S_\yg)-\sum_i\grad\sg_i=\opn{div}(w_o\otimes w_c+w_c\otimes w_o+w_c\otimes w_c)+\sum_i\sq{\opn{div}\pa{w_{o,i}\otimes w_{o,i}-\hg_i^2\pint S_\yg-\xfr13\opn{tr}S_i}}. \xtag{PDecEO}\label{eq:PDecEO}\]
We now rewrite the first three terms using the definition of $w$, \eqref{eq:PDecEO} (to rewrite $\s{E}_O$), and the fact that $D_te^{i\lg_{q+1}k\per\Fg_i}=0$ (to rewrite $\s{E}_T$):
\begin{align*}
\s{E}_N={}&\sum_{i,k}\s{R}(\!(\fD\Fg_i^{-1}b_{i,k}A_k+\lg_{q+1}^{-1}c_{i,k})\per\grad ve^{i\lg_{q+1}k\per\Fg_i}) \\
\s{E}_T={}&\sum_{i,k}\s{R}(\!(D_t\fD\Fg_i^{-1}b_{i,k}A_k+\fD\Fg_i^{-1}D_tb_{i,k}A_k+\lg_{q+1}^{-1}D_tc_{i,k})e^{i\lg_{q+1}k\per\Fg_i}) \\
\s{E}_O={}&\s{R}\opn{div}(w_o\otimes w_c+w_c\otimes w_o+w_c\otimes w_c)+\sum_i\s{R}\opn{div}(w_{o,i}\otimes w_{o,i}-S_i) \\
{}={}&\s{R}\opn{div}(w_o\otimes w_c+w_c\otimes w_o+w_c\otimes w_c)+\sum_{i,k}\s{R}(\opn{div}(\sg_i\fD\Fg_i^{-1}C_k(\tilde S_i)\fD^T\Fg_i^{-1})e^{i\lg_{q+1}k\per\Fg_i}),
\end{align*}
where the $C_k$ are as defined in \kcref{thm:lemma:MikFlows}. We now note that the leading order terms are
\begin{align*}
\s{E}_N\stp L\coloneq{}&\sum_{i,k}\s{R}(\!(\fD\Fg_i^{-1}b_{i,k}A_k)\per\grad ve^{i\lg_{q+1}k\per\Fg_i}) \\
\s{E}_T\stp L\coloneq{}&\sum_{i,k}\s{R}(\!(D_t\fD\Fg_i^{-1}b_{i,k}+\fD\Fg_i^{-1}D_tb_{i,k})A_ke^{i\lg_{q+1}k\per\Fg_i}) \\
\s{E}_O\stp L\coloneq{}&\s{R}\opn{div}(w_o\otimes w_c+w_c\otimes w_o)+\sum_{i,k}\s{R}(\opn{div}(\sg_i\fD\Fg_i^{-1}C_k(\tilde S_i)\fD^T\Fg_i^{-1})e^{i\lg_{q+1}k\per\Fg_i})\eqcolon\s{E}_O\stp{L,1}+\s{E}_O\stp{L,2}.
\end{align*}
We start by estimating $\s{E}_O\stp{L,1}$. Since $\s{R}\opn{div}$ is Calderón-Zygmund, we have that
\[\norm{\s{R}\opn{div}(w_o\otimes w_c+w_c\otimes w_o)}_\ag\lsim\norm{w_o}_\ag\norm{w_c}_0+\norm{w_o}_0\norm{w_c}_\ag.\]
From \eqref{eq:Pbik}, \eqref{eq:PGrad}, and \eqref{eq:Pcik}, we can conclude that
\[\norm{w_o}_N\lsim\rg_i^{\xfr12}\lg_{q+1}^N\qquad\norm{w_c}_N\lsim\rg_i^{\xfr12}\lg_{q+1}^{N-1}\ell_{q,i}^{-1}.\]
By interpolation, this lets us conclude that
\[\norm{\s{E}_O\stp{L,1}}_\ag\lsim\rg_i\ell_{q,i}^{-1}\lg_{q+1}^{\ag-1}.\]
To estimate the other leading terms, we start by using \kcref{thm:lemma:SPL} on all three:
\begin{align*}
\norm{\s{E}_N\stp L}_\ag\lsim{}&\sum_{i,k}\pa{\xfr{\norm{\fD\Fg_i^{-1}b_{i,k}\per\grad v}_0}{|k\lg_{q+1}|^{1-\ag}}+\xfr{\norm{\fD\Fg_i^{-1}b_{i,k}\per\grad v}_{N+\ag}+\norm{\fD\Fg_i^{-1}b_{i,k}\per\grad v}_0\norm{\Fg_i}_{N+\ag}}{|k\lg_{q+1}|^{N-\ag}}} \\
\norm{\s{E}_T\stp L}_\ag\lsim{}&\sum_{i,k}\xfr{\norm{D_t\fD\Fg_i^{-1}b_{i,k}+\fD\Fg_i^{-1}D_tb_{i,k}}_0}{|k\lg_{q+1}|^{1-\ag}} \\
&{}+\sum_{i,k}\xfr{\norm{D_t\fD\Fg_i^{-1}b_{i,k}+\fD\Fg_i^{-1}D_tb_{i,k}}_{N+\ag}+\norm{D_t\fD\Fg_i^{-1}b_{i,k}+\fD\Fg_i^{-1}D_tb_{i,k}}_0\norm{\Fg_i}_{N+\ag}}{|k\lg_{q+1}|^{N-\ag}} \\
\norm{\s{E}_O\stp{L,2}}_\ag\lsim{}&\sum_{i,k}\xfr{\norm{\opn{div}(\sg_i\fD\Fg_i^{-1}C_k(\tilde S_i)\fD^T\Fg_i^{-1})}_0}{|k\lg_{q+1}|^{1-\ag}} \\
&{}+\sum_{i,k}\xfr{\norm{\opn{div}(\sg_i\fD\Fg_i^{-1}C_k(\tilde S_i)\fD^T\Fg_i^{-1})}_{N+\ag}+\norm{\opn{div}(\sg_i\fD\Fg_i^{-1}C_k(\tilde S_i)\fD^T\Fg_i^{-1})}_0\norm{\Fg_i}_{N+\ag}}{|k\lg_{q+1}|^{N-\ag}}.
\end{align*}
To estimate $\s{E}_N\stp L$, we combine \eqref{eq:PGrad}, \eqref{eq:Pbik}, and \eqref{eq:Pv0} with a Leibniz inequality:
\begin{multline*}
\sum_{i,k}\pa{\xfr{\norm{\fD\Fg_i^{-1}b_{i,k}\per\grad v}_0}{|k\lg_{q+1}|^{1-\ag}}+\xfr{\norm{\fD\Fg_i^{-1}b_{i,k}\per\grad v}_{N+\ag}+\norm{\fD\Fg_i^{-1}b_{i,k}\per\grad v}_0\norm{\Fg_i}_{N+\ag}}{|k\lg_{q+1}|^{N-\ag}}} \\
{}\lsim\sum_{i,k}\xfr{\norm{\fD\Fg_i^{-1}}_0\norm{b_{i,k}}_0\norm{\grad v}_0}{|k\lg_{q+1}|^{1-\ag}}+\sum_{i,k}\xfr{\norm{\fD\Fg_i^{-1}}_0\norm{b_{i,k}}_0\norm{\grad v}_0\norm{\Fg_i}_{N+\ag}}{|k\lg_{q+1}|^{N-\ag}}+\sum_{i,k}\xfr{\norm{\fD\Fg_i^{-1}}_\ag\norm{b_{i,k}}_\ag\norm{\grad v}_{N+\ag}}{|k\lg_{q+1}|^{N-\ag}} \\
{}+\sum_{i,k}\xfr{\norm{\fD\Fg_i^{-1}}_{N+\ag}\norm{b_{i,k}}_\ag\norm{\grad v}_\ag+\norm{\fD\Fg_i^{-1}}_\ag\norm{b_{i,k}}_{N+\ag}\norm{\grad v}_\ag}{|k\lg_{q+1}|^{N-\ag}} \\
{}\lsim\xfr{\rg_i^{\xfr12}\dg_q^{\xfr12}\lg_q^{1+\ag}}{\lg_{q+1}^{1-\ag}}+\xfr{\ell_{q,i}^{-N-3\ag}\rg_i^{\xfr12}\dg_q^{\xfr12}\lg_q^{1+\ag}}{\lg_{q+1}^{N-\ag}}.
\end{multline*}
The above holds for any $N$. If we choose $N$ to be the $\lbar N$ from Section \ref{Equestimates}, by \eqref{eq:PRho} and \eqref{eq:LambdaEll.2}, we conclude that
\[\norm{\s{R}_N\stp L}_\ag\lsim\xfr{\dg_{q+1}^{\xfr12}\dg_q^{\xfr12}\lg_q}{\lg_{q+1}^{1-2\ag}}.\]
Notice how the leading order term here is the one that does not depend on $N$ thanks to the $\ell_q$ in \eqref{eq:LambdaEll.2}. This is also true of $\s{E}_T\stp L$ and $\s{E}_O\stp{L,2}$.
\nipar $\s{E}_T\stp L$ is estimated in a similar manner, using \eqref{eq:PGrad}, \eqref{eq:PMDGrad}, \eqref{eq:Pbik}, \eqref{eq:PMDbik}, and \eqref{eq:LambdaEll.2}.
\nipar As for $\s{E}_O\stp{L,1}$, we first ensure that adding a derivative, whichever factor it lands on, costs at most $\ell_{q,i}^{-1}$. This ensures that the leading term is the first one, because of that gain of $\ell_q$ mentioned above. We then estimate the leading term.
\nipar By \eqref{eq:PGrad} and \eqref{eq:Pbik}, differentiating $b_{i,k}$ or $\fD\Fg_i^{-1}$ costs $\ell_{q,i}^{-1}$, and by \eqref{eq:SigmaI}, differentiating $\sg_i$ does not cost anything, so we are left with showing that $C_k(\tilde S_i)$ scales like $\ell_{q,i}^{-N}$. Thanks to \eqref{eq:Compos2}, \eqref{eq:MikCk}, and \eqref{eq:PStildeN}, we have that
\[\norm{C_k(\tilde S_i)}_N\lsim\norm{C_k}_1\norm{\fD\tilde S_i}_{N-1}+\norm{\grad C_k}_{N-1}\norm{\tilde S_i}_1^N\lsim|k|^{-6}\ell_{q,i}^{-N}. \sxtag{PCkStildei}\]
We then use \eqref{eq:PSTr}, \eqref{eq:PGrad}, and the above estimate on $C_k(\tilde S_i)$ to estimate the leading term:
\begin{multline*}
\sum_{i,k}\xfr{\norm{\opn{div}(\sg_i\fD\Fg_i^{-1}C_k(\tilde S_i)\fD^T\Fg_i^{-1})}_0}{\abs{k\lg_{q+1}}^{1-\ag}} \\
{}\lsim\sum_{i,k}\xfr{\norm{\sg_i}_1\norm{\fD\Fg_i^{-1}}_0^2\norm{C_k(\tilde S_i)}_0+2\norm{\sg_i}_0\norm{\fD\Fg_i^{-1}}_1\norm{C_k(\tilde S_i)}_0\norm{\fD^T\Fg_i^{-1}}_0+\norm{\sg_i}_0\norm{\fD\Fg_i^{-1}}_0^2\norm{C_k(\tilde S_i)}_1}{|k\lg_{q+1}|^{1-\ag}} \\
{}\lsim\rg_i\ell_{q,i}^{-1}\lg_{q+1}^{-1}.
\end{multline*}
\nipar We thus obtain that
\begin{align*}
\norm{\s{E}_N}_\ag\lsim{}&\xfr{\dg_q^{\xfr12}\dg_{q+1}^{\xfr12}\lg_q}{\lg_{q+1}^{1-2\ag}} \\
\norm{\s{E}_T}_\ag\lsim{}&\xfr{\dg_q^{\xfr12}\dg_{q+1}^{\xfr12}\lg_q}{\lg_{q+1}^{1-5\ag}} \\
\norm{\s{E}_O}_\ag\lsim{}&\xfr{\rg_i}{\ell_{q,i}\lg_{q+1}^{1-\ag}}.
\end{align*}
The relation \eqref{eq:DLRel} easily yields that the above terms satisfy \eqref{eq:PErr} for $a\dgreat1$ sufficiently large, since
\begin{align*}
\xfr{\dg_{q+1}^{\xfr12}\dg_q^{\xfr12}\lg_q}{\lg_{q+1}^{1-3\ag}}\leq{}&\xfr{\dg_{q+1}^{\xfr12}\dg_q^{\xfr12}\lg_q}{\lg_{q+1}^{1-5\ag}}\leq\dg_{q+2}\lg_{q+1}^{-6\ag} \\
\xfr{\rg_i}{\ell_{q,i}\lg_{q+1}^{1-\ag}}\leq{}&\xfr{\Lg\zg_{q+1}^{\xfr{1-\gg}{2}}\zg_q^{\xfr12}\lg_q^{1+\ag}\ell_q^{-\ag}}{\lg_{q+1}^{1-\ag}}\lsim\xfr{\dg_q^{\xfr12}\dg_{q+1}^{\xfr12}\lg_q}{\lg_{q+1}^{1-3\ag-\bg\gg}}\lsim\dg_{q+2}\lg_{q+1}^{-6\ag}. \xtag{DLRelx}\label{eq:DLRelx}
\end{align*}
Coming to $\s{E}_D$, which is not present in \cite{BDLSzV}, we estimate it as follows:
\[\norm{\s{R}(\!\xfrl w)}_0\lsim\norm{\s{R}w}_{2\EXP+\eg}\lsim\norm{\s{R}w}_0^{1-2\EXP-\eg}\norm{\s{R}w}_1^{2\EXP+\eg}. \sxtag{PED1}\]
At this point, we use \kcref{thm:lemma:SPL} to obtain that
\begin{align*}
\norm{\s{R}w_o}_0\lsim{}&\norm{\s{R}w_o}_\ag\lsim\sum_{i,k}\pa{\xfr{\norm{\fD\Fg_i^{-1}b_{i,k}}_0}{|k|^{1-\ag}}+\xfr{\norm{\fD\Fg_i^{-1}b_{i,k}}_{N+\ag}+\norm{\fD\Fg_i^{-1}b_{i,k}}_0\norm{\Fg_i}_{N+\ag}}{|k|^{N-\ag}}} \\
\lsim{}&\dg_{q+1}^{\xfr12}\per\2\sum_{k\neq0}\pa{\xfr{1}{\lg_{q+1}^{1-\ag}|k|^{7-\ag}}+\xfr{\ell_q^{-N-\ag}}{\lg_{q+1}^{N-\ag}|k|^{N-\ag+7}}}\lsim\xfr{\dg_{q+1}^{\xfr12}}{\lg_{q+1}^{1-\ag}}, \sxtag{PED2}
\end{align*}
the last step being due to \eqref{eq:LambdaEll.2}. We also note that
\[\norm{\s{R}w_o}_1=\max_i\norm{\s{R}\pd_iw_o}_0.\]
Proceeding on the $\s{R}\pd_iw_o$ as we did on $\s{R}w_o$ then yields
\begin{align*}
\norm{\s{R}w_o}_1\lsim{}&\max_i\norm{\s{R}\pd_tw_o}_\ag\lsim\sum_{i,k}\pa{\xfr{\norm{\fD\Fg_i^{-1}b_{i,k}}_1}{|k|^{1-\ag}}+\xfr{\norm{\fD\Fg_i^{-1}b_{i,k}}_{N+1+\ag}+\norm{\fD\Fg_i^{-1}b_{i,k}}_1\norm{\Fg_i}_{N+\ag}}{|k|^{N-\ag}}} \\
\lsim{}&\dg_{q+1}^{\xfr12}\lg_{q+1}^\ag. \sxtag{PED3}
\end{align*}
Such estimates analogously also hold for $\s{R}w_c$, and thus for $\s{R}w$. Thus, by \eqref{eq:PED1}-\eqref{eq:PED3}
\[\norm{\s{R}\xfrl w}_0\lsim\dg_{q+1}^{\xfr12}\lg_{q+1}^\ag\lg_{q+1}^{2\EXP+\eg-1}.\]
In particular, for $a\dgreat1$ large enough, \eqref{eq:PErr} is satisfied if
\[2\EXP+\eg-1+\ag-\bg<-2b\bg-6\ag\iff7\ag+\eg<1+\bg-2\EXP-2b\bg. \sxtag{MPSRTermCond}\]
Since $\EXP<\bg$, and $2b\bg<1-\bg$ by \eqref{eq:bBeta}, we have that $1+\bg-2\EXP-2b\bg>0$. Thus, \eqref{eq:MPSRTermCond} above holds for $\ag,\eg$ sufficiently small.
\begin{center}
\tbold{Step 5: Estimates on the new Reynolds term $\s{E}\stp2$.}
\end{center}
Now we turn to $\s{E}\stp2$. Consider the decomposition
\begin{align*}
\abs{\s{E}\stp2}={}&\xfr13\abs{\fint_{\T^3}{}\2|\tilde v|^2-|v|^2-\opn{tr}S_\yg} \\
{}\leq{}&\xfr13\abs{\fint_{\T^3}\2|w_o|^2-\opn{tr}S_\yg}+\xfr13\abs{\fint_{\T^3}\22w\per v}+\xfr13\abs{\fint_{\T^3}\22w_c\per w_o+|w_c|^2}, \xtag{PDecompE2}\label{eq:PDecompE2}
\end{align*}
and proceed as in \cite[Proposition 6.2]{BDLSzV}. In the case of the first term, we will estimate the whole tensor, and therefore the trace. For the other terms, only the trace will be handled.
\nipar Concerning the first term in \eqref{eq:PDecompE2}, thanks to \eqref{eq:Psgisg}, $\sum\int\sg_i=\int\sg_\yg$, so that two cancellations occur:
\begin{align*}
\int w_o\otimes w_o-S_\yg\diff x={\e}&\sum_{i,k\neq0}\int\sg_i\fD\Fg_i^{-1}C_k(\tilde S_i)\fD\Fg_i^{-T}e^{i\lg_{q+1}k\per\Fg_i}\diff x+\int\pa{\sum_i\sg_i-\sg_\yg}\opn{Id}\diff x \\
{}={\e}&\sum_{i,k\neq0}\int Z_{i,k}e^{i\lg_{q+1}k\per\Fg_i}\diff x, \xtag{PStep5.1}\label{eq:PStep5.1}
\end{align*}
where we write $Z_{i,k}\coloneq\sg_i\fD\Fg_i^{-1}C_k(\tilde S_i)\fD\Fg_i^{-T}$. Using \eqref{eq:SPL-Int}, \eqref{eq:MikCk}, and \eqref{eq:LambdaEll.2}, we obtain that
\[\abs{\xints{\T^3}{}\2\sum_{i,k\neq0}Z_{i,k}e^{i\lg_{q+1}k\per\Fg_i}}\lsim\sum_{i,k\neq0}\xfr{\norm{Z_{i,k}}_{\lbar N}+\norm{Z_{i,k}}_0\norm{\Fg_i}_{\lbar N}}{|\lg_{q+1}k|^{\lbar N}}\lsim\sum_{k\neq0}\xfr{\dg_{q+1}\ell_q^{-\lbar N}}{\lg_{q+1}^{\lbar N}|k|^{\lbar N}}\lsim\xfr{\dg_{q+1}\ell_q}{\lg_{q+1}} \xtag{PStep5.2}\label{eq:PStep5.2}.\]
The second inequality above is easily justified by using \eqref{eq:SigmaI}, \eqref{eq:PDistGrad}, and \eqref{eq:PCkStildei} to estimate $Z_{i,k}$ as follows:
\begin{align*}
\norm{Z_{i,k}}_N\lsim{}&\norm{\sg_i}_N\norm{\fD_i^{-1}}_0^2\norm{C_k(\tilde S_i)}_0+2\norm{\sg_i}_0\norm{\fD\Fg_i^{-1}}_N\norm{C_k(\tilde S_i)}\norm{\fD\Fg_i^{-1}}_0+\norm{\sg_i}_0\norm{\fD\Fg_i^{-1}}_0^2\norm{C_k(\tilde S_i)}_N \\
\lsim{}&|k|^{-6}\dg_{q+1}\ell_q^{-N}.
\end{align*}
\nipar To estimate the second term in \eqref{eq:PDecompE2}, observe that
\[w\per v=\sum_{i,k}(\!(\fD\Fg_i)^{-1}b_{i,k}+\lg_{q+1}^{-1}c_{i,k})\per ve^{i\lg_{q+1}k\per\Fg_i},\]
so that, combining \kcref{thm:lemma:SPL}, \eqref{eq:Pbik}, \eqref{eq:Pcik}, \eqref{eq:Pv0}, \eqref{eq:PDistGrad}, and \eqref{eq:LambdaEll.2}, we obtain that
\begin{align*}
\abs{\fint2w\per v\diff x}\lsim{}&\sum_{i,k}\xfr{\norm{(\!(\fD\Fg_i)^{-1}b_{i,k}+\lg_{q+1}^{-1}c_{i,k})\per v}_{\lbar N}+\norm{(\!(\fD\Fg_i)^{-1}b_{i,k}+\lg_{q+1}^{-1}c_{i,k})\per v}_0\norm{\fD\Fg_i}_{\lbar N}}{\lg_{q+1}^{\lbar N}|k|^{\lbar N}} \\
{}\lsim{}&\dg_{q+1}^{\xfr12}\per\dg_q^{\xfr12}\lg_q^{1+\ag}\per\ell_q\lg_{q+1}^{-1}\lsim\dg_{q+1}^{\xfr12}\dg_q^{\xfr12}\lg_q^\ag\lg_{q+1}^{-1}. \xtag{PStep5.3}\label{eq:PStep5.3}
\end{align*}
Concerning the third term in \eqref{eq:PDecompE2}, note that the estimates on $w_c$ are always no coarser than those for $w_o$, so if we estimate $\int w_o\per w_c$ well, the whole term is estimated well. To this end, we observe that
\begin{align*}
\abs{\int w_o\per w_c\diff x}\leq{}&\sum_i\e[2]\sum_{0\neq k\neq l}\abs{\int(\fD\Fg_i)^{-1}b_{i,k}\lg_{q+1}^{-1}c_{i,l-k}e^{i\lg l\per\Fg}\diff x} \\
{}\lsim{}&\sum_i\e[2]\sum_{0\neq k\neq l\neq0}\e[10]\xfr{\norm{(\fD\Fg_i)^{-1}b_{i,k}\lg_{q+1}^{-1}c_{i,l-k}}_{\lbar N}+\norm{(\fD\Fg_i)^{-1}b_{i,k}\lg_{q+1}^{-1}c_{i,l-k}}_0\norm{\fD\Fg_i}_{\lbar N}}{\lg_{q+1}^{\lbar N}|l|^{\lbar N}} \\
&{}+\sum_i\sum_{k\neq0}\abs{\int(\fD\Fg_i)^{-1}b_{i,k}\lg_{q+1}^{-1}c_{i,-k}\diff x}\eqcolon I+II,
\end{align*}
where we used \kcref{thm:lemma:SPL} in the case $l\neq0$, as well as the fact that the $w_{o,i}$ and $w_{c,i}$ have disjoint support so we do not have products of the form $b_{i,k}c_{j,l-k}$ for $i\neq j$. The term $I$ is easily estimated as $\dg_{q+1}\ell_q^{-\lbar N}\lg_{q+1}^{-\lbar N}\lsim\dg_{q+1}\ell_q\lg_{q+1}^{-1}$, so that $I$ satisfies the same estimate as the second term in \eqref{eq:PDecompE2}. As for $II$, \eqref{eq:Pbik} and \eqref{eq:Pcik} easily yield $II\lsim\rg_i\ell_{q,i}^{-1}\lg_{q+1}^{-1}$. Therefore
\[\abs{\int w_o\per w_c\diff x}\lsim\xfr{\dg_{q+1}\ell_q}{\lg_{q+1}}+\xfr{\rg_i\ell_q^{-1}}{\lg_{q+1}}. \xtag{PStep5.4}\label{eq:PStep5.4}\]
Combining \eqref{eq:PDecompE2}-\eqref{eq:PStep5.4} with the fact that $\int|w_c|^2$ also satisfies \eqref{eq:PStep5.4}, we arrive at
\[\abs{\s{E}}\stp2\lsim\xfr{\dg_{q+1}\ell_q}{\lg_{q+1}}+\xfr{\dg_{q+1}^{\xfr12}\dg_q^{\xfr12}\lg_q^{1+\ag}\ell_q}{\lg_{q+1}}+\xfr{\rg_i\ell_{q,i}^{-1}}{\lg_{q+1}}.\]
By \eqref{eq:DLRel}, we thus conclude that, for $a\dgreat1$ sufficiently large, $\s{E}\stp2$ satisfies \eqref{eq:PErr}. Combining with the fact (obtained in the previous step) that $\pint{\s{E}}\stp1$ satisfies \eqref{eq:PErr}, we thus conclude that \eqref{eq:PErr} holds.
\begin{center}
\tbold{Step 6: Estimates on $\pd_t\opn{tr}\s{E}$}
\end{center}
Observe that $\pint{\s{E}}\stp1$ is traceless, whereas $\s{E}\stp2$ is a function of $t$ only. In order to estimate the time derivative of $\s{E}\stp2$, observe that, since $v$ is solenoidal, for every $F=F(x,t)$
\[\Der{t}\xints{\T^3}{}\2F=\xints{\T^3}{}\2D_tF,\]
where $D_t=\pd_t+v\per\grad$. Therefore, using again the decomposition in \eqref{eq:PDecompE2}, we have that
\begin{align*}
\abs{\Der{t}\xints{\T^3}{}\2|\tilde v|^2-|v|^2-\opn{tr}S_\yg}\leq{}&\abs{\xints{\T^3}{}\2\opn{tr}\sq{D_t\pa{\sum_{i,k\neq0}\sg_i\fD\Fg_i^{-1}C_k(\tilde S_i)\fD\Fg_i^{-T}e^{i\lg_{q+1}k\per\Fg_i}}}} \\
&{}+\abs{\xints{\T^3}{}\2D_t(2w_c\per w_o+|w_c|^2)}+\abs{\xints{\T^3}{}\2D_t(2v\per w)}. \xtag{PDecTrDer}\label{eq:PDecTrDer}
\end{align*}
Let us first estimate $\|D_tw_o\|_0$. Recall from \eqref{eq:MDInvGrad} that $D_t(\fD\Fg_i)^{-1}=\fD v(\fD\Fg_i)^{-1}$, which, combined with the fact that $D_te^{i\lg_{q+1}k\per\Fg_i}=0$, yields
\begin{align*}
D_tw_o={}&\sum_{i,k\neq0}D_t\pa{\rad{\sg_i}a_k(\tilde S_i)}\fD\Fg_i^{-1}A_ke^{i\lg_{q+1}k\per\Fg_i} \\
&{}+\sum_{i,k\neq0}\rad{\sg_i}a_k(\tilde S_i)\fD v\fD\Fg_i^{-1}A_ke^{i\lg_{q+1}k\per\Fg_i} \\
{}={}&\sum_{i,k\neq0}\fD\Fg_i^{-1}D_tb_{i,k}e^{i\lg_{q+1}k\per\Fg_i}+\sum_{i,k\neq0}\fD v\fD\Fg_i^{-1}b_{i,k}e^{i\lg_{q+1}k\per\Fg_i}.
\end{align*}
First notice that, by using \eqref{eq:Pv0}, \eqref{eq:PDistGrad}, and \eqref{eq:Pbik}, we obtain that
\[\norm{\fD v\fD\Fg_i^{-1}b_{i,k}}_0\lsim\xfr{\dg_{q+1}^{\xfr12}\dg_q^{\xfr12}\lg_q^{1+\ag}}{|k|^6}.\]
As for the coefficients $\fD\Fg_i^{-1}D_tb_{i,k}$, combining \eqref{eq:PDistGrad} and \eqref{eq:PMDbik} gives
\[\norm{\fD\Fg_i^{-1}D_tb_{i,k}}_0\lsim\tg_q^{-1}\dg_{q+1}^{\xfr12}|k|^{-6}=\xfr{\dg_{q+1}^{\xfr12}\dg_q^{\xfr12}\lg_q\ell_q^{-4\ag}}{|k|^6}.\]
Therefore
\[\norm{D_tw_o}_0\lsim\dg_{q+1}^{\xfr12}\dg_q^{\xfr12}\lg_q\ell_q^{-4\ag}.\]
Observing that
\[D_tw_c=\sum_{i,k}\lg_{q+1}^{-1}D_tc_{i,k}e^{i\lg_{q+1}k\per\Fg_i},\]
which follows from $D_te^{i\lg_{q+1}k\per\Fg_i}=0$ seen above, \eqref{eq:PMDcik} implies
\[\norm{D_tw_c}_0\lsim\dg_{q+1}^{\xfr12}\dg_q^{\xfr12}\lg_q\ell_q^{-1-4\ag}\lg_{q+1}^{-1}.\]
Combining with $\|w_o\|_0+\|w_c\|_0\lsim\dg_{q+1}^{\xfr12}$ and using \eqref{eq:DLRel}-\eqref{eq:LambdaEll}, we obtain that
\begin{align*}
\abs{\xints{\T^3}{}D_t(w_o\otimes w_c+w_c\otimes w_o+w_c\otimes w_c)}\lsim{}&\norm{D_tw_o}_0\norm{w_c}_0+\norm{w_o}_0\norm{D_tw_c}_0+\norm{D_tw_c}_0\norm{w_c}_0 \\
\lsim{}&\dg_{q+1}\dg_q^{\xfr12}\lg_q\ell_q^{-4\ag}=\dg_{q+1}^{\xfr12}(\dg_{q+1}^{\xfr12}\dg_q^{\xfr12}\lg_q)\ell_q^{-4\ag} \\
{}\lsim{}&\dg_{q+1}^{\xfr12}(\dg_{q+2}\lg_{q+1}^{1-10}\ag\lg_{q+1}^{4\ag}=\dg_{q+2}\dg_{q+1}^{\xfr12}\lg_{q+1}^{1-6\ag}.
\end{align*}
The second term of \eqref{eq:PDecTrDer} is thus estimated. We then similarly decompose the first term in \eqref{eq:PDecTrDer} as
\begin{align*}
D_t\sq{\sum_{i,k\neq0}\sg_i\fD\Fg_i^{-1}C_k(\tilde S_i)\grad\Fg_i^{-1}e^{i\lg_{q+1}k\per\Fg_i}}={}&\sum_{i,k\neq0}D_t\sg_i\fD\Fg_i^{-1}C_k(\tilde S_i)\grad\Fg_i^{-1}e^{i\lg_{q+1}k\per\Fg_i} \\
&{}+\sum_{i,k\neq0}\sg_i\fD v\fD\Fg_i^{-1}C_k(\tilde S_i)\grad\Fg_i^{-1}e^{i\lg_{q+1}k\per\Fg_i} \\
&{}+\sum_{i,k\neq0}\sg_i\fD\Fg_i^{-1}D_t[C_k(\tilde S_i)]\grad\Fg_i^{-1}e^{i\lg_{q+1}k\per\Fg_i} \\
&{}+\sum_{i,k\neq0}\sg_i\fD\Fg_i^{-1}C_k(\tilde S_i)\grad\Fg_i^{-1}\grad ve^{i\lg_{q+1}k\per\Fg_i}.
\end{align*}
In order to estimate this, we still need to estimate $D_t[C_k(\tilde S_i)]$ and $D_t\sg_i$. To obtain the former, we first use \eqref{eq:Compos2}:
\begin{align*}
\norm{D_t(C_k(\tilde S_i)\!)}_N\leq{}&\norm{(\fD_2C_k)(\tilde S_i)}_N\norm{D_t\tilde S_i}_0+\norm{(\fD_2C_k)(\tilde S_i)}_0\norm{D_t\tilde S_i}_N \\
{}\lsim{}&(\norm{C_k}_{N+1}\norm{\tilde S_i}_1^N+\norm{C_k}_2\norm{\tilde S_i}_N)\norm{D_t\tilde S_i}_0+\norm{C_k}_1\norm{D_t\tilde S_i}_N,
\end{align*}
We then use \eqref{eq:MikCk}, \eqref{eq:PStildeN}, and \eqref{eq:PMDStilde} to conclude that
\[\norm{D_t(C_k(\tilde S_i)\!)}_N\lsim|k|^{-6}\tg_q^{-1}\ell_{q,i}^{-N}. \sxtag{PMDCkStildei}\]
Coming to $D_t\sg_i$, we claim that
\[\norm{D_t\sg_i}_N\lsim\dg_{q+1}\tg_q^{-1}\ell_{q,i}^{-N}. \sxtag{PMDSgi}\]
To obtain \eqref{eq:PMDSgi}, we set
\begin{align*}
h(t)\coloneq{}&\sum_j\int\hg_j^2(x,t)\diff x \\
D_t\sg_i={}&\xfr{|\T^3|\yg^2\sg}{h}2\hg_iD_t\hg_i+|\T^3|\hg_i^2\pd_t\pa{\xfr{\yg^2\sg}{h}}\eqcolon I+II.
\end{align*}
We first estimate the term $I$. Recalling \eqref{eq:hgi6:NormDer}, \eqref{eq:hgi5:c0}, $\yg\leq1$, and \eqref{eq:PSTr}, we conclude that
\[\norm{I}_N\lsim\xfr{|\T^3|\yg^2\sg}{h}\pa{\norm{\hg_i}_N\norm{D_t\hg_i}_0+\norm{\hg_i}_0\norm{D_t\hg_i}_N}\lsim\dg_{q+1}\tg_q^{-1}\ell_{q,i}^{-N}.\]
As for the second term, we already see that, since the only factor depending on $x$ is $\hg_i^2$ which, by \eqref{eq:hgi6:NormDer}, satisfies $\|\hg_i^2\|_N\lsim1$ for all $N$, the estimates for $II$ will only depend on $N$ via an $a$-independent constant, thus making it sufficient to estimate $\pd_t(\yg^2\sg h^{-1})$ in $\s{C}^0$. To that end, we rewrite it as
\[\pd_t\pa{\xfr{\yg^2\sg}{h}}=\xfr{2\yg\yg'\sg}{h}+\xfr{\yg^2\pd_t\sg}{h}-\xfr{\yg^2\sg h'}{h^2}\eqcolon T_1+T_2+T_3.\]
To estimate $T_1$, we recall \eqref{eq:PSTr}, \eqref{eq:PCutDer}, and \eqref{eq:hgi5:c0}:
\[\norm{T_1}_0\leq\xfr{2\dg_q^{\xfr12}\lg_q\per4\dg_{q+1}}{c_0}\lsim\tg_q^{-1}\dg_{q+1}.\]
Coming to $T_2$, by \eqref{eq:PSTr}, \eqref{eq:PSTrDer}, $\yg\leq1$, and \eqref{eq:hgi5:c0}, we obtain that
\[\norm{T_2}_0\leq\xfr{4C\dg_{q+1}\dg_q^{\xfr12}\lg_q}{c_0}\lsim\dg_{q+1}\tg_q^{-1},\]
where $C$ is the implicit constant in \eqref{eq:PSTrDer}. Finally, to estimate $T_3$, we use \eqref{eq:PSTr}, $\yg\leq1$, and \eqref{eq:PDerh}:
\[\norm{T_3}_0\leq\xfr{4K\dg_{q+1}\tg_q^{-1}}{c_0^2}\lsim\dg_{q+1}\tg_q^{-1},\]
where $K$ is the implicit constant in \eqref{eq:PDerh}. The estimate \eqref{eq:PMDSgi} is thus proved. By \eqref{eq:PMDSgi}, \eqref{eq:PGrad}, \eqref{eq:PCkStildei}, \eqref{eq:SigmaI}, \eqref{eq:Pv0}, and \eqref{eq:PMDCkStildei}, we conclude that
\[\norm{D_t\sq{\sum_{i,k\neq0}\sg_i\fD\Fg_i^{-1}C_k(\tilde S_i)\grad\Fg_i^{-1}e^{i\lg_{q+1}k\per\Fg_i}}}_0\lsim\dg_{q+1}\tg_q^{-1}=\dg_{q+1}\dg_q^{\xfr12}\lg_q\ell_q^{-4\ag}\lsim\dg_{q+1}\dg_{q+1}^{\xfr12}\lg_{q+1}^{1-6\ag},\]
the last step being exactly as done above. Finally, to estimate the term involving $D_t(w\per v)$, we first note that
\[\int D_t(v\per w)=\int D_tv\per w+\int v\per D_tw=-\int(\grad p+\xfrl v+\opn{div}R)\per w+\int v\per D_tw, \xtag{PDecomp\int D_t(v\per w)$}\label{eq:PDecompIntDtVDotW}\]
using \eqref{eq:FNSR} in the last step. To estimate the second term of \eqref{eq:PDecompIntDtVDotW}, we write
\[v\per D_tw=\sum_{i,k\neq0}h_{i,k}e^{i\lg_{q+1}k\per\Fg_i},\]
where
\[h_{i,k}\coloneq v\per D_t[(\fD\Fg_i)^{-1}b_{i,k}+\lg_{q+1}^{-1}c_{i,k}]=v\per[D_t(\fD\Fg_i)^{-1}b_{i,k}+(\fD\Fg_i)^{-1}D_tb_{i,k}+\lg_{q+1}^{-1}D_tc_{i,k}].\]
By \kcref{thm:lemma:HNI}, we obtain that
\begin{align*}
\norm{h_{i,k}}_N\lsim{}&\norm{v}_N\pa{\norm{D_t\fD\Fg_i^{-1}}_0\norm{b_{i,k}}_0+\norm{\fD\Fg_i^{-1}}_0\norm{D_tb_{i,k}}_0+\xfr{1}{\lg_{q+1}}\norm{D_tc_{i,k}}_0} \\
&{}+\norm{v}_0\bigg(\norm{D_t\fD\Fg_i^{-1}}_N\norm{b_{i,k}}_0+\norm{D_t\fD\Fg_i^{-1}}_0\norm{b_{i,k}}_N \\
&{\hsp{1.3cm}}+\norm{\fD\Fg_i}_N\norm{D_tb_{i,ki}}_0+\norm{\fD\Fg_i^{-1}}_0\norm{D_tb_{i,k}}_N+\xfr{1}{\lg_{q+1}}\norm{D_tc_{i,k}}_N\bigg).
\end{align*}
Thus, using \eqref{eq:Pv00}, \eqref{eq:Pv0} (in the form $\|v\|_{N+\ag}\lsim\dg_q^{\xfr12}\lg_q^{1+\ag}\ell_{q,i}^{1-N}\lsim\tg_q^{-1}\ell_q^{-N}$), \eqref{eq:Pbik}-\eqref{eq:Pcik}, \eqref{eq:PMDbik}-\eqref{eq:PMDcik}, and \eqref{eq:PGrad}-\eqref{eq:PMDGrad}, we conclude that
\[\norm{h_{i,k}}_N\lsim\dg_{q+1}^{\xfr12}\tg_q^{-1}\ell_q^{-N}=\dg_{q+1}^{\xfr12}\dg_q^{\xfr12}\lg_q\ell_q^{-4\ag-N}.\]
With \kcref{thm:lemma:SPL}, the above estimate yields that $\int v\per D_tw$ satisfies \eqref{eq:PErrDer}.
\nipar To deal with the first term of \eqref{eq:PDecompIntDtVDotW}, we first note that, since $\opn{div}w=0$, $\int\grad p\per w=0$. Integrating by parts, the term $\int\opn{div}R\per w$ can be estimated as follows:
\[\abs{\int\opn{div}R\per w\diff x}\leq\|R\|_1\|w\|_0\lsim\Lg\varrho^{1+\gg}\ell_q^{-2\ag}\ell_{q,i}^{-1}\per\dg_{q+1}^{\xfr12}\lsim\Lg^{\xfr12}\zg_{q+1}^{\xfr{1+\gg}{2}}\dg_q^{\xfr12}\lg_q^{1+\ag}\ell_q^{-3\ag}\per\dg_{q+1}^{\xfr12},\]
where we used \eqref{eq:PRStrong} and \eqref{eq:PRho}. To conclude that the first term in \eqref{eq:PDecompIntDtVDotW} satisfies \eqref{eq:PErrDer}, we would require
\[\zg_{q+1}^{\xfr{1+\gg}{2}}\zg_q^{\xfr12}\lg_q^{1+\ag}\ell_q^{-3\ag}\lsim\zg_{q+2}\lg_{q+1}^{1-6\ag}. \xtag{PRelParExtra}\label{eq:PRelParExtra}\]
For $\ag,\gg$ sufficiently small, this follows from
\[-b\bg-\bg+1<-2b^2\bg+b\iff1-\bg-2b\bg<b(1-\bg-2b\bg)\iff1-\bg-2b\bg>0,\]
which in turn follows from \eqref{eq:bBeta}.
\begin{center}
\tbold{Step 7: $\s{T}_p$ and its derivative}
\end{center}
The time derivative $\pd_t\s{T}_p$ is readily estimated as
\begin{align*}
\abs{\pd_t\s{T}_p}={}&\abs{\xints{\T^3}{}\2\xfrl*(2v+w)\per\xfrl* w\diff x} \\
{}\leq{}&\abs{\xints{\T^3}{}\22\xfrl*v\per\xfrl*w\diff x}+\xints{\T^3}{}\2\abs{\xfrl*w}^2\diff x \\
{}\lsim{}&(2\|v\|_{\EXP+\eg}\|w\|_0^{1-\EXP-\eg}\|w\|_1^{\EXP+\eg}+\|w\|_0^{2-2\EXP-2\eg}\|w\|_1^{2\EXP+2\eg}).
\end{align*}
By \eqref{eq:PFracv0}, we have that $\|v\|_{\EXP+\eg}\lsim\Lg^{\xfr12}$. As for $w$, we have that $\|w\|_N\lsim\dg_{q+1}^{\xfr12}\lg_{q+1}^N$ (cfr. Step 3 above). Thus, recalling that $\EXP+\eg<\bg$
\[\abs{\pd_t\s{T}_p}\lsim2\Lg^{\xfr12}\per\Lg^{\xfr12}\lg_{q+1}^{\EXP+\eg-\bg}+\Lg\lg_{q+1}^{2\EXP+2\eg-2\bg}\lsim\Lg\lg_{q+1}^{\EXP+\eg-\bg}.\]
Since this is exactly \eqref{eq:PHypoTrace}, the proposition is proved.
\end{qeddim}
\xbegin{oss}[The fractional dissipation term][thm:oss:FracDissTerm]
Note that \eqref{eq:PHypoTrace} is stronger than \eqref{eq:PErrDer}, since
\[\Lg\lg_{q+1}^{\EXP+\eg-\bg}=\Lg^{\xfr12}\dg_{q+1}^{\xfr12}\lg_{q+1}^{\EXP+\eg}\lsim\dg_{q+1}^{\xfr12}\dg_{q+2}\lg_{q+1}^{1-6\ag}.\]
Indeed, this inequality follows from $\EXP+\eg-\bg<1-6\ag-\bg-2b\bg$ which, for $\ag,\eg$ sufficiently small, follows from \eqref{eq:bBeta} and the fact $\EXP<\bg$.
\nipar However, $\s{T}_p$ is only estimated as follows:
\[|\s{T}_p(t)|\lsim t\Lg\lg_{q+1}^{\EXP+\eg-\bg}.\]
To ensure that this satisfies \eqref{eq:PErr} for any $q\geq0$, we would require
\[0<\EXP+\eg-\bg<-2b^2\bg-3b\ag\iff 3b\ag<\bg-\EXP-\eg-2b^2\bg.\]
Seen as the above right-hand side is, in general, negative, we cannot require it. Thus, in general, $\s{T}_p$ only satisfies \eqref{eq:PErr} if the $q$ in the statement is sufficiently large, which is why we separated $\s{T}_p$ from the other Reynolds terms.
\nipar However, for $t\lsim\Lg^{-1}\lg_{q+1}^{\bg-\eg-\EXP}$, we can contrast the growth of $\Lg\lg_{q+1}^{\EXP+\eg-\bg}$ with the smallness of the time, meaning that $\s{T}_p$ only satisfies \eqref{eq:PErr} for a short period of time, or if $q$ is sufficiently large. \end{oss}
\xbegin{oss}[\texorpdfstring{$\s{C}^0$}{C0} estimate on the Reynolds stress][thm:oss:C1Reyn]
The requirement \eqref{eq:PRStrong} is only used to obtain \eqref{eq:PErrDer}, meaning we only need it on $\opn{supp}S$, since $S=0\implies\s{E}=0$.
\xend{oss}

\sect{From strict to adapted subsolutions}\label{S2AProof}
The aim of this section is to prove \kcref{thm:propo:S2A} (p. \pageref{thm:propo:S2A}). The proof closely follows the arguments of \cite[Section 8]{DRSz}. Each stage contains a localized gluing step performed using \kcref{thm:propo:Gluing}, and a perturbation step performed using \kcref{thm:propo:MPS}.
\nipar\begin{qeddim*}[\emph{\kcref{thm:propo:S2A}}]
\begin{center}
\tbold{Step 1: Setting the parameters of the scheme}
\end{center}
Let $\strsub$ be a smooth strict subsolution and let $0<\ad\bg<\bg<\xfr13,\ng>0$. Choose $b>1$ according to \eqref{eq:bBeta}, furthermore let $\tilde\eg>0$ such that:
\[b(1+\tilde\eg)<\xfr{1-\bg}{2\bg}. \xtag{S2AbEpstilde}\label{eq:S2AbEpstilde}\]
Then, let $\tilde\dg,\tilde\gg>0$ be the constants given by \kcref{thm:cor:S2Suse1}, and choose $0<\ag<1$ and $0<\gg<\hat\gg<\tilde\gg$ so that:
\bi
\item The inequalities \eqref{eq:DLRel}, \eqref{eq:LambdaEll} are satisfied by both the pairs $(\ag,\gg)$ and $(\ag,\hat\gg)$;
\item The other conditions in Sections \ref{Gluing} and \ref{MPS}, namely \eqref{eq:GParRel1}-\eqref{eq:GParRel2} and consequently \eqref{eq:PParam1}, \eqref{eq:MPSRTermCond}, and \eqref{eq:PRelParExtra}, are satisfied by both the pairs $(\ag,\gg)$ and $(\ag,\hat\gg)$;
\item Condition \eqref{eq:LambdaEll.2} can hold for both pairs $(\ag,\gg)$ and $(\ag,\hat\gg)$; since $\hat\gg>\gg$, relation \eqref{eq:LambdaEllq2log} reduces this to:
\[(b-1)(1-\bg(b+1)\!)-\hat\gg\bg b^2-2\ag b>0; \xtag{S2ALambdaEll2}\label{eq:S2ALambdaEll2}\]
\item The following conditions holds:
\begin{align*}
\ng>{}&\xfr{1-3\bg+\ag}{2\bg} \xtag{S2AnuBeta}\label{eq:S2AnuBeta} \\
\xfr{\ag}{\bg}<b\hat\gg<\xfr{3\ag}{2\bg},\qquad\qquad0<{}&b\gg<\hat\gg-\xfr{\ag}{\bg},\qquad3\ag>2b\bg\gg. \xtag{S2AalphaBetabGammas}\label{eq:S2AalphaBetabGammas}
\end{align*}
\ei
Having fixed $b,\bg,\ag,\gg,\hat\gg$, we may choose $\lbar N\in\N$ so that \eqref{eq:LambdaEll.2} is also valid. For $a\dgreat1$ sufficiently large (to be determined) we then define $(\lg_q,\dg_q)$ as in \eqref{eq:DefParams}. Thus we are in the setting of Section \ref{Equestimates}.
\begin{center}
\tbold{Step 2: From strict to strong subsolution}
\end{center}
We apply \kcref{thm:cor:S2Suse1} to obtain from $\strsub$ a strong subsolution $(v_0,p_0,R_0)$ with $\dg=\dg_1$ such that the properties from \eqref{eq:U1TrEst} to \eqref{eq:U1TrDer} hold. By \eqref{eq:U1TrEst}-\eqref{eq:U1TrDer}, $(v_0,p_0,R_0)$ satisfies
\begin{align*}
\xfr34\dg_1\leq{}&\rg_0\leq\xfr54\dg_1 \xtag{S2ARho0}\label{eq:S2ARho0} \\
\norm{\pint R_0(t)}_0\leq{}&\Lg\varrho_0^{1+\hat\gg} \xtag{S2AStrong}\label{eq:S2AStrong} \\
\norm{v_0}_{H^{-1}}\leq{}&\lg_0^{-1} \xtag{S2AvH-1}\label{eq:S2AvH-1} \\
\abs{\pd_t\rg_0}\leq{}&\dg_1\dg_0^{\xfr12}\lg_0. \norm{v_0}_{1+\ag}\leq{}&\dg_0^{\xfr12}\lg_0^{1+\ag} \xtag{S2Av0}\label{eq:S2Av0} \\
\abs{\pd_t\rg_0}\leq{}&\dg_1\dg_0^{\xfr12}\lg_0. \xtag{S2ADerRho0}\label{eq:S2ADerRho0}
\end{align*}
\begin{center}
\tbold{Step 3: Inductive construction of $(v_q,p_q,R_q)$}
\end{center}
Starting from $(v_0,p_0,R_0)$, we show how to inductively construct a sequence $\{(v_q,p_q,R_q)\}_{q\in*\N}$ of smooth strong subsolutions with:
\[R_q(x,t)=\rg_q(t)\opn{Id}+\pint R_q(x,t)\]
which satisfy the following properties:
\ben[label={$(\alph*_q)$}]
\item For all $t\in[0,T]$
\[\xints{\T^3}{}\pa[b]{|v_q|^2+\opn{tr}R_q}\diff x=\xints{\T^3}{}\pa[b]{|v_0|^2+\opn{tr}R_0}\diff x;\]
\item For all $t\in[0,T]$
\[\norm{\pint R_q(t)}_0\leq\Lg\varrho_q^{1+\gg};\]
\item If $2^{-j}T<t\leq2^{-j+1}T$ for some $j=1,\dotsc,q$, then
\[\xfr38\dg_{j+1}\leq\rg_q\leq4\dg_j;\]
\item For all $t\leq2^{-q}T$:
\[\norm{\pint R_q(t)}_0\leq\Lg\varrho_q^{1+\hat\gg},\qquad\qquad\xfr34\dg_{q+1}\leq\rg_q\leq\xfr54\dg_{q+1};\]
\item If $2^{-j}T<t\leq2^{-j+1}T$ for some $j=1,\dotsc,q$, then
\begin{align*}
\norm{v_q}_{1+\ag}\leq{}&M\dg_j^{\xfr12}\lg_j^{1+\ag} \\
\abs{\pd_t\rg_q}\lsim{}&\dg_{j+1}\dg_j^{\xfr12}\lg_j, \\
\intertext{whereas if $t\leq2^{-q}T$}
\norm{v_q}_{1+\ag}\leq{}&M\dg_q^{\xfr12}\lg_q^{1+\ag} \\
\abs{\pd_t\rg_q}\lsim{}&\dg_{q+1}\dg_q^{\xfr12}\lg_q.
\end{align*}
\item For all $t\in[0,T]$ and $q\geq1$:
\[\norm{v_q-v_{q-1}}_{H^{-1}}\leq M\dg_q^{\xfr12}(\zg_q^{\xfr\gg2}\ell_{q-1}^\ag+\ell_{q-1}^{-1}\lg_q^{-1})\qquad\norm{v_q-v_{q-1}}_0\leq M\dg_q^{\xfr12}.\]
\item $\|v_q\|_{\EXP+\eg}\Big)\leq M(1+\Lg^{\xfr12}\sum\limits_{i=0}^q\lg_i^{\EXP+\eg-\bg}\Big)$.
\een
Thanks to our choice of parameters in Step 1 above, $(v_0,p_0,R_0)$ satisfies \eqref{eq:S2ARho0}-\eqref{eq:S2ADerRho0}, and thus the inductive assumptions $(a_0)$-$(g_0)$ (the last condition can be deduced from \eqref{eq:U1Fracv}).
\nipar Suppose then $(v_q,p_q,R_q)$ is a smooth strong subsolution which satisfies $(a_q)$-$(g_q)$. The construction of $(v_{q+1},p_{q+1},R_{q+1})$ consists of two steps: first a localized gluing step performed using \kcref{thm:propo:Gluing} to get from $(v_q,p_q,R_q)$ to a smooth strong subsolution $(\lbar v_q,\lbar p_q,\lbar R_q)$, then a localized perturbation step done using \kcref{thm:propo:MPS} to get $(v_{q+1},p_{q+1},R_{q+1})$ from $(\lbar v_q,\lbar p_q,\lbar R_q)$.
\nipar We apply \kcref{thm:propo:Gluing} with
\[[T_1,T_2]=[0,2^{-q}T].\]
Then $T_2-T_1\geq4\tg_q$, if $a\dgreat1$ is sufficiently large. Moreover, by $(d_q)$, $(e_q)$, and $(g_q)$, $(v_q,p_q,R_q)$ fulfils the requirements of \kcref{thm:propo:Gluing} on $[T_1,T_2]$ with parameters $\ag,\hat\gg>0$.
\nipar Then, by \kcref{thm:propo:Gluing}, we obtain a smooth strong subsolution $(\lbar v_q,\lbar p_q,\lbar R_q)$ on $[0,T]$ such that $(\lbar v_q,\lbar p_q,\lbar R_q)$ is equal to $(v_q,p_q,R_q)$ on $[2^{-q}T,T]$, and on $[0,2^{-q}T]$ satisfies
\[\begin{aligned}
\norm{\lbar v_q-v_q}_\ag\lsim{}&\Lg^{\xfr12}\lbar\varrho_q^{\xfr{1+\hat\gg}{2}}\ell_q^\ag \\
\norm{\lbar v_q}_{1+\ag}\lsim{}&\dg_q^{\xfr12}\lg_q^{1+\ag} \\
\norm{\lbar v_q}_{\EXP+\eg}\lsim{}&1+\sum_{i=0}^{q+1}\dg_i^{\xfr12}\lg_i^{\EXP+\eg} \\
\norm{\pint{\lbar R}_q}_0\leq{}&\Lg\lbar\varrho_q^{1+\hat\gg}\ell_q^{-2\ag} \\
\xfr58\dg_{q+1}\leq{}&\lbar\rg_q\leq\xfr32\dg_{q+1} \\
\abs{\pd_t\lbar\rg_q}\lsim{}&\dg_{q+1}\dg_q^{\xfr12}\lg_q.
\end{aligned}
\xtag{S2AvqEst1}\label{eq:S2AvqEst1}\]
Moreover, on $[0,t_{\lbar n}]$ one has that
\[\begin{aligned}
\norm{\lbar v_q}_{N+1+\ag}\lsim{}&\dg_q^{\xfr12}\lg_q^{1+\ag}\ell_q^{-N} \\
\norm{\pint{\lbar R}_q}_{N+\ag}\lsim{}&\Lg\lbar\varrho_q^{1+\hat\gg}\ell_q^{-N-2\ag} \\
\norm{(\pd_t+\lbar v_q\per\grad)\pint{\lbar R}_q}_{N+\ag}\lsim{}&\Lg\lbar\varrho_q^{1+\hat\gg}\ell_q^{-N-6\ag}\dg_q^{\xfr12}\lg_q.
\end{aligned}
\xtag{S2AvqEst2}\label{eq:S2AvqEst2}\]
and
\[\pint{\lbar R}_q\equiv0\qquad t\in\bigcup_{i=0}^{\lbar n}J_i. \xtag{S2ASuppRq}\label{eq:S2ASuppRq}\]
Recalling \kcref{thm:defi:Subintervs+} and \eqref{eq:6.4} observe that
\[\sq{0,\xfr342^{-q}T}\sbs[0,t_{\lbar n}], \xtag{S2AIntervIncl}\label{eq:S2AIntervIncl}\]
provided $a\dgreat1$ is chosen sufficiently large (e.g. so that $\xfr53\tg_q<\xfr142^{-q}T$). Then, choose a cut-off function $\yg_q\in\Cinf_C([0,\xfr342^{-q}T];[0,1])$ such that
\[\yg_q(t)=\case{
1 & t\leq2^{-(q+1)}T \\
0 & t>\xfr342^{-q}T
} \xtag{S2APsiq}\label{eq:S2APsiq}\]
and such that $|\yg'_q(t)|\lsim2^q$. By choosing $a\dgreat1$ sufficiently large, we may assume that
\[\abs{\yg'_q(t)}\leq\xfr12\dg_q^{\xfr12}\lg_q \xtag{S2APsiq2}\label{eq:S2APsiq2}\]
for all $q$. Then, set
\[S_\yg\coloneq\yg_q^2(\lbar R_q-\dg_{q+2}\opn{Id})=\yg_q^2S.\]
Using \eqref{eq:S2APsiq2}, \eqref{eq:S2AalphaBetabGammas}, \eqref{eq:S2AvqEst1}-\eqref{eq:S2AIntervIncl}, and the easy observation that $\lbar\rg_q\lsim\lbar\rg_q-\dg_{q+2}$, we see that $S_\yg$ and $(\lbar v_q,\lbar p_q,\lbar R_q)$ satisfy the assumptions of \kcref{thm:propo:MPS} on the interval $[0,t_{\lbar n}]$ with parameters $\ag,\hat\gg>0$. We have that
\[\sg_\yg=\yg_q^2(\lbar\rg_q-\dg_{q+2})=\yg_q^2\sg.\]
Recalling \kcref{thm:oss:C1Reyn}, since $\opn{supp}S_\yg\sbse[t_{\ubar n},t_{\lbar n}]$ where \eqref{eq:GGStrong2} holds, we can apply \kcref{thm:propo:MPS}, thus obtaining a new subsolution $(v_{q+1},p_{q+1},\lbar R_q-S_\yg-\s{E}_{q+1})$ with
\begin{align*}
\norm{v_{q+1}-\lbar v_q}_0+\ell_q\lg_{q+1}\norm{v_{q+1}-\lbar v_q}_{H^{-1}}\hsp{.5cm}& \\
{}+\lg_{q+1}^{-1-\ag}\norm{v_{q+1}-\lbar v_q}_{1+\ag}\hsp{.5cm}& \\
{}+\lg_{q+1}^{-\EXP-\eg}\norm{v_{q+1}-\lbar v_q}_{\EXP+\eg}\leq{}&M\dg_{q+1}^{\xfr12} \\
\xints{\T^3}{}\2|v_{q+1}|^2-\opn{tr}S-\opn{tr}\s{E}_{q+1}={}&\e\xints{\T^3}{}\2|\lbar v_q|^2\qquad t\in[0,T],
\end{align*}
and such that the estimates \eqref{eq:PErr} and \eqref{eq:PErrDer} hold for $\s{E}_{q+1}$.
Let
\[R_{q+1}\coloneq\lbar R_q-S_\yg-\s{E}_{q+1}.\]
We claim that $(v_{q+1},p_{q+1},R_{q+1})$ is a smooth strong subsolution satisfying $(a_{q+1})$-$(g_{q+1})$. Notice that $(a_{q+1})$ is satisfied by construction. Since $(v_{q+1},p_{q+1},R_{q+1})=(v_q,p_q,R_q)$ for $t\geq2^{-q}T$, we may restrict $t$  to $[0,2^{-q}T]$ in the following, so that in particular \eqref{eq:S2AvqEst1} holds.
\nipar Let us now prove $(b_{q+1})$. On the one hand
\begin{align*}
\norm{\pint R_{q+1}}_0={}&\norm{(1-\yg_q^2)\pint{\lbar R}_q-\pint{\s{E}}_{q+1}}_0 \\
{}\leq{}&(1-\yg_q^2)\Lg\lbar\varrho_q^{1+\hat\gg}\ell_q^{-2\ag}+\dg_{q+2}\lg_{q+1}^{-3\ag}\one_{\{\yg_q>0\}}, \xtag{S2ARq+1est}\label{eq:S2ARq+1est}
\end{align*}
on the other hand
\begin{align*}
\rg_{q+1}={}&(1-\yg_q^2)\Lg\lbar\varrho_q+\yg_q^2\dg_{q+2}+\xfr13\opn{tr}\s{E}_{q+1} \\
{}\geq{}&(1-\yg_q^2)\lbar\rg_q+\yg_q^2\dg_{q+2}-\dg_{q+2}\lg_{q+1}^{-3\ag}\one_{\{\yg_q>0\}}. \xtag{S2Argq+1est}\label{eq:S2Argq+1est}
\end{align*}
The proof of $(b_{q+1})$ thus reduces to assessing whether there exists a suitable $\gg$ such that
\[(1-\yg_q^2)\Lg^{-\hat\gg}\lbar\rg_q^{1+\hat\gg}\ell_q^{-2\ag}+\dg_{q+2}\lg_{q+1}^{-3\ag}\one_{\{\yg_q>0\}}\leq\Lg^{-\gg}[(1-\yg_q^2)\lbar\rg_q+\yg_q^2\dg_{q+2}-\dg_{q+2}\lg_{q+1}^{-3\ag}]^{1+\gg}. \xtag{S2AStrongBar}\label{eq:S2AStrongBar}\]
To this end set
\begin{align*}
F(s)\coloneq{}&(1-s)\Lg\lbar\varrho_q^{1+\hat\gg}\ell_q^{-2\ag}+\dg_{q+2}\lg_{q+1}^{-3\ag} \\
G(s)\coloneq{}&(1-s)\lbar\rg_q+s\dg_{q+2}-\dg_{q+2}\lg_{q+1}^{-3\ag} \\
H(s)\coloneq{}&\Lg^{-\gg}G^{1+\gg}(s)-F(s),
\end{align*}
and observe that \eqref{eq:S2AStrongBar} is equivalent to $H(\yg_q^2)\geq0$ if $\yg_q>0$, and follows from this inequality also in case $\yg_q=0$. In particular, \eqref{eq:S2AStrongBar} follows from:
\ben[label=(\roman*)]
\item $H(0)\geq0$ and $H(1)\geq0$;
\item $H'(0)\leq0$ and $H'(1)\leq0$.
\item $H''(s)\geq0$.
\een
We note next that, since $2b\bg\hat\gg<3\ag$
\[\dg_{q+2}\lg_{q+1}^{-3\ag}\lsim\Lg\lbar\varrho_q^{1+\hat\gg},\]
so that we have the estimates
\[F(0)\lsim\Lg\lbar\varrho_q^{1+\hat\gg}\ell_q^{-2\ag},\qquad\qquad G(0)\gtrsim\lbar\rg_q.\]
It is also clear that $G(s)\leq\lbar\rg_q$.
\nipar It is then easy to check that the requirement $H(0)\geq0$, i.e. $F(0)\leq\Lg^{-\gg}G^{1+\gg}(0)$, amounts to $\Lg\lbar\varrho_q^{1+\hat\gg}\ell_q^{-2\ag}\lsim\Lg\lbar\varrho_q^{1+\gg}$, i.e. $\lbar\varrho_q^{\hat\gg-\gg}\ell_q^{-2\ag}\lsim1$. Hence, since $\ell_q^{-1}\leq\lg_{q+1}$ by \eqref{eq:LambdaEll} and $\lbar\varrho_q\gsim\zg_{q+1}$ by $(d_{q+1})$, $H(0)\geq0$ follows from
\[\hat\gg-\xfr{\ag}{\bg}>\gg, \xtag{S2Amu}\label{eq:S2Amu}\]
provided $a\dgreat1$ is sufficiently large to absorb geometric constants. The relation \eqref{eq:S2Amu} follows from \eqref{eq:S2AalphaBetabGammas} since $b>1$.
\nipar The next requirement, $H(1)\geq0$, i.e. $\lg_{q+1}^{-3\ag}\lsim\zg_{q+2}^\gg(1-\lg_{q+1}^{-3\ag})^{1+\gg}$, requires $3\ag>2b\bg\gg$ as found in \eqref{eq:S2AalphaBetabGammas}, since $1-\lg_{q+1}^{-3\ag}\geq\xfr12$ for $a$ sufficiently large.
\nipar The following condition, $H'(0)\leq0$, can be rewritten as
\[-\Lg\lbar\varrho_q^{1+\hat\gg}\ell_q^{-2\ag}\geq(1+\gg)(\lbar\varrho_q-\zg_{q+2}\lg_{q+1}^{-3\ag})^\gg(\dg_{q+2}-\lbar\rg_q)\sse\Lg\lbar\varrho_q^{1+\hat\gg}\ell_q^{-2\ag}\lsim(\lbar\varrho_q-\zg_{q+2}\lg_{q+1}^{-3\ag})^\gg(\lbar\rg_q-\dg_{q+2}).\]
Noting that $\lbar\rg_q\gsim\dg_{q+1}\dgreat\dg_{q+2}$ by $(d_{q+1})$, and therefore $\lbar\rg_q-\dg_{q+2}\geq\xfr12\lbar\rg_q$ for $a$ sufficiently large, the above reduces to
\[\lbar\varrho_q^{\hat\gg-\gg}\ell_q^{-2\ag}\lsim1\slse\lg_{q+1}^{2\ag-2\bg\hat\gg+2\bg\gg}=\zg_{q+1}^{\hat\gg-\gg}\lg_{q+1}^{2\ag}\lsim1,\]
which follows from condition \eqref{eq:S2Amu} deduced above.
\nipar We then need the condition $H'(1)\leq0$, which can be rewritten as
\[-\Lg\lbar\varrho_q^{1+\hat\gg}\ell_q^{-2\ag}\geq(1+\gg)\zg_{q+2}^\gg(1-\lg_{q+1}^{-3\ag})^\gg(\dg_{q+2}-\lbar\rg_q),\]
which similarly follows from \eqref{eq:S2Amu}.
\nipar The last condition, $H''\geq0$, follows from the fact that $F''\equiv0$ and $G''\equiv0$, and thus $H''=\Lg^{-\gg}(1+\gg)\gg G^{\gg-1}G'^2$ is positive.
\nipar Thus, our choice of $\ag,\gg,\hat\gg$ in \eqref{eq:S2AalphaBetabGammas} guarantees that \eqref{eq:S2AStrongBar} holds, which yields $(b_{q+1})$.
\nipar Consider now $(c_{q+1})$, where we only need to consider the case $j=q+1$, i.e. the estimate on $[2^{-q-1}T,2^{-q}T]$. Using \eqref{eq:S2Argq+1est}, the fact that $\lbar\rg_q\geq\dg_{q+2}$ for $a$ large enough, and \eqref{eq:S2AvqEst1}, we see that
\[\dg_{q+2}(1-\lg_{q+1}^{-3\ag})\leq\rg_{q+1}(t)\leq\lbar\rg_q(t)+\dg_{q+2}\lg_{q+1}^{-3\ag}\leq\xfr32\dg_{q+1}+\dg_{q+2}\lg_{q+1}^{-3\ag}.\]
Therefore $(c_{q+1})$ holds, provided $a\dgreat1$ is sufficiently large.
\nipar Similarly, concerning $(d_{q+1})$, observe that for $t\leq2^{-(q+1)}T$ we have that $\yg_q(t)=1$, so that
\[\dg_{q+2}(1-\lg_{q+1}^{-3\ag})\leq\rg_{q+1}\leq\dg_{q+2}(1+\lg_{q+1}^{-3\ag}).\]
Moreover, using \eqref{eq:S2ARq+1est} and the fact that $\yg_q=1$ for $t\leq2^{-(q+1)}T$
\[\norm{\pint R_{q+1}}_0\leq\dg_{q+2}\lg_{q+1}^{-3\ag}\leq\Lg\pa{\xfr34\zg_{q+2}}^{1+\hat\gg},\]
where we used the fact that $2b\bg\hat\gg<3\ag$ and chose $a\dgreat1$ sufficiently large. Therefore $(d_{q+1})$, i.e. $\|\pint R_{q+1}\|_0\lsim\Lg\varrho_{q+1}^{1+\hat\gg}$ and $\xfr34\dg_{q+2}\leq\rg_{q+1}\leq\xfr54\dg_{q+2}$, holds.
\nipar Concerning $(e_{q+1})$, it is once more enough to restrict to $t\leq2^{-q}T$, i.e. the case $j=q+1$. From \eqref{eq:S2AvqEst1} and \eqref{eq:PDist1} we deduce that
\begin{align*}
\norm{v_{q+1}}_{1+\ag}\leq{}&\norm{v_{q+1}-\lbar v_q}_{1+\ag}+\norm{\lbar v_q}_{1+\ag} \\
{}\leq{}&\xfr M2\dg_{q+1}^{\xfr12}\lg_{q+1}^{1+\ag}+C\dg_q^{\xfr12}\lg_q^{1+\ag} \\
{}\leq{}&M\dg_{q+1}^{\xfr12}\lg_{q+1}^{1+\ag},
\end{align*}
where $C$ is the implicit constant in \eqref{eq:S2AvqEst1}, which can be absorbed by choosing $a\dgreat1$ sufficiently large. The estimate on $|\pd_t\rg_{q+1}|$ similarly follows from the trace estimate of \eqref{eq:S2AvqEst1} and \eqref{eq:PErrDer}. $(e_{q+1})$ is thus proved.
\nipar $(f_{q+1})$ follows from \eqref{eq:S2AvqEst1}, \eqref{eq:PDist0}, and \eqref{eq:PDistH-1}.
\nipar Finally, $(g_{q+1})$ easily follows from \eqref{eq:PFracDist} and \eqref{eq:S2AvqEst1}.
\begin{center}
\tbold{Step 4: Convergence to an adapted subsolution}
\end{center}
We have thus obtained a sequence $(v_q,p_q,R_q)$ satisfying $(a_q)$-$(g_q)$.
\nipar From $(f_q)$ it follows that $(v_q,p_q)$ is a Cauchy sequence in $\s{C}^0$. Indeed, it is clear for $\{v_q\}$, and concerning $\{p_q\}$ we may use \eqref{eq:FNSR} to write
\[\Dg(p_{q+1}-p_q)=-\opn{div}\opn{div}(\pint R_{q+1}-\pint R_q+(v_{q+1}-v_q)\otimes v_q+v_{q+1}\otimes(v_{q+1}-v_q)\!),\]
and apply Schauder estimates (\kcref{thm:lemma:Schauder}). Similarly, $\{R_q\}$ also converges in $\s{C}^0$. Indeed, from the definition and using \eqref{eq:GGStrong}, \eqref{eq:GUStrong}, \eqref{eq:PSStrong}, \eqref{eq:PErr}, and $(b_q)$, we have that
\begin{align*}
\norm{R_{q+1}-R_q}_0={}&\norm{\lbar R_q-R_q-S_\yg-\s{E}_{q+1}}_0 \\
{}\leq{}&\norm{\lbar R_q}_0+\norm{R_q}_0+\norm{S_\yg}_0+\norm{\s{E}_{q+1}}_0 \\
{}\lsim{}&\dg_{q+1}.
\end{align*}
For all $t>0$ there exists $q(t)\in\N$ such that
\[(v_q,p_q,R_q)(\per,t)=(v_{q(t)},p_{q(t)},R_{q(t)})(\per,t)\qquad\VA q\geq q(t),\]
thus $(v_q,p_q,R_q)$ converges uniformly to a strong subsolution $\adsub$ satisfying
\[\norm{\ad R}_0\leq\Lg\ad\varrho^{1+\gg},\]
and, using \eqref{eq:U1EC} and $(a_q)$
\[\xints{\T^3}{}\pa[b]{|\ad v|^2+\opn{tr}\ad R}\diff x=\e\xints{\T^3}{}\pa[b]{|\str v|^2+\opn{tr}\str R}\diff x\qquad\VA t\in[0,T].\]
Furthermore, using \eqref{eq:U1vH-1} and $(f_q)$
\begin{align*}
\norm{\ad v-\str v}_{H^{-1}}\leq{}&\norm{v_0-\str v}_{H^{-1}}+\norm{v_0-\ad v}_{H^{-1}} \\
{}\lsim{}&\dg_1\lg_0^{-1}+\sum_{q=0}^\8\norm{v_{q+1}-v_q}_{H^{-1}} \\
{}\lsim{}&\dg^{\xfr12}\zg_q^{\xfr\gg2}\ell_q^\ag,
\end{align*}
leading to \eqref{eq:S2Av} for $a$ sufficiently large. Using $(f_q)$ and the fact that $\str v,\strng v$ are smooth and thus bounded in $\s{C}^0$, \eqref{eq:S2AvC0} is proved similarly:
\begin{align*}
\norm{\ad v-\str v}_{\s{C}^0}\leq{}&\norm{v_0-\str v}_{\s{C}^0}+\norm{v_0-\ad v}_{\s{C}^0} \\
{}\lsim{}&1+\sum_{q=0}^\8\norm{v_{q+1}-v_q}_0 \\
{}\lsim{}&1+\dg_1^{\xfr12}.
\end{align*}
Concerning the initial datum, from $(e_q)$ and $(f_q)$ we obtain by interpolation that $\hat v(\per,0)\in\s{C}^{\ad\bg}$, and from $(d_q)$ we obtain that $\ad R(\per,0)=0$.
\nipar Finally, we verify conditions \eqref{eq:ASv}, and \eqref{eq:ASTr} for being a $\s{C}^{\ad\bg}$-adapted subsolution. Let $t>0$. Then there exists $q\in\N$ such that $t\in[2^{-q}T,2^{-q+1}T]$. By $(c_q)$ and $(e_q)$
\begin{align*}
\xfr38\dg_{q+1}\leq{}&\ad\rg\leq4\dg_q \\
\norm{\ad v}_{1+\ag}\leq{}&M\dg_q^{\xfr12}\lg_q^{1+\ag}.
\end{align*}
Therefore $\ad\rg^{-1}\geq\xfr14\dg_q^{-1}$, and hence, using \eqref{eq:DefParams} and \eqref{eq:S2AnuBeta}, we deduce that
\begin{align*}
\norm{\ad v}_{1+\ag}\leq{}&\Lg^{\xfr12}\ad\varrho^{-(1+\ng)},
\end{align*}
for $a\dgreat1$ sufficiently large. Similarly, using $(e_q)$ and \eqref{eq:S2AnuBeta}, we deduce that
\[\abs{\pd_t\ad\rg}\lsim\dg_{q+1}\dg_q^{\xfr12}\lg_q=\Lg^{\xfr32}\lg_q^{1-\bg}\lg_{q+1}^{-2\bg}\sim\Lg^{\xfr32}\lg_q^{1-\bg-2b\bg}=\Lg^{\xfr32}\zg_q^{-\xfr{1}{2\bg}(1-\bg-2b\bg)}\leq\Lg^{\xfr32}\zg_q^{1-\xfr{1-\bg}{2\bg}}\lsim\Lg^{\xfr32}\ad\varrho^{-\ng}.\]
Finally, a word about the term
\[\ad{\s{T}}\coloneq\sum(\s{T}_g\stp q+\s{T}_d\stp q),\]
where $\s{T}_g\stp q$ and $\s{T}_d\stp q$ are the extra trace terms from the $q$th gluing and perturbation steps. We have that $|\pd_t\s{T}_g\stp q|+|\pd_t\s{T}_d\stp q|\lsim\Lg\lg_q^{\EXP+\eg-\bg}$, thus proving \eqref{eq:S2ADerT}. However, adding $\hat{\s{T}}$ into $\ad R$ could compromise the adaptedness of $\adsub$ by rendering \eqref{eq:ASv}-\eqref{eq:ASTr} invalid, which is why we keep it separated and deal with it in the final argument. The estimate \eqref{eq:S2ADerT} implies that
\[\abs{\ad{\s{T}}(t)}\lsim\sum t\Lg\lg_q^{\EXP+\eg-\bg}.\]
To be able to make it as small as we desire, we must contrast the $a$-growth of the $q=0$ and $q=1$ terms of this sum. This is easily achieved by requiring $t\leq\Lg^{-1}\lg_0^{\bg-\EXP-\eg-\ig}$ for $\ig$ arbitrarily small. In any case, calling $t_s$ the maximal time where $\ad{\s{T}}$ can be estimated with small quantities, we have that
\[\lim_{a\to\8}t_s=0,\]
since we meed $t_s\dg_0^{\nicefrac12}\lg_0^{\EXP+\eg}$ to be small.
\end{qeddim*}

\sect{From adapted subsolutions to solutions}\label{A2SProof}
The aim of this section is to prove \kcref{thm:propo:A2S} (p. \pageref{thm:propo:A2S}). The proof closely follows the arguments of \cite[Section 9]{DRSz}. We now start from an adapted subsolution and, by a convex integration scheme, build a sequence of strong subsolutions which converge to a solution of the fractional Navier-Stokes equation. As in \kcref{thm:propo:S2A}$\e[3]$, the convex integration scheme needs the localized gluing and perturbation arguments of \kcref{thm:propo:Gluing} (in the form of \kcref{thm:oss:MultiGluing}) and \kcref{thm:propo:MPS}. However, the choice of the cut-off functions will be, as in \DSz, dictated by the shape of the trace part of the Reynolds stress, and not fixed a priori as in \kcref{thm:propo:S2A}. Before we start the proof, a remark needs to be made about starting the chain of \kcref{thm:propo:Gluing} and \kcref{thm:propo:MPS} with worse estimates.
\xbegin{oss}[Worse starting estimate][thm:oss:WorseStartEst]
In \kcref{thm:propo:Gluing}, if we replace \eqref{eq:GUStrong} with
\[\norm{\pint R}_0\leq\Lg\varrho_q^{1+\gg}\ell_q^{-\xfr2b\ag},\]
as we will need to do below, the estimates \eqref{eq:GMR}, \eqref{eq:GGDistN}, \eqref{eq:GGDist2}, \eqref{eq:GGDist}, \eqref{eq:GGv}, \eqref{eq:GGv2}, \eqref{eq:GGStrong}, \eqref{eq:GGStrong2}, \eqref{eq:GGTrDist}, \eqref{eq:GGPreMatDer}, and \eqref{eq:GGPreTrDist} will be worsened by a factor $\ell_q^{-\xfr2b\ag}$. In fact, we can gain a factor $\ell_q^\ag$ in \eqref{eq:GGStrong} and \eqref{eq:GGStrong2}, and a factor $\ell_q^\ag\lg_q^\ag$ in \eqref{eq:GGTrDist} and \eqref{eq:GGPreMatDer}. To keep the inductive estimates on the velocity gap $\|v_{q+1}-v_q\|_0$ and $\|v_{q+1}-v_q\|_{H^{-1}}$, the velocity $\|v_{q+1}\|_0$, and the derivative of the trace $|\pd_t\rg_q|$, we will need
\[\begin{aligned}
\Lg^{\xfr12}\varrho_q^{\xfr{1+\gg}{2}}\ell_q^{(1-\xfr2b)\ag}\lsim{}&\dg_{q+1}^{\xfr12} \\
\dg_q^{\xfr12}\lg_q\ell_q^{-\xfr2b\ag}\lsim{}&\dg_{q+1}^{\xfr12}\lg_{q+1} \\
\varrho_{q,i}^\gg\ell_q^{-(2+\xfr2b)\ag}\lsim{}&1,
\end{aligned}\xtag{WorseEstFix}\label{eq:WorseEstFix}\]
all of which can easily be deduced by assuming $2\ag<\bg\gg$ and $\ag<\xfr29$. The former assumption also yields \eqref{eq:GGDistA2S}, which will allow us to bound the $H^{-1}$ norm of $v_{q+1}-v_q$ sufficiently tightly. If we then start the perturbation step of \kcref{thm:propo:MPS} from estimates that we can obtain from the modified output estimates mentioned above, we can get the same output estimates from \kcref{thm:propo:MPS}.
\end{oss}
\begin{qeddim*}[\emph{\kcref{thm:propo:A2S}}]
\begin{center}
\tbold{Step 1: Setting the parameters in the scheme}
\end{center}
Let $\adsub$ be a $\s{C}^{\ad\bg}$-adapted subsolution on $[0,T]$, with $\Wg=\Lg$, satisfying the ``strong'' condition $|\pint{\ad R\,\,\,}\!\!\!|\leq\Lg\ad\varrho^{1+\gg}$ for some $\gg>0$ and \eqref{eq:ASv} and \eqref{eq:ASTr} for some $\ag,\ng>0$ as in \kcref{thm:defi:Adapted} of adapted subsolution, with
\[\xfr{1-\ad\bg}{2\ad\bg}<1+\ng<\xfr{1-\sol\bg}{2\sol\bg}.\]
Fix $b>0$ so that
\[b^2(1+\ng)<\xfr{1-\sol\bg}{2\sol\bg},\qquad2\sol\bg(b^2-1)<1. \xtag{A2SbNu}\label{eq:A2SbNu}\]
Observe that both the strongness condition \eqref{eq:Strong} and the adaptedness conditions \eqref{eq:ASv}-\eqref{eq:ASTr} remain valid for any $\hat\gg<\gg$ and $\ag'\leq\ag$ (cfr. \kcref{thm:oss:StrengthParams}). Then, we may assume that $\ag,\hat\gg>0$ are sufficiently small, so that $\adsub$ satisfies \eqref{eq:Strong} for some $\hat\gg>0$ and \eqref{eq:ASv}-\eqref{eq:ASTr} for some $\ag,\ng>0$, and furthermore choose $\gg$ so that
\[2\ag<\sol\bg\hat\gg<\sol\bg\gg<3\ag\qquad b\sol\bg\hat\gg<3\ag. \xtag{A2SGammas}\label{eq:A2SGammas}\]
For the reasons discussed in \kcref{thm:oss:WorseStartEst} above, and for another technical reason we will see below, we require
\[2\ag<\ad\bg\gg<3\ag. \xtag{A2SAgGg}\label{eq:A2SAgGg}\]
Finally, having fixed $b,\ad\bg,\sol\bg,\ag,\gg,\hat\gg$, we may choose $\lbar N\in\N$ so that \eqref{eq:LambdaEll.2} holds. For $a\dgreat1$ sufficiently large (to be determined) we then define $(\lg_q,\dg_q)$ as in \eqref{eq:DefParams} (using $\sol\bg$). Thus, we are in the setting of Section \ref{Equestimates}.
\begin{center}
\tbold{Step 2: Conditions on $(v_0,p_0,R_0)$ and the inductive construction of $(v_q,p_q,R_q)$}
\end{center}
Differently from \cite[Section 9]{DRSz}, we can take $(v_0,p_0,R_0)=\adsub$, since we are assuming $\ad\rg\leq\xfr54\dg_1=\xfr54\dg$, which is $a$-independent. We do have some estimates to verify for $(v_0,p_0,R_0)$, namely that, wherever $\rg_0\geq\dg_{q+2}$
\[\begin{aligned}
\norm{v_0}_{1+\ag}\leq{}&\dg_q^{\xfr12}\lg_q^{1+\ag} \\
\abs{\pd_t\rg_0}\leq{}&\rg_0\dg_q^{\xfr12}\lg_q.
\end{aligned}
\xtag{A2Sv0}\label{eq:A2Sv0}\]
Indeed, where $\rg_0\geq\dg_{q+2}$, \eqref{eq:A2SbNu} easily yields
\begin{align*}
\Lg^{\xfr12}\varrho_0^{-(1+\ng)}\lsim{}&\Lg^{\xfr12}\lg_q^{2\sol\bg b^2(1+\ng)}\leq\dg_q^{\xfr12}\lg_q \\
\Lg^{\xfr32}\varrho_0^{-\ng}\lsim{}&\Lg^{\xfr32}\lg_q^{2\sol\bg b^2\ng}\leq\dg_{q+2}\dg_q^{\xfr12}\lg_q,
\end{align*}
provided $a\dgreat1$ is sufficiently large. These two relations, combined with \eqref{eq:ASv} and \eqref{eq:ASTr}, yield \eqref{eq:A2Sv0}.
\nipar Start from $(v_0,p_0,R_0)$, we will inductively construct a sequence $(v_q,p_q,R_q)$ of smooth strong subsolutions for $q=1,2,\dotsc$, with
\[R_q(x,t)=\rg_q(t)\opn{Id}+\pint R_q(x,t),\]
satisfying the following properties:
\ben[label=$(\Alph*_q)$]
\item For all $t\in[0,T]$
\[\xints{\T^3}{}\pa[b]{|v_q|^2+\opn{tr}R_q}\diff x=\xints{\T^3}{}\pa[b]{|v_0|^2+\opn{tr}R_0}\diff x; \xtag{A}\label{eq:A2Sa}\]
\item For all $t\in[0,T]$
\[\rg_q\leq\xfr52\dg_{q+1}; \xtag{B}\label{eq:A2Sb}\]
\item For all $t\in[0,T]$
\[\norm{\pint R_q}_0\leq\case{
\Lg\varrho_q^{1+\hat\gg}\ell_q^{-\xfr2b\ag} & \rg_q\geq2\dg_{q+2} \\[.5em]
\Lg\varrho_q^{1+\hat\gg} & \xfr32\dg_{q+2}\leq\rg_q\leq2\dg_{q+2} \\[.5em]
\Lg\varrho_q^{1+\gg} & \rg_q\leq\xfr32\dg_{q+2}
}; \xtag{C}\label{eq:A2Sc}\]
\item If $\rg_q\geq\dg_{j+2}$ for some $j\geq q$, then
\begin{align*}
\norm{v_q}_{1+\ag}\leq{}&M\dg_j^{\xfr12}\lg_j^{1+\ag} \xtag{D.1}\label{eq:A2Sd1} \\
\abs{\pd_t\rg_q}\leq{}&\rg_q\dg_j^{\xfr12}\lg_j; \xtag{D.2}\label{eq:A2Sd2}
\end{align*}
\item For all $t\in[0,T]$ and $q\geq1$
\[\norm{v_q-v_{q-1}}_{H^{-1}}\lsim(\zg_q^{\xfr\gg2}\ell_q^{\xfr\ag2}+\dg_q^{\xfr12}\lg_q^{-1})\qquad\norm{v_q-v_{q-1}}_0\lsim\dg_q^{\xfr12}. \xtag{E}\label{eq:A2Se}\]
\item $\|v_q\|_{\EXP+\eg}\leq M\pa{1+\Lg^{\xfr12}\sum\limits_{i=0}^q\lg_i^{\EXP+\eg-\sol\bg}}$.
\een
Thanks to our choice of parameters in Step 1 above, $(v_0,p_0,R_0)$ satisfies \eqref{eq:A2Sv0}, and therefore our inductive assumptions $(A_0)$-$(F_0)$.
\nipar Suppose now $(v_q,p_q,R_q)$ satisfies $(A_q)$-$(F_q)$ above. Let
\[J_q\coloneq\br{t\in[0,T]:\rg_q(t)>\xfr32\dg_{q+2}},\qquad K_q\coloneq\br{t\in[0,T]:\rg_q(t)\geq2\dg_{q+2}}.\]
Being (relatively) open in $[0,T]$, $J_q$ is a disjoint, possibly countable, union of (relatively) open intervals $(T_1\stp i,T_2\stp i)$. Let
\[\s{I}_q\coloneq\br{i:(T_1\stp i,T_2\stp i)\cap K_q\neq\0},\]
and let $t_0\in(T_1\stp i,T_2\stp i)\cap K_q$ for some $i\in\s{I}_q$. Since $K_q$ is compact, we may assume that the open interval $(T_1\stp i,t_0)$ is contained in $J_q\ssm K_q$. Using \eqref{eq:A2Sd2}, we then have that
\begin{align*}
\xfr32\dg_{q+2}={}&\rg_q(T_1\stp i)\geq\rg_q(t_0)-|T_1\stp i-t_0|\sup_{J_q}\abs{\pd_t\rg_q} \\
{}\geq{}&2\dg_{q+2}-2\dg_{q+2}\dg_q^{\xfr12}\lg_q|T_1\stp i-t_0|,
\end{align*}
hence
\[\abs{T_1-t_0}\geq\xfr14(\dg_q^{\xfr12}\lg_q)^{-1}=\xfr{\ell_q^{-4\ag}}{4}\tg_q>4\tg_q, \xtag{A2S.$T_1\stp i-t_0$}\label{eq:A2ST1it0}\]
provided $a\dgreat1$ is chosen sufficiently large. A similar estimate holds for $T_2\stp i$. Therefore $T_2\stp i-T_1\stp i>4\tg_q$ for any $i\in\s{I}_q$, so that $\s{I}_q$ is a finite index set.
\nipar Next, we apply \kcref{thm:propo:Gluing} (in the form of \kcref{thm:oss:MultiGluing}), keeping \kcref{thm:oss:WorseStartEst} in mind, to $(v_q,p_q,R_q)$ on this disjoint union of intervals $\bigcup_{\s{I}_q}J_{q,i}$. Since $\rg_q>\xfr32\dg_{q+2}$, from $(A_q)$-$(F_q)$ and \eqref{eq:A2SbNu}-\eqref{eq:A2SGammas} we see that the assumptions of \kcref{thm:propo:Gluing} on $(v_q,p_q,R_q)$ hold with parameter $\hat\gg$. Then we obtain $(\lbar v_q,\lbar p_q,\lbar R_q)$ such that, on $J_q$
\begin{align*}
\norm{\lbar v_q(t)-v_q(t)}_\ag\lsim{}&\dg_{q+1}^{\xfr12}\ell_q^{(\xfr12+\xfr2b-\xfr2b\one_{K_q})}\lsim\dg_{q+1}^{\xfr12}\ell_q^{\xfr\ag2} \tag{\text{From }\eqref{eq:GGDistA2S}} \\
\norm{\lbar v_q}_{1+\ag}\lsim{}&\dg_q^{\xfr12}\lg_q^{1+\ag}\ell_q^{-\xfr2b\ag\one_{K_q}}\lsim\dg_{q+1}^{\xfr12}\lg_{q+1}^{1+\ag} \tag{\text{From }\eqref{eq:GGv}} \\
\norm{\pint{\lbar R}_q}_0\leq{}&\lbar\rg_q^{1+\hat\gg}\ell_q^{-2\ag+(1-\xfr2b)\ag\one_{K_q}} \tag{\text{From }\eqref{eq:GGStrong}} \\
\xfr78\rg_q\leq{}&\Lg\lbar\varrho_q\leq\xfr98\rg_q \tag{\text{From }\eqref{eq:GGTr}} \\
\abs{\pd_t\lbar\rg_q}\lsim{}&\lbar\rg_q\dg_q^{\xfr12}\lg_q. \tag{\text{From }\eqref{eq:GGDer}}
\end{align*}
Moreover, recalling \eqref{eq:6.4}, for any $i\in\s{I}_q$ we have the following additional estimates valid for $t\in[T_1\stp i+2\tg_q,T_2\stp i-2\tg_q]\cap J_q$:
\[\begin{aligned}
\norm{\lbar v_q}_{N+1+\ag}\lsim{}&\dg_q^{\xfr12}\lg_q^{1+\ag}\ell_q^{-N-\xfr2b\ag\one_{K_q}} \\
\norm{\pint{\lbar R}_q}_{N+\ag}\lsim{}&\lbar\Lg\lbar\varrho_q^{1+\hat\gg}\ell_q^{-N-2\ag+(1-\xfr2b)\ag\one_{K_q}} \\
\norm{(\pd_t+\lbar v_q\per\grad)\pint{\lbar R}_q}_{N+\ag}\lsim{}&\Lg\lbar\varrho_q^{1+\hat\gg}\ell_q^{-N-6\ag+(1-\xfr2b)\ag\one_{K_q}}\dg_q^{\xfr12}\lg_q,
\end{aligned} \xtag{A2SEst}\label{eq:A2SEst}\]
and
\[\opn{supp}\pint{\lbar R}_q\sbs\T^3\x\bigcup_iI_i, \xtag{A2Ssupp}\label{eq:A2Ssupp}\]
where $\{I_i\}_i$ are the intervals defined in \eqref{eq:DeftiIi}.
\nipar Let us choose a cut-off function $\yg_q\in\Cinf_c(J_q;[0,1])$ such that
\begin{align*}
\opn{supp}\yg_q\sbs{}&\bigcup_{i\in*\s{I}_q}\pa{T_1\stp i+2\tg_q,T_2\stp i-2\tg_q} \xtag{A2SsuppPsiq}\label{eq:A2SsuppPsiq} \\
K_q\sbs{}&\br{\yg_q=1} \xtag{A2SKqPsiq}\label{eq:A2SKqPsiq} \\
\abs{\yg'_q}\lsim{}&\dg_q^{\xfr12}\lg_q. \xtag{A2SDerPsiq}\label{eq:A2SDerPsiq}
\end{align*}
Such a choice is made possible by \eqref{eq:A2ST1it0}. We then want to apply \kcref{thm:propo:MPS} (using \kcref{thm:oss:WorseStartEst} above where $\rg_q\geq2\dg_{q+2}$) to $(\lbar v_q,\lbar p_q,\lbar R_q)$ with
\[S_\yg\coloneq\yg_q^2(\lbar R_q-\dg_{q+2}\opn{Id}),\]
hence $\sg_\yg=\yg_q^2(\lbar\rg_q-\dg_{q+2})$. Using \eqref{eq:A2SDerPsiq}, \eqref{eq:A2SGammas}, \eqref{eq:A2SbNu}, \eqref{eq:A2SEst}, \eqref{eq:A2Ssupp}, \eqref{eq:GGTr}-\eqref{eq:GGDer}, and $(A_q)$-$(F_q)$, we see that $S$ and $(\lbar v_q,\lbar p_q,\lbar R_q)$ satisfy the required assumptions on the interval $[T_1\stp i+2\tg_q,T_2\stp i-2\tg_q]$ with parameters $\ag,\hat\gg>0$. In particular, \eqref{eq:PRStrong} (or its worsened form discussed in \kcref{thm:oss:WorseStartEst}) follows from \eqref{eq:A2SEst}, since we only need it on $\opn{supp}\yg_q\sbse[T_i\stp1+2\tg_q,T_i\stp2-2\tg_q]$.
\nipar\kcref{thm:propo:MPS} gives then a new subsolution $(v_{q+1},p_{q+1},\lbar R_q-S_\yg-\s{E}_{q+1})$ with
\begin{align*}
\norm{v_{q+1}-\lbar v_q}_0+\norm{v_{q+1}-\lbar v_q}_{H^{-1}}\lg_{q+1}\hsp{.5cm}& \\
{}+\lg_{q+1}^{-1-\ag}\norm{v_{q+1}-\lbar v_q}_{1+\ag}\hsp{.5cm}& \\
{}+\lg_{q+1}^{-\EXP-\eg}\norm{v_{q+1}-\lbar v_q}_{\EXP+\eg}\leq{}&M\dg_{q+1}^{\xfr12} \tag{\text{From \eqref{eq:PDist0} and \eqref{eq:PDist1}}} \\
\xints{\T^3}{}\2\abs{v_{q+1}}^2-S_\yg-\s{E}_{q+1}={}&\xints{\T^3}{}\2\abs{\lbar v_q}^2\qquad t\in[0,T]. \tag{\text{From }\eqref{eq:PPIC}}
\end{align*}
and such that $\s{E}_{q+1}$ satisfies \eqref{eq:PErr}-\eqref{eq:PErrDer}. Let
\[R_{q+1}\coloneq\lbar R_q-S_\yg-\s{E}_{q+1},\]
\renewcommand{\theoss}{\thesection.\arabic{oss}}%
We claim that $(v_{q+1},p_{q+1},R_{q+1})$ is a smooth strong subsolution satisfying $(A_{q+1})$-$(F_{q+1})$. Notice that $(A_{q+1})$ is satisfied by construction. By definition of $S_\yg$, one has that
\begin{align*}
\rg_{q+1}={}&\lbar\rg_q(1-\yg_q^2)+\yg_q^2\dg_{q+2}-\xfr13\opn{tr}\s{E}_{q+1} \\
\pint R_{q+1}={}&\pint{\lbar R}_q(1-\yg_q^2)-\pint{\s{E}}_{q+1}.
\end{align*}
For $t\in K_q$, condition $(B_{q+1})$ follows easily from \eqref{eq:PErr} and the fact that $K_q\sbs\{\yg_q=1\}$. For $t\nin J_q$, we have that
\[\rg_{q+1}=\rg_q\leq\xfr32\dg_{q+2}<\xfr52\dg_{q+2}.\]
For $t\in J_q\ssm K_q$, we have that $\lbar\rg_q\leq\xfr98\rg_q\leq\xfr98\per2\dg_{q+2}=\xfr94\dg_{q+2}$, which means
\[\rg_{q+1}\leq\xfr94\dg_{q+2}\pa{1-\xfr49\yg_q^2+\xfr59\lg_{q+1}^{-6\ag}}\leq\xfr54\dg_{q+2}\pa{\xfr95+\lg_1^{-6\ag}},\]
and if $\lg_1^{-6\ag}\leq\xfr15$, which is a matter of choosing $a$ large enough, we have $(B_{q+1})$.
\nipar Note that, by the construction of $\rg_{q+1}$, we have that $J_q\sbse K_{q+1}$, since on the whole of $J_q$ we have that $\rg_{q+1}\sim\dg_{q+2}\dgreat\dg_{q+3}$. This is the reason why we required $\hat\gg<\gg$ and used the larger $\gg$ outside of $J_q$ in $(C_q)$: to make sure $(C_{q+1})$ was automatically verified outside $J_q$. This is in stark contrast to what happened in Section \ref{S2AProof}, where the perturbation regions $P_q\coloneq[0,2^{-q}T]$ satisfied the opposite inclusion $P_{q+1}\sbse P_q$, and where we consequently required $\hat\gg>\gg$ to ensure the weaker ``strong condition'' $(b_{q+1})$ in $P_q\ssm P_{q+1}$, while the stronger $(d_{q+1})$ only held in $P_{q+1}$, where $\yg_q=1$.
\nipar By the above paragraph, in verifying conditions $(C_{q+1})$-$(D_{q+1})$, it suffices to restrict to the case when $\rg_{q+1}\geq2\dg_{q+3}$ and $j=q+1$, respectively.
\nipar The argument showing $(C_{q+1})$ for $t\in K_{q+1}$ is similar to the proof of $(b_{q+1})$ in Step 3 of \kcref{thm:propo:S2A} above. On the one hand
\begin{align*}
\norm{\pint R_{q+1}}_0={}&\norm{(1-\yg_q^2)\pint{\lbar R}_q-\pint{\s{E}}_{q+1}}_0 \\
{}\leq{}&(1-\yg_q^2)\Lg\lbar\varrho_q^{1+\hat\gg}\ell_q^{-2\ag+(1-\xfr2b)\ag\one_{\yg_q=1}}+\dg_{q+2}\lg_{q+1}^{-6\ag},
\end{align*}
on the other hand
\begin{align*}
\rg_{q+1}={}&(1-\yg_q^2)\Lg\lbar\varrho_q+\yg_q^2\dg_{q+2}+\xfr13\opn{tr}\s{E}_{q+1} \\
{}\geq{}&(1-\yg_q^2)\lbar\rg_q+\yg_q^2\dg_{q+2}-\dg_{q+2}\lg_{q+1}^{-6\ag}.
\end{align*}
So where $\yg_q=1$ we have the condition since, for $a$ large enough, we can guarantee
\[\dg_{q+2}\lg_{q+1}^{-6\ag}\lsim\dg_{q+2}\zg_{q+2}^{\hat\gg}(1-\lg_{q+1}^{-6\ag})^{1+\hat\gg}\iff\lg_{q+1}^{-6\ag}\lsim\zg_{q+2}^{\hat\gg}(1-\lg_{q+1}^{-6\ag})^{1+\hat\gg}, \xtag{A2Sblah}\label{eq:A2Sblah}\]
since $6\ag>2b\sol\bg\hat\gg$ is required in \eqref{eq:A2SGammas}. If $\yg_q\neq1$, however, we need
\[(1-\yg_q^2)\Lg^{-\hat\gg}\lbar\rg_q^{1+\tilde\gg}\ell_q^{-2\ag}+\dg_{q+2}\lg_{q+1}^{-6\ag}\one_{\{\yg_q>0\}}\leq\Lg^{-\hat\gg}[(1-\yg_q^2)\lbar\rg_q+\yg_q^2\dg_{q+2}-\dg_{q+2}\lg_{q+1}^{-6\ag}]^{1+\hat\gg}\ell_{q+1}^{-\xfr2b\ag}. \xtag{A2SStrongBar}\label{eq:A2SStrongBar}\]
To this end set
\begin{align*}
F(s)\coloneq{}&(1-s)\Lg\lbar\varrho_q^{1+\tilde\gg}\ell_q^{-2\ag}+\dg_{q+2}\lg_{q+1}^{-6\ag} \\
G(s)\coloneq{}&(1-s)\lbar\rg_q+s\dg_{q+2}-\dg_{q+2}\lg_{q+1}^{-6\ag}=\lbar\rg_q+s(\dg_{q+2}-\lbar\rg_q)-\dg_{q+2}\lg_{q+1}^{-6\ag} \\
H(s)\coloneq{}&\Lg^{-\tilde\gg}G^{1+\tilde\gg}(s)\ell_{q+1}^{-\xfr2b\ag}-F(s),
\end{align*}
and, just like in \kcref{thm:propo:S2A}, deduce that $H(\yg_q^2)\geq0$ by proving that:
\ben[label=(\roman*)]
\item $H(0)\geq0$ and $H(1)\geq0$;
\item $H'(0)\geq0$ and $H'(1)\geq0$.
\item $H''(s)\geq0$.
\een
To this end, we first obtain the estimates
\[\dg_{q+2}\lg_{q+1}^{-6\ag}\lsim\Lg\lbar\varrho_q^{1+\hat\gg},\qquad\quad F(0)\lsim\Lg\lbar\varrho_q^{1+\hat\gg}\ell_q^{-2\ag},\qquad G(0)\gtrsim\lbar\rg_q,\qquad\quad G(s)\leq\lbar\rg_q.\]
The first one follows from \eqref{eq:A2Sblah}, \eqref{eq:GGTr}, and the fact we are working for $t\in J_q$.. The second one follows from the first one. The fourth one is obvious, since $\lbar\rg_q\geq\xfr78\xfr32\dg_{q+2}\xfr{21}{16}\dg_{q+2}>\dg_{q+2}$. For the third one, we reduce it to $\dg_{q+2}\lg_{q+1}^{-6\ag}\lsim\lbar\rg_q$, and then it follows from the first estimate, since $\Lg\lbar\varrho_q^{1+\hat\gg}\leq\lbar\rg_q$. We then prove (i)-(v) as follows.
\bi
\item It is easy to check that the two parts of (i) amount to
\[\Lg\lbar\varrho_q^{1+\hat\gg}\ell_q^{-2\ag}\lsim\Lg\lbar\varrho_q^{1+\hat\gg}\ell_{q+1}^{-\xfr2b\ag},\qquad\dg_{q+2}\lg_{q+1}^{-6\ag}\leq\Lg^{-\hat\gg}[\dg_{q+2}(1-\lg_{q+1}^{-6\ag})]^{1+\hat\gg}\ell_{q+1}^{-\xfr2b\ag};\]
the first one follows from $\ell_q\sim\ell_{q+1}^{\nicefrac1b}$; the second one follows from \eqref{eq:A2SGammas} and the following relations, which hold for $a$ sufficiently large:
\[1-\lg_{q+1}^{-6\ag}\geq\xfr12\iff{}\lg_{q+1}^{-6\ag}\leq\xfr12,\qquad\qquad\lg_{q+1}^{-6\ag}\leq\zg_{q+2}^{\hat\gg}\ell_{q+1}^{-\xfr2b\ag}2^{-1-\hat\gg};\]
\item The requirements (ii) can be rewritten as
\[\lbar\rg_q^{1+\hat\gg}\ell_q^{-2\ag}\geq(1+\hat\gg)(\lbar\rg_q-\dg_{q+2})\ell_{q+1}^{-\xfr2b\ag}\max\{[\lbar\rg_q-\dg_{q+2}\lg_{q+1}^{-6\ag}]^{\hat\gg},\dg_{q+2}^{\hat\gg}(1-\lg_{q+1}^{-6\ag})^{\hat\gg}\},\]
which easily follows for sufficiently small $\hat\gg$ and sufficiently large $a$, since $\ell_q^{-2\ag}\sim\ell_{q+1}^{-\xfr2b\ag}$;
\item Note that $G''=0$ because $G$ is linear in $s$, and the same is true of $F''$, meaning that (iii) is simply
\begin{align*}
0\leq{}&\Lg^{-\hat\gg}\hat\gg(1+\hat\gg)G^{\hat\gg-1}(s)G'^2(s)\ell_{q+1}^{-\xfr2b\ag} \\
{}={}&\Lg^{-\hat\gg}\hat\gg(1+\hat\gg)[(1-s)\lbar\rg_q+\dg_{q+2}(1-\lg_{q+1}^{-6\ag})]^{\hat\gg-1}(\lbar\rg_q-\dg_{q+2})^2\ell_{q+1}^{-\xfr2b\ag},
\end{align*}
which is obvious, since all those factors are positive.
\ei
We have thus obtained $(C_{q+1})$.
\nipar The velocity estimate in $(D_{q+1})$ for $j=q+1$ follows from \eqref{eq:PDist1} and \eqref{eq:GGv2}. The trace estimate in $(D_{q+1})$ follows from \eqref{eq:PErrDer} and \eqref{eq:GGDer}. Finally, $(E_{q+1})$ follows precisely as $(f_{q+1})$ in the proof of \kcref{thm:propo:S2A} in Section \ref{S2AProof} above, and $(F_{q+1})$ is obtained just like $(g_{q+1})$. Keep in mind \kcref{thm:oss:WorseStartEst} above for all of these.
\nipar Thus, the inductive step is proved.
\nipar Finally, the convergence of $\{v_q\}$ to a solution of the hypodissipative Navier-Stokes equations as in the statement of \kcref{thm:propo:A2S} (i.e. the one we are proving) follows easily from the sequence of estimates in $(A_q)$-$(F_q)$, analogously to Step 4 of \kcref{thm:propo:S2A} proved in Section \ref{S2AProof} above.
\nipar The Navier-Stokes term $\sol{\s{T}}$ will be handled in the same way as $\ad{\s{T}}$ was dealt with in \kcref{thm:propo:S2A}, giving us once more that the maximal time $t_s$ of ``smallness'' of $\sol{\s{T}}$ must satisfy
\[\lim_{a\to\8}t_s=0.\]
\begin{center}
\tbold{Step 3: From one to infinitely many}
\end{center}
Looking at the details of the scheme, we realize that replacing $\adsub$ with
\[(v_0,p_0,R_0)\coloneq(\ad v,\ad p,\ad R+{}^e/_3),\]
as described in the statement of \kcref{thm:propo:A2S}, clearly retains condition \eqref{eq:AIV}, since the initial datum is not changed. It does not necessarily preserve conditions \eqref{eq:ASv}, and \eqref{eq:ASTr}. Those, however, are only needed to obtain the conditions $(D_0)$. If we then show that the conditions $(A_0)-(F_0)$ (and thus also $(D_0)$) are maintained with such a perturbation, we need not worry about losing \eqref{eq:ASv} and \eqref{eq:ASTr}.
\nipar Conditions about the velocity are clearly preserved, and $(A_0)$ and $(E_0)$ are vacuous, so all we need is
\begin{align*}
\rg_0'\leq{}&\xfr52\dg\tag{$B_0$} \\
\norm{\pint R_0'}_0\leq{}&\case{
\Lg\varrho_0'^{1+\hat\gg}\ell_q^{-\xfr2b\ag} & \rg_0'\geq2\dg_2 \\
\Lg\varrho_0'^{1+\hat\gg} & \xfr32\dg_2\leq\rg_0'\leq2\dg_2 \\
\Lg\varrho_0'^{1+\gg} & \rg_0'\leq\xfr32\dg_2
} \tag{$C_0$} \\
\abs{\pd_t\rg_0'}\leq{}&\rg_0'\dg_0^{\xfr12}\lg_0. \tag{$D_0$.2, i.e. \eqref{eq:A2Sd2}}
\end{align*}
Concerning $(B_0)$, the proposition assumes $\ad\rg\leq{}^{_{^{\scriptstyle5}}}\2/{\2\,}_2\dg$, so that the condition is preserved by requiring \eqref{eq:CondE1}. Since $\ad\rg(0)=0$, $e$ has the possibility to vary in a neighborhood of $t=0$ without becoming negative.
\nipar $(D_0.3)$ boils down to the following condition on $e$:
\[\case{
\abs{\pd_te}\leq e\rad{\dg_0}\lg_0 & e>0 \\
\abs{\pd_te}+\abs e\rad{\dg_0}\lg_0\leq\ad\rg\rad{\dg_0}\lg_0-\abs{\pd_t\ad\rg} & \text{otherwise}
}.\]
To keep things simple, we require \eqref{eq:CondE2} and \eqref{eq:CondE3}.
\nipar Coming to $(C_0)$, we first assume $e>0$, which immediately yields, by the properties of $\adsub$, that
\[\norm{\pint R_0}_0\leq\case{
\Lg(\ad\varrho+\Lg^{-1}e)^{1+\hat\gg}\ell_q^{-\xfr2b\ag} & \ad\rg\geq2\dg_2 \\
\Lg(\ad\varrho+\Lg^{-1}e)^{1+\hat\gg} & \xfr32\dg_2\leq\ad\rg\leq2\dg_2 \\
\Lg(\ad\varrho+\Lg^{-1}e)^{1+\gg} & \ad\rg\leq\xfr32\dg_2
}.\]
Our goal is to obtain that
\[\norm{\pint R_0}_0\leq\case{
\Lg(\ad\varrho+\Lg^{-1}e)^{1+\hat\gg}\ell_q^{-\xfr2b\ag} & \ad\rg+e\geq2\dg_2 \\
\Lg(\ad\varrho+\Lg^{-1}e)^{1+\hat\gg} & \xfr32\dg_2\leq\ad\rg+e\leq2\dg_2 \\
\Lg(\ad\varrho+\Lg^{-1}e)^{1+\gg} & \ad\rg+e\leq\xfr32\dg_2
}.\]
\bi
\item We first note that $\ad\rg+e\leq{}^3/_2\dg_2\implies\ad\rg\leq{}^3/_2\dg_2$, so in this case we have the desired estimate for $\|\pint R_0\|_0$;
\item If $\ad\rg\leq{}^3/_2\dg_2$ but ${}^3/_2\dg_2\leq\ad\rg+e\leq2\dg_2$, since $\hat\gg<\gg$ and $\|\pint R_0\|_0\leq\Lg(\ad\varrho+\Lg^{-1}e)^{1+\gg}$, we have the desired estimate;
\item If $\ad\rg\leq{}^3/_2\dg_2$ but $\ad\rg+e\geq2\dg_2$, the desired estimate is even looser than in the previous item;
\item If ${}^3/_2\dg_2\leq\ad\rg\leq2\dg_2$, then either $\ad\rg+e$ also satisfies this bound, in which case we have the desired estimate, or $\ad\rg+e\geq2\dg_2$, in which case the desired estimate is looser;
\item Finally, if $\ad\rg\geq2\dg_2$, then so is $\ad\rg+e$, meaning again we have the desired estimate.
\ei
Thus, we need no additional conditions to obtain $(C_0)$. Summing up, the conditions we must impose on $e$ are precisely \eqref{eq:CondE1}-\eqref{eq:CondE3}. The proof is complete.
\end{qeddim*}

\appendix
\sect{Proof of \texorpdfstring{\kcref{thm:lemma:ExistStrictSubsol}}{thm:lemma:ExistStrictSubsol}}\label{ExistStrictSubsol}
\begin{qeddim*}[\kcref{thm:lemma:ExistStrictSubsol}]
Fix $\rg\in\Cinf_c(B_1(0)\!)$ a standard mollification kernel in space, and define:
\[\rg_\eg(x)\coloneq\eg^{-3}\rg(x\eg^{-1}).\]
To ensure the regularity of the initial datum, we consider the smoothed datum
\[w_0\coloneq w\ast\rg_{\hg_0},\]
where
\[\hg_0\coloneq\max\br{\hg:\|w\ast\rg_\hg-w\|_{L^2}\leq\xfr\dg3\wedge\xints{\T^3}{}\relax[\abs{w_0}^2-\abs{w}^2](x)\diff x\leq\xfr23\dg}. \xtag{ESSDefHg0}\label{eq:ESSDefHg0}\]
By \kcref{thm:teor:ExistLeraySol}, there exists a solution $(\tilde v,\tilde p)$ with initial datum $w_0$, where $\tilde p$ can be recovered uniquely once we impose $\int\tilde p=0$.
\nipar We now fix a standard mollification kernel in time $\xg\in\Cinf_c(\!(-1,0)\!)$ and, with $\rg_\eg,\rg$ as defined above, we define
\begin{align*}
\xg_\eg(t)\coloneq{}&\eg^{-1}\xg(t\eg^{-1}) \\
v(x,t)\coloneq{}&\xints{t}{t+\eg}\2(\tilde v\ast\rg_\eg)(x,s)\xg_\eg(t-s)\diff s, \\
p(x,t)\coloneq{}&\xints{t}{t+\eg}\2(\tilde p\ast\rg_\eg)(x,s)\xg_\eg(t-s)\diff s, \\
R(x,t)\coloneq{}&\lbar{\tilde v\otimes\tilde v}-v\otimes v,
\end{align*}
where
\[\lbar f=\int_t^{t+\eg}(f\ast\rg_\eg)(x,s)\xg_\eg(t-s)\diff s.\]
By construction and since $(\tilde v,\tilde p)$ solves \eqref{eq:FNS}, $(v,p,R)$ is a smooth solution of \eqref{eq:FNSR}, i.e.
\[\begin{sistema}
\pd_tv+\opn{div}(v\otimes v)+\grad p+\xfrl v=-\opn{div}R \\
\opn{div}v=0
\end{sistema}.\]
By using Jensen's inequality on $\T^3\x[t,t+\eg]$ with the measure $\rg_\eg(x-y)\xg_\eg(t-s)\diff x\diff s$, we conclude that
\[R=\lbar{\tilde v\otimes\tilde v}-v\otimes v\geq0. \xtag{ExStrictSubsolTensEst}\label{eq:ExStrictSubsolTensEst}\]
Coming to the initial datum, we have that
\[v|_{t=0}=\xints0\eg\!(\!(\tilde v-w_0)\ast\rg_\eg)(x,s)\xg_\eg(-s)\diff s+w_0\ast\rg_\eg.\]
Taking the $L^2$ norm, we can easily obtain that
\[\norm{v|_{t=0}-w_0}_{L^2(\T^3)}\leq\sup_{\mcl{t\in*[0,\eg]}}\norm{\tilde v(\per,t)-w_0}_{L^2(\T^3)}+\norm{w_0\ast\rg_\eg-w_0}_{L^2(\T^3)}\eqcolon\sup_{\mcl{t\in*[0,\eg]}}I_t+II_\eg.\]
$II_\eg$ can be made as small as we desire by choosing $\eg$ small enough. Let $\eg_0$ be the maximal parameter such that $II_\eg<\xfr\dg3$. As for $\sup I_t$, using $\tilde v_t(x)\coloneq\tilde v(x,t)$, we can obtain that
\begin{align*}
I_t^2={}&\xints{\T^3}{}\2\abs{\tilde v_t-w_0}^2\diff x=\xints{\T^3}{}\2\abs{\tilde v_t}^2-\abs{w_0}^2\diff x-2\xints0t\xints{\T^3}{}\2\ang{\pd_t\tilde v,w_0}\diff x\diff s \\
\mathbin{\overset\ast\leq}&2\xints0t\xints{\T^3}{}\2-(\tilde v_s\otimes\tilde v_s)\mathbin:\fD w_0+\ang{\tilde v_s,\xfrl w_0}\diff x\diff s \\
\mathbin{\overset\bullet\leq}&2\rad{C(1-2\EXP)}\xints0t\norm{\tilde v_t}_{L^2}[w_0]_1\diff s+2\xints0t\norm{\tilde v_t}^2_{L^2}\norm{\fD w_0}_{L^\8}\diff s \\
{}\leq{}&2t\norm{Dw_0}_{L^\8}\norm{w_0}_{L^2}(\rad{C(1-2\EXP)}+2\norm{w_0}_{L^2})\leq K(w)t\norm{\fD w_0}_{L^\8}.
\end{align*}
In $\ast$, we used the fact that $\|\tilde v_t\|_{L^2}^2\leq\|w_0\|_{L^2}^2$, i.e. \eqref{eq:energy}, as well as the fact that $(\tilde v,\tilde p)$ is a solution of \eqref{eq:FNS} and the fact $\opn{div}w_0=0$. In $\bullet$, we used \kcref{thm:teor:FrLHoeldNorms}, choosing $\eg=1-2\EXP$. In the last step, we used that $\|w_0\|_{L^2}\leq\|w\|_{L^2}$. This becomes arbitrarily small if we choose $t$ appropriately small, and since $t\in[0,\eg]$ that translates to $\eg$ small enough. Since $Dw_0=D\rg_{\hg_0}\ast w=\hg_0^{-4}D\rg(\hg_0^{-1}\per)\ast w$, Hölder's inequality yields
\[\norm{Dw_0}_{L^\8}\leq\hg_0^{-4}\norm{D\rg(\hg_0^{-1}x)}_{L^2_x(B_{\hg_0})}\norm{w}_{L^2}\leq\hg_0^{-4}\hg_0^2\|D\rg\|_{\s{C}^0}\norm{w}_{L^2}=C(w)\hg_0^{-2},\]
so that, to ensure $\sup_{[0,\eg]}I_t\leq{}^\dg/_2$, we choose
\[\eg\leq\xfr{\dg\hg_0^2}{3C(w)K(w)}\eqcolon\tilde\eg.\]
Choosing $\eg\coloneq\min\{\eg_0,\tilde\eg,\hg_0\}$ thus yields
\[\norm{\lbar v|_{t=0}-w}_{L^2}\leq\dg\qquad\norm{\lbar v|_{t=0}-w_0}_{L^2}\leq\xfr23\dg.\]
We have thus obtained \eqref{eq:ExStrictSubsolIV}. As for \eqref{eq:ExStrictSubsolEI}, we first notice that
\[\xints{\T^3}{}\2(\rg_\eg\ast f)(x,t)\diff x=\xints{\T^3}{}\2f(x,t)\diff t, \xtag{ESSConv}\label{eq:ESSConv}\]
for any $t\in[0,T]$ and any function $f$. Thus, by the definition of $R$, we have that
\[\xints{\T^3}{}\2\abs v^2(x,t)+\opn{tr}R(x,t)\diff x=\xints{\T^3}{}\2\lbar{\abs{\tilde v}^2}(x,t)\diff x=\xints{\T^3}{}\2(\xg_\eg\ast\abs{\tilde v}^2)(x,t)\diff x,\]
We have thus reduced \eqref{eq:ExStrictSubsolEI} to the following inequality:
\[\xfr12\xints{\T^3}{}\2(\xg_\eg\ast\abs{\tilde v}^2)(x,t)+\xints0t\xints{\T^3}{}\2\abs{\xfrl*v}^2(x,s)\diff x\diff s\leq\xints{\T^3}{}\2\abs w^2(x)\diff x+\dg.\]
Since $(\tilde v,\tilde p)$ satisfies \eqref{eq:energy}, we can see that
\begin{multline*}
\xfr12\xints{\T^3}{}\2(\xg_\eg\ast\abs{\tilde v}^2)(x,t)+\xints0t\xints{\T^3}{}\2\abs{\xfrl*v}^2(x,s)\diff x\diff s\leq\xfr12\xints{\T^3}{}\2\abs w^2\diff x+\xfr12\xints{\T^3}{}\2\abs{w_0}^2-\abs{w}^2\diff x \\
{}+\xfr12\xints{\T^3}{}\relax[(\xg_\eg\ast\abs{\tilde v}^2)-\abs{\tilde v}^2](x,t)\diff x+\xints0t\xints{\T^3}{}\2\sq{\abs{\xfrl*v}^2-\abs{\xfrl*\tilde v}^2}(x,s)\diff x\diff s \\
{\hsp{-2cm}}\eqcolon\xfr12\xints{\T^3}{}\2\abs w^2\diff x+I+II+III.
\end{multline*}
Our desired estimate \eqref{eq:ExStrictSubsolEI} will then follow from
\[I\leq\xfr\dg3\qquad II\leq\xfr\dg3\qquad III\leq\xfr\dg3. \xtag{ESSEITerms}\label{eq:ESSEITerms}\]
The first of these relations follows from \eqref{eq:ESSDefHg0}.
\nipar The second relation in \eqref{eq:ESSEITerms} follows, for $\eg$ sufficiently small, from the fact that, since $\tilde v\in L^2_tH^\EXP_x\sbse L^2_tL^2_x$, $|\tilde v|^2\in L^1_tL^1_x$, so $\xg_\eg\ast|\tilde v|^2\to|\tilde v|^2$ in $L^1_tL^1_x$ for $\eg\to0$.
\nipar Coming to the third relation, we first rewrite the integral of the first integrand:
\begin{align*}
\xints0t\xints{\T^3}{}\2\abs{\xfrl*\lbar v}^2(x,s)\diff x\diff s={}&\xints0t\xints{\T^3}{}\abs{\xints{s}{s+\eg}\xints{\T^3}{}\2\rg_\eg(x-y)\xfrl*\tilde v(y,\tg)\xg_\eg(t-\tg)\diff y\diff\tg}^2\diff x\diff s \\
{}\leq{}&\xints0t\xints{s}{s+\eg}\xints{\T^3}{}\2\abs{\xfrl*\tilde v}^2(x,\tg)\xg_\eg(t-\tg)\diff x\diff\tg\diff s,
\end{align*}
where we used Jensen's inequality in the first step and \eqref{eq:ESSConv} in the second one. Therefore, the remaining term is rewritten as:
\[III=\xints0t\xints{\T^3}{}\sq{\pa{\xg_\eg\ast\abs{\xfrl*\tilde v}^2}(x,s)-\abs{\xfrl*\tilde v}^2(x,s)}\diff x\diff s.\]
We now note that, since $\tilde v\in L^2_tH^\EXP_x$, we have that $\xfrl*\tilde v\in L^2_tL^2_x$, and thus $|\xfrl*\tilde v|^2\in L^1_tL^1_x$. Therefore, as before, $\xg_\eg\ast|\xfrl*\tilde v|^2-|\xfrl*\tilde v|^2\to0$ in $L^1_tL^1_x$, and the third relation of \eqref{eq:ESSEITerms} reduces to an opportune choice of $\eg$.
\nipar Summing up, $(v,p,R)$ is a smooth solution of \eqref{eq:FNSR}, which satisfies \eqref{eq:ExStrictSubsolEI} and \eqref{eq:ExStrictSubsolIV}, and $R\geq0$ by \eqref{eq:ExStrictSubsolTensEst}. The proof of \kcref{thm:lemma:ExistStrictSubsol} is thus complete.
\end{qeddim*}

\sect*{Acknowledgements}
The author wishes to thank his PhD advisors, Prof. Sara Daneri and Prof. Camillo De Lellis, for guiding this work with many useful conversations.

\end{document}